\documentclass[a4paper,11pt,final]{amsart}
\usepackage{amssymb,amsthm,amsmath,amscd}
\usepackage[all]{xy}
\newcommand{\nc}{\newcommand}
\renewcommand{\frak}{\mathfrak}
\providecommand{\cal}{\mathcal}
\renewcommand{\bold}{\mathbf}

\numberwithin{equation}{subsection}

\newcommand{\pfname}{Proof.}

\newenvironment{pf}{\vskip-\lastskip\vskip\medskipamount{\it\pfname}}%
                      {$\square$\vskip\medskipamount\par}
\newenvironment{pfof}[1]{\vskip-\lastskip\vskip\medskipamount{\it
    Proof of #1.}}%
                      {$\square$\vskip\medskipamount\par}

\newtheorem{thm}[subsubsection]{Theorem}

\newtheorem{cor}[subsubsection]{Corollary}
\newtheorem{corollary}[subsubsection]{Corollary}

\newtheorem{prop}[subsubsection]{Proposition}
\newtheorem{proposition}[subsubsection]{Proposition}

\newtheorem{lemma}[subsubsection]{Lemma} 

\theoremstyle{definition}
\newtheorem{definition}[subsubsection]{Definition} 

\theoremstyle{remark}

\newtheorem{remark}[subsubsection]{Remark}
\newtheorem{remarks}[subsubsection]{Remarks}

\newtheorem{example}[subsubsection]{Example} 
\newtheorem{examples}[subsubsection]{Examples} 

\nc{\Theorem}[1]{Theorem~\ref{#1}}
\nc{\Th}[1]{({\sl Th.}~\ref{#1})}
\nc{\Theorems}[2]{Theorems~\ref{#1} and ~\ref{#2}}
\nc{\Prop}[1]{({\sl Prop.}~\ref{#1})}
\nc{\Proposition}[1]{Proposition~\ref{#1}}
\nc{\Propositions}[2]{Propositions~\ref{#1} and ~\ref{#2}}
\nc{\Props}[2]{({\sl Props.}~\ref{#1} and ~\ref{#1})}
\nc{\Cor}[1]{({\sl Cor.}~\ref{#1})}
\nc{\Corollary}[1]{Corollary~\ref{#1}}
\nc{\Definition}[1]{Definition~\ref{#1}}
\nc{\Defn}[1]{({\sl Def.}~\ref{#1})}
\nc{\Lemma}[1]{Lemma~\ref{#1}} 
\nc{\Lem}[1]{({\sl Lem.} ~\ref{#1})} 
\nc{\Eq}[1]{equation~(\ref{#1})}
\nc{\Equation}[1]{Equation~(\ref{#1})}
\nc{\Section}[1]{Section~\ref{#1}}
\nc{\Sec}[1]{({\sl Sec.} ~\ref{#1})} 
\nc{\Chapter}[1]{Chapter~\ref{#1}}
\nc{\Chapt}[1]{({\sl Ch.}~\ref{#1})}

\nc{\Ex}[1]{{\sl Ex.}~\ref{#1}}
\nc{\Exa}[1]{{\sl Example}~\ref{#1}}
\nc{\Example}[1]{{\sl Example}~\ref{#1}}
\nc{\Examples}[1]{{\sl Examples}~\ref{#1}}
\nc{\Exercise}[1]{{\sl Exercise}~\ref{#1}}

\nc{\Rem}[1]{({\sl Rem.}~\ref{#1})}
\nc{\Remark}[1]{{\sl Remark}~\ref{#1}}
\nc{\Note}[1]{{\sl Note}~\ref{#1}}

\nc{\Conjecture}[1]{Conjecture~\ref{#1}}
\nc{\Claim}[1]{Claim~\ref{#1}}

\nc \Proof{{  \it Proof. }}

\theoremstyle{definition}

\theoremstyle{remark}

\nc{\Sats}[1]{Sats~\ref{#1}}
\nc{\Sa}[1]{({\sl Sa.}~\ref{#1})}
\nc{\Kor}[1]{({\sl Kor.}~\ref{#1})}
\nc{\Korollarium}[1]{Korollarium~\ref{#1}}
\nc{\Exe}[1]{{\sl Exempel}~\ref{#1}}
\nc{\Anm}[1]{{\sl Anmärkning}~\ref{#1}}
\nc{\Not}[1]{{\sl Not}~\ref{#1}}

\nc{\Formodan}[1]{Förmodan~\ref{#1}}
\nc{\Pastaende}[1]{Påstående~\ref{#1}}

                      {$\square$\vskip\medskipamount\par}

\nc \Ab{{\bold A}}
\nc \Gb{{\bold G}}
\nc \Qb{{\bold Q}}
\nc \Rb{{\bold R}} \nc \Cb{{\bold C}} 
\nc \Fb{{\bold F}}
\nc \Pb{{\bold P}}
\nc \Zb{{\bold Z}} 
\nc \Nb{{\bold N}} 

\nc \mf{\frak m} \nc \mh{\hat{\m}} 
\nc \nf{\frak n}
\nc \Of{\frak O}
\nc \rf{\frak r}
\nc \hf{\frak h} 
\nc \qf{\frak q} 
\nc \bfr{\frak b} 
\nc \kfr{\frak k} 
\nc \pfr{\frak p} 
\nc \af{\frak a }
\nc \cf{\frak c }
\nc \sfr{\frak s} 
\nc \ufr{\frak u} 
\nc \g{\frak g} 
\nc \gA{\g_{\Ao}} 
\nc \lfr{\frak l}
\nc \afr{\frak a}
\nc \gfh{\hat {\frak g}}
\nc \gl{\frak { gl }}
\nc \Sl{\frak {sl}}

\newcommand{\on}{\operatorname}

\nc \Ann {\on {Ann}}
\nc \Coker {\on {Coker}}
\nc \Co{\on C}
\nc \Ker {\on {Ker}}
\nc \No {\on N}
\nc \Oo {\on O}
\nc \ch {\on {ch}}
\nc \rk {\on {rk}}
\nc \Top {\on {T}}
\nc \Tr {\on {Tr}}
\nc \tr {\on {tr}}
\nc \trdeg {\on {tr.deg}}
\nc \codim {\on {codim }}
\nc \coht {\on {coht}}
\nc \coh {\on {coh}}
\nc \qcoh {\on {qcoh}}
\nc \grade {\on {grade}}
\nc \hto {\on {ht}}
\nc \depth {\on {depth}}
\nc \prof {\on {prof}}
\nc \ind{\on {ind}}
\nc \prodo{\on {prod}}
\nc \coind{\on {coind}}
\nc \Con{\on {Con}}
\nc \Crit{\on {Crit}}
\nc \Der{\on {Der}}
\nc \Char{\on {Char}}
\nc \Ch{\on {Ch}}

\nc \Ext{\on {Ext}}
\nc \Eo{\on {E}}
\nc \End{\on {End}}
\nc \ad{\on {ad}}
\nc \Ad{\on {Ad}}
\nc \gr{\on {gr}}
\nc \Fo{\on {F}}
\nc \Gr{\on {Gr}}
\nc \Go{\on {G}}
\nc \Glo{\on {Gl}}
\nc \Ho{\on {H}}
\nc \hol{\on {hol}}
\nc \supp{\on {supp}}
\nc \ssupp{\on {s-supp}}
\nc \singsupp{\on {singsupp}}
\nc \msupp{\on {msupp}}
\nc \spec{\on {spec}}
\nc \Max{\on {Max}}
\nc \Mod{\on {Mod}}
\nc \Rad {\on {Rad}}
\nc \rad {\on {rad}}
\nc \rank {\on {rank}}
\nc \Slo{\on {SL}}
\nc \soc {\on {soc}}
\nc \Irr {\on {Irr}}
\nc \Imo {\on {Im}}
\nc \SSo{\on {SS}}
\nc \lub{\on {lub}}
\nc \gldim{\on {gl.d.}}
\nc \pdo{\on {p.d.}} 
\nc \ido{\on {i.d.}} 
\nc \dSSo{\dot {\SSo}}
\nc \So{\on S}
\nc \I{\on I}
\nc \Io{\on I}
\nc \Jo{\on J}
\nc \jo{\on j}
\nc \Ko{\on K}
\nc \PBW{\Ac_{PBW}}
\nc \Ro{\on R}
\nc \To{\on T}
\nc \Ao{\on A}

\nc \Do{{\on D}}
\nc \Bo{\on B}
\nc \Po{\on P}
\nc \Qo{\on Q}
\nc \Zo{\on Z}
\nc \U{\on U}
\nc \wt{\on {wt}}
\nc \Uh{\hat {\U}}
\nc \T{\on T}
\nc \Lo{\on L}
\nc{\dop}{\on d}
\nc{\eo}{\on e}
\nc{\ado}{\on{ad}}
\nc{\Tot}{\on{Tot}}
\nc{\Aut}{\on{Aut}}

\nc{\overrightleftarrows}[2]{\overset{#1}{\underset{#2}{\rightleftarrows}}}

\nc{\CCF}{\cal{CF}}
\nc{\CDF}{\cal{DF}}
\nc{\CHC}{\check{\cal C}}
\nc{\CHom}{\cal{H}om}
\nc{\Cone}{\on{Cone}}
\nc{\dec}{\on{dec}}
\nc{\Diff}{\on{Diff}}
\nc{\dirlim}{\underset{\rightarrow}{\on{lim}}}
\nc{\dpar}{\partial}
\nc{\GL}{\on{GL}}
\nc{\CGr}{\cal{G}r}
\nc{\pr}{\on{pr}}
\nc{\semid}{|\!\!\!\times}
\nc{\Hom}{\on{Hom}}
\nc \RHom{\on {RHom}}
\nc{\Id}{\on{Id}}
\nc{\id}{\on{id}}
\nc{\Ima}{\on{Im}}
\nc{\invtimes}{\underset{\leftarrow}{\otimes}}
\nc{\invlim}{\underset{\leftarrow}{\on{lim}}}
\nc{\Lie}{\on{Lie}}
\nc{\res}{\on{res}}
\nc{\re}{\on{Re }}
\nc{\Pic}{\on{Pic }}
\nc{\LPic}{\on{LPic }}
\nc{\Sch}{\on{Sch}}
\nc{\Set}{\on{Set}}
\nc{\Spec}{\on{Spec}}
\nc{\mSpec}{\on{mSpec}}
\nc{\Specb}{\bold {Spec}}
\nc{\Projb}{\bold {Proj}}
\nc{\Specan}{\on{Specan}}
\nc{\Spo}{\on{Sp}}
\nc{\Spf}{\on{Spf}}
\nc{\sym}{\on{sym}}
\nc{\symm}{\on{symm}}
\nc{\rop}{\on{r}}
\nc{\Td}{\on{Td}}
\nc{\Tor}{\on{Tor}}

\nc{\Artin}{\cal{A}rtin}
\nc{\Dgcoalg}{\cal{D}gcoalg}
\nc{\Dglie}{\cal{D}glie}
\nc{\Ens}{\cal{E}ns}
\nc{\Fsch}{\cal{F}sch}
\nc{\Groupoids}{\cal{G}roupoids}
\nc{\Holie}{\cal{H}olie}
\nc{\Mor}{\cal{M}or}

\nc \Kc{{\Cal K}}
\nc \Lc{{\cal L} }
\nc \lc{{\cal l}} 
\nc \CC{{\cal C}} 
\nc \Cc{{\cal C} }
\nc \Pc{{\cal P} }
\nc \Dc{{\cal D}}
\nc \Ac{{\cal A}} 
\nc \Bc{{\cal B}}
\nc \Ec{{\cal E}}
\nc \Mc{{\cal M}} \nc \hM{\hat{\Mc}} \nc \bM{\bar {\Mc}} \nc
\hbM{\hat{\bar \Mc}}  
\nc \Nc{{\cal N}}
\nc \Hc{{\cal H}} 
\nc \Ic{{\cal I}} 
\nc \Oc{{\cal O} }
\nc \Och{\hat{\cal O}} 
\nc \Sc{{\cal S} }
\nc \Tc{{\cal T}} 
\nc \Vc{{\cal V}} 
\nc{\CA}{{\cal A}}
\nc{\CB}{{\cal B}}
\nc{\C}{{\cal F}}
\nc{\Gc}{{\cal G}}
\nc{\CH}{{\cal H}}
\nc{\CI}{{\cal I}}
\nc{\CM}{{\cal M}}
\nc{\CN}{{\cal N}}
\nc{\CO}{{\cal O}}
\nc{\Rc}{{\cal R}}
\nc{\CT}{{\cal T}}
\nc{\CU}{{\cal U}}
\nc{\CV}{{\cal V}}
\nc{\CZ}{{\cal Z}}

\nc{\fa}{\frak{a}}
\nc{\fA}{\frak{A}}
\nc{\fg}{\frak{g}}
\nc{\fh}{\frak{h}}
\nc{\fI}{\frak{I}}
\nc{\fK}{\frak{K}}
\nc{\fm}{\frak{m}}
\nc{\fP}{\frak{P}}
\nc{\fS}{\frak{S}}
\nc{\ft}{\frak{t}}
\nc{\fX}{\frak{X}}
\nc{\fY}{\frak{Y}}

\nc{\bF}{\bar{F}}
\nc{\bCP}{\bar{\cal{P}}}
\nc{\bm}{\mbox{\bf{m}}}
\nc{\bT}{\mbox{\bf{T}}}
\nc{\hB}{\hat{B}}
\nc{\hC}{\hat{C}}
\nc{\hP}{\hat{P}}
\nc{\htest}{\hat P}

\nc{\nen}{\newenvironment}
\nc{\ol}{\overline}
\nc{\ul}{\underline}
\nc{\ra}{\rightarrow}
\nc{\lla}{\longleftarrow}
\nc{\lra}{\longrightarrow}
\nc{\Lra}{\Longrightarrow}
\nc{\Lla}{\Longleftarrow}
\nc{\Llra}{\Longleftrightarrow}
\nc{\hra}{\hookrightarrow}
\nc{\iso}{\overset{\sim}{\lra}}

\nc{\dsize}{\displaystyle}
\nc{\sst}{\scriptstyle}
\nc{\tsize}{\textstyle}

\nen{exa}[1]{\label{#1}{\bf Example.\ } }{}

\nen{rem}[1]{\label{#1}{\em Remark.\ } }{}
\nen{note}[1]{\label{#1}{\em Note.\ } }{}
\nen{exer}[1]{\label{#1}{\em Exercise.\ } }{}

\begin{document}

\pagestyle{plain}


\title{Smooth modules over Lie algebroids I}

\author{Rolf K{\"al}lstr{\"o}m}

\address{Department of Mathematics, Royal Institute of Technology, 100
  44 Stockholm, Sweden} \email{rok@math.kth.se}  \date{\today} \maketitle

\tableofcontents
\section{Introduction} 
Let $(\Oc_X , \g_X , \alpha)$ be a Lie algebroid \Defn{Lie algebroid}
on a space $X$ where $\alpha : \g_X \to T_X$ is a homomorphism to the
tangent sheaf; a space is either a scheme locally of finite type over
an algebraically closed field $k$ of characteristic $0$ or a complex
analytic space ($k=\Cb$). Let $M$ be a coherent module over $\g_X$
\Defn{gmodule}, $\g'_X$ a sub-algebroid of $\g_X$, and $x$ a point in
$X$. We say that $M$ is smooth along $\g'_X$ at $x$ if its stalk $M_x$
contains an $\Oc_x$-module $M_x^0$ of finite type which generates
$M_x$ over $\g_x$ and satisfies $\g'_x \cdot M^0_x \subseteq M_x^0$.
In this paper and its continuation \cite{kallstrom:2smooth} we shall
study this notion of smoothness, treating the algebraic and complex
analytic theory simultaneously.  Our main results will be described a
little later.

Let $S\subset X$ be a subspace and $T_X(I_S)\subset T_X$ be the
sub-algebroid of sections that preserve the ideal $I_S$ of $S$.
Particularly interesting torsion free (over $\Oc_X$) $T_X$-modules are
those that are smooth along $T_X(I)$. When $S$ is a divisor and $X$ is
normal then $M$ has this property if and only if it has regular
singularities along $S$ in the sense of P. Deligne \cite{deligne:eq};
there are similar results for torsion modules.
 
Deligne's modestly titled but very influential work [loc. \ cit.]
contains the first treatment of regular singularities of connections
on complex spaces. Together with M. Kashiwara's work
\cite{kashiwara:master} it forms a starting point for the theory of
$\Dc_X$-modules (= $T_X$-modules), leading to a more uniform treatment
of singularities.  Some important architects of the analytic theory are
Kashiwara, M.  Sato, T.  Kawai, and Z. Mebkhout, while the algebraic
theory is mainly due to J.  Bernstein; see the text-books
\cite{bjork:analD}, \cite{borel:Dmod}, and \cite{bei-ber:jantzen}.
Lie algebroids also play a prominent role in differential geometry,
where one of the first uses is due to M.  Atiyah
\cite{atiyah:connections}; see the survey \cite{mackenzie-kirill}.
 
 The intrinsic interest in Lie algebroids should be rather obvious,
 but since the paper does not deal much with specific applications we
 briefly mention some connections with other areas.
 
 {\it Quantum gauge theory}\/: The dictionary from the terminology of
 physics to mathematics is that gauge fields, field strengths, and the
 Yang-Mills equations (the field equations) correspond to sections,
 sections of the curvature, and the Bianchi identity for the Hodge
 dual of the curvature of Lie algebroids; and quantum gauge fields act
 on modules over Lie algebroids. Gauge theories first attracted
 differential geometers and analysts, considering the geometry of
 principal bundles and solutions of the Yang-Mills equations. Now an
 interesting case of gauge theory on curves occurs in conformal field
 theory, studying Lie algebroids on curves and Lie algebroids on
 moduli spaces of curves. The subject has therefore attracted the
 attention of a great many algebraic geometers; see
 \cite{beilinson-schechtman}, \cite{Tsuchiya-ueya}, and
 \cite{ueno:conformal}.

 {\it Cohomology}\/: Let $M$ be a regular holonomic $T_X$-module on a
 non-singular space $X$.  Then the local cohomology sheaves
 $H^i_{[S]}(M)$ also are regular holonomic $T_X$-modules, even when
 they are not coherent over $\Oc_X$; see e.g.  \cite{lyubeznik} for
 one use of the $T_X$-module structure on $H^i_{[S]}(M)$.  As for the
 cohomology of $S$ itself a nice version is intersection cohomology,
 occurring as the homology of certain objects in the category of
 perverse sheaves, which by the Riemann-Hilbert correspondence is
 equivalent to the category of regular holonomic $T_X$-modules whose
 support belong to $S$; see \cite{kashiwara-ka},
 \cite{mebkhout:equivalence}, \cite{borel:Dmod} and \cite{be-be-de}.
 Similar facts hold for modules over any transitive Lie algebroid.

 {\it Localisation of representations}\/: Often one is concerned with
 homo\-mor\-phisms of Lie algebras over $k$ of the type $\alpha :\g
 \to \Gamma(X, T_X)$: The action of a Lie group on a variety $G\times
 X \to X$ is one source of such $\alpha$; other interesting sources
 occur in moduli problems related to curves, where $\alpha$ is the
 inverse of a Kodaira-Spencer mapping; see [loc.\ cit.].  The
 homomorphism $\alpha$ gives the $\Oc_X$-module
 $\g^0_X:=\Oc_X\otimes_k \g$ a structure of Lie algebroid, and letting
 $\g_X$ be the push-out of $\g^0_X$ by a character of $\Ker (\g^0_X
 \to T_X)$ one has a homomorphism $\alpha: \Dc(\g) \to \Dc (\g_X)$
 from the enveloping algebra of $\g$ to the enveloping ring of
 differential operators of $\g_X$.  Now the $\g_X$-module
 $M_X:=\Dc(\g_X)\otimes_{\Dc(\g)}M$ serves as a bridge to a geometric
 description of a $\g$-module $M$ using the Riemann-Hilbert
 correspondence. For example, when $X$ is a flag variety the
 localisation functor $M \mapsto M_X$ is essential in the proof of the
 Kazhdan-Lusztig algorithm for the Jordan-H{\"o}lder multiplicities of
 Verma modules over a complex semi-simple Lie algebra; see
 \cite{beilinson-be}, \cite{brylinski-ka}, \cite{bei-ber:jantzen},
 \cite{kashiwara:rep}.

 {\it Families}\/: Coherent sheaves of Lie algebras $\g_X$ over
 $\Oc_X$ (the Lie bracket being $\Oc_X$-bilinear) are obvious examples
 of Lie algebroids. We think of these as families of Lie algebras
 $\g(x) :=\g_x/\mf_x\g_x$ over the fields $k(x)= \Oc_x/\mf_x$, where a
 trivial family is one which is isomorphic to an extension of scalars
 $\g_X = \Oc_X\otimes_k \g$ of a Lie algebra $\g$ over a field $k$;
 more interesting cases occur when the Lie bracket $[\cdot , \cdot]_x$
 for $\g(x)$ degenerates at special points $x\in X$ (contraction).
 Similarly, a sheaf of modules $M$ over a sheaf of Lie algebras $\g_X$
 we regard as a family of modules, resulting in a fibre module
 $M(x):=M/\mf_x M$ over the $k(x)$-Lie algebra $\g (x)$ at each point
 $x$.  An interesting decreasing filtration of $M(x)$ by
 $\g(x)$-modules $\cdots \subset M_n(x) \subset M_{n-1}(x)\subset
 \cdots \subset M_0(x) =M(x)$ is defined using an ideal $I \subseteq
 \mf_x$ and a section $S\in Hom_{\g_X}(M\otimes_{\Oc_X}M, \Oc_X)$
 (where usually the $\g_X$-module $\Oc_X$ is trivial; $S$ is a
 ``Shapovalov form''), setting $M_n(x) := \{m\in M : S(m,M)\subseteq
 I^n \}\mod \mf_x$.  Jantzen's filtration of Verma modules is an
 important example; see \cite{jantzen:mod}, \cite{bei-ber:jantzen}.

The content of the sections is as follows:

\Section{def-section}: is mostly a survey containing general
definitions and basic facts about Lie algebroids, complementing and
recalling parts of \cite{bei-ber:jantzen}.  There are no really new
results and we include more material than is absolutely needed for our
main results, but it should be helpful to have this information
collected.

\Section{Modules}: \Proposition{bfacts} states that it is equivalent that a space $X$ is
non-singular, that the $T_X$-module $\Oc_X$ contains no coherent
proper submodules, and that the Lie algebroid $T_X$ contains no
coherent ideals; this should be well-known. The main part of this
section has to do with operations on Lie algebroids and their derived
categories of complexes of modules, generalizing standard operations
in $\Dc$-module theory.

\Section{sec-smooth}:  \Propositions{point}{globalsmooth} are  relations between
global and local smoothness. Limiting ourselves to morphisms $\pi : Y
\to X$ of non-singular spaces we first show that the derived pull-back
functor $\pi^+$ preserves the derived category $D^b_s$ of bounded
complexes of $\g$-modules with smooth homology modules
\Prop{preservation}, and we also give a partial converse
\Prop{inv-equiv}. \Theorem{direct-smooth} gives a condition when the
direct image functor preserves $D^b_s$; it will later be used to prove
that the direct image of a module with regular singularities again has
regular singularities. We record the important fact implied by the
Artin-Rees lemma, that smoothness of a module is well reflected in its
completion \Prop{complete}. \Theorem{gaga-eq} contains the expected
equivalences implied by GAGA; its \Corollary{gaghom} is a statement
about cohomology groups of certain sheaves which are not
$\Oc$-modules.

\Section{prolonging}:  
We say that an $\Oc_X$-module $M$ is pure if all its non-zero coherent
sub-modules have the same dimension of their support; $M$ is torsion
free if it is pure and $\dim \supp M = \dim X$. Let $M$ be a
coherent $\Oc_X$-module with support $V\subset X$, and $j: \Omega \to
X$ an open inclusion of spaces such that $\codim_V (V\cap (X\setminus
\Omega))\geq 2$. Then if $j^*(M)$ is pure we have $j_*j^*(M) \in \coh
(\Oc_X)$ \Cor{coherent}.  In the algebraic case 
\Corollary{coherent} follows from A. Grothendieck's finiteness theorem
\cite[Exp. VIII, p. 13, Th. 3.1]{grothendieck:cohomologie} while in the
complex analytic case it follows from a corresponding result by Y.-T.
Siu and G. Trautmann \cite{siu:local},\cite{trautmann}, as stated in
\Theorem{groth-coherent}. We include a proof of this theorem in the
algebraic case which perhaps is somewhat more direct than that of
Grothendieck.

Using \Corollary{coherent} we can prove our first main result in
\Theorem{relative}, stating: Let $\pi :Y \to X$ be a proper morphism
of spaces, $S\subset X$ a (closed) subspace, and $\g'_Y \subset \g_Y$
an inclusion of Lie algebroids on $Y$. Suppose that $M\in \coh(\g_Y)$
contains, locally in $X$, an $\Oc_Y$-lattice $M^0$, i.e. $M=
\Dc(\g_Y)M^0$ and $M^0\in \coh (\Oc_X)$ (such $M^0$ always exist when
$X$ is a noetherian variety), and let $ \cup_{i=1}^k V_i$ be an
irreducible decomposition of the support of $M$.  Assume that $V_i\cap
V_j \cap \pi^{-1}(S)= \emptyset$ for all pairs $i\neq j$ and that
$\dim \pi (V_i) \geq \dim S +2$ when $V_i\cap \pi^{-1}(S) \neq
\emptyset$.  Assume moreover that the maximal coherent
$\g_X$-submodule of $M$ with support in $V_i$ is pure, $i=1,\ldots,
k$, and that the canonical homomorphisms
\begin{displaymath}
  \pi^*\pi_*(M)_y\to M_y \quad \text{ and } \quad  \pi^*\pi_*(\g'_Y)_y\to \g'_y
\end{displaymath}
are surjective when $y\in \pi^{-1}(S)$. Then if $M$ is point-wise
smooth along $\g'_Y$ in $Y \setminus \pi^{-1}(S)$, it follows that it
is point-wise smooth in $Y$. 

Letting $Y=X$ and $\pi=\id$ we obtain a non-relative version: Let $M$
be a coherent $\g_X$-module, torsion free over $\Oc_X$, and $\Gamma$ a
closed subset such that $\codim_X \Gamma \geq 2$. Then if $M$ is
point-wise smooth along a sub-algebroid $\g'_X \subseteq \g_X$ at each
point in $X\setminus \Gamma$ it follows that $M$ is smooth along
$\g'_X$ in all of $X$ \Cor{pure}; \Corollary{purestcor} is a similar
result for pure modules. Moreover, if $M$ is point-wise smooth along
$\g'_X$ in the complement of a divisor $S\subset X$, then the points
in $S$ where $M$ is smooth form a closed set, so one gets as an added
dividend that it suffices to check smoothness in one point in every
irreducible component of the non-singular locus of $S$
\Cor{point-cor}.  This implies Deligne's theorem that a $T_X$-module,
coherent and torsion free over $\Oc_X(*S)$, has regular singularities
along $S$ (in the sense of pull-backs to curves) if this holds at all
points in a dense subset of the non-singular locus of $S$.

\Section{regularity}:  
\Theorem{fibre:1} is a useful smoothness criterion, and follows from a
straightforward application of Gabber's integrability theorem.
\Theorem{curve-test} is a curve test applicable to torsion free
modules for any pair of Lie algebroids $\g'_X \subset \g_X$.  The
remaining part of the section is concerned with $\g_X$-modules with
regular singularities, where $\g'_X$ is determined by a sub-space in
$X$. When $M$ is a pure coherent $\g_X$-module with support $V$ we say
that it has regular singularities along a divisor $S\subset V \subset
X$ if it is smooth along $\g_X(I)\cap \g_X(J)$; here $I$ and $J$ are
the ideals of $S$ and $V$, and $\g_X(I) = \alpha^{-1}(T_X(I))$. We
refer to \Section{reg-sing} for the general case when the module is
not pure.  From now assume that the spaces $X$ and $Y$ are {\it
  non-singular}, and the Lie algebroids $\g_X$ and $\g_Y$ are locally
free and transitive.  Let $D_{rs}^b(\g_X)$ be the derived category of
bounded complexes of $\g_X$-modules $M^\bullet$ whose homology
$H^{\bullet}(M^\bullet)$ is a $\g_X$-module with regular
singularities.

We have $D^b_{rs}(\g_X)\subset D^b_{\hol}(\g_X)$\Prop{holonomic} where
$D^b_{\hol}(\g_X)$ is the category of complexes with holonomic
homology, and if $S\subset X$ is a subspace, then
$R\Gamma_{[S]}(M^\bullet), M(*S)^\bullet\in D^b_{rs}(\g_X)$ if and
only if $M\in D^b_{rs}(\g_X)$ \Th{bernstein}; this corresponds to a
theorem by Bernstein.

Let $\pi : Y \to X$ be a morphism of spaces. The inverse
image functor preserves the property of having regular singularities,
$\pi^!(D^b_{rs}(\g_X))\subset D^b_{rs}(\g_Y)$ \Th{inv-rs}, and if
$\pi$ is proper, the direct image functor has the same good behaviour
\Th{directreg}
 \begin{displaymath}
   \pi_+(D^b_{rs}(\g_Y))\subset D^b_{rs}(\g_Y).
 \end{displaymath}
 
 Let $D^b_{crs}(\g_X)\subset D^b_{\coh}(\g_X)$ be the sub-category of
 complexes $M^\bullet$ such that the pull back $\pi^+(M^{\bullet})\in
 D^b_{rs}(\g_C)$ when $\pi : C\to X$ is an embedding of a curve (curve
 regular complexes).  We have \Th{equivalent}
\begin{displaymath}
  D^b_{rs}(\g_X) = D^b_{crs}(\g_X).
\end{displaymath}
As a consequence any sub-quotient of the homology of a complex in
$D^b_{crs}(\g_X)$ is curve regular \Cor{full-reg}. That
$D^b_{crs}(T_X)$ has this property is known but the previous proofs
are quite hard; see \cite{borel:Dmod} and \cite{bjork:analD}.

In the remaining part we study the sub-category $\bar D^b_{rs}(\g_X)
\subset D^b_{rs}(\g_X)$ of completely regular complexes on a
quasi-projective algebraic manifold $X$.

A locally free and transitive Lie algebroid $(\bar X , \g_{\bar X},
\bar \alpha)$ is a {\it completion} of a Lie algebroid $(X, \g_{X},
\alpha)$ when there exists an open embedding $j: X \to \bar X$, where
$\bar X$ is projective, such that $\g_X=j^*(\g_{\bar X})$ and $\alpha
= j^*(\bar \alpha): \g_X\to T_X$.  We say that $M^\bullet \in \bar
D^b_{rs}(\g_X)$ if there exists such a completion so that
$j_+(M^\bullet)\in \bar D^b_{rs}(\g_{\bar X})$.   This definition is
intrinsic, independent of the choice of completion \Lem{intrinsic},
and $\bar D^b_{rs}(T_X)$ coincides with the category of regular
holonomic $T_X$-modules as defined in \cite{borel:Dmod}
\Cor{complete-equiv}.  We have $\pi^!(\bar D^b_{rs}(\g_X))\subset \bar
D_{rs}^b(\g_Y)$ \Cor{inv-complete}. Conversely, if $\pi(Y)$ contains
the support of $M\in \bar D^b_{\coh}(\g_X)$ and $\pi^!(M)\in \bar
D^b_{rs}(\g_Y)$, then $M\in \bar D^b_{rs}(\g_X)$
\Prop{inv-equiv-complete}, a result that can be used to prove that
certain equivariant $\g_X$-modules are completely regular.  That
\begin{displaymath}
   \pi_+ (\bar D^b_{rs}(\g_Y)) \subset \bar D^b_{rs}(\g_X)
\end{displaymath}
follows from \Theorem{directreg} when $\pi$ is proper, and it holds
more generally when $\pi$ can be factorized $\pi= p\circ j$ where $j$
is a completion and $p$ is proper \Prop{completedirect}.  Let us
compare to the proof of \cite[VII, Th.  12.2]{borel:Dmod}, that
$\pi_+(\bar D^b_{crs}(T_Y))\subset \bar D^b_{crs}(T_Y)$. This proof is
inspired by Deligne's proof \cite{deligne:eq} that the Gauss-Manin
connection has regular singularities; one first has to prepare by
constructing a generating class of objects, the standard modules, and
one considers only completions $Y \to \bar Y$ such that $\bar Y
\setminus Y$ is a divisor with normal crossing singularities (regular
completions), which exist by Hironaka's theorem.  Our proof is also
based on Deligne's proof, but otherwise quite different. We do not
need that fact that standard modules form a generating class, and our
use of Hironaka's theorem is instead for taking care of modules whose
support is not non-singular in codimension $1$; one may also note that
it is unnecessary to only consider regular completions.

I want to thank Rikard B{\"o}gvad for useful suggestions to improve
the disposition of this work.

\section{Lie algebroids} 
\label{def-section}
\subsection{Definitions} \label{def-section-def}
By a {\it variety}\/ we shall mean a separated scheme locally of finite type
over an algebraically closed field $k$ of characteristic $0$.  A {\it
  complex analytic space} \/ is a separated complex analytic space
($k=\Cb$). By {\it space}\/ we mean either a variety or a complex
analytic space.

\begin{remark} Our terminology is a little unconventional since a
  variety usually is an integral separated scheme of finite
  type over $k$.  Note also that by Lefshetz' principle one can
  assume that all varieties of finite type are defined over
  $k=\Cb$, but not invoking transcendental methods to prove
  algebraic results we prefer also in such cases to work over a
  specific algebraically closed field $k$ of characteristic
  $0$.
\end{remark}

Let $\coh (\Oc_X) \subset \Mod (\Oc_X)$ be the subcategory of coherent
sheaf of modules over the structure sheaf $\Oc_X$ on a space $X$; we
abbreviate sheaf of rings and sheaf of modules by ring and module. The
ring $\Oc_X$ and its $\Oc_X$-module of $k$-linear derivations $T_X$
both belong to $\coh (\Oc_X)$. Let $M$ be an $\Oc_X$-module. Then
$M_x$ is the stalk at a point $x$ in $X$, and if $\Omega$ is an open
subset of $X$ the vector space of sections over $\Omega$ is $\Gamma
(\Omega , M)$ or $M(\Omega)$; the restriction of $M$ to $\Omega$ is
$M_\Omega$.  When we say that $m$ is a vector in $M$ we mean that $m$
is a section that is defined in some open subset of $X$.  Recall that
a space $(\Oc_X, X)$ is integral when the rings $\Oc_X(\Omega)$ are
integral domains, and this is the same as $(\Oc_X, X)$ being
irreducible and reduced.  For a locally principal (Weil) divisor
$S\subset X$, $\Oc_X(*S)$ is the coherent ring defined locally by
$\Oc_X[1/f]$ when $f$ is a local generator of $S$, and we put $M(*S)=
\Oc_X(*S)\otimes_{\Oc_X}M$.  Let $S\subset X$ be a closed subspace
(algebraic or complex analytic) and $j: X\setminus S \to X$ the
canonical open immersion. We let $\Gamma_S(M)\subseteq M$ be the
maximal submodule whose support belongs to $S$ and
$\Gamma_{[S]}(M)=H^0_{[S]}(M)= \{m\in M: I^n_S\cdot m =0\}$, where
$I_S$ is the ideal of $S$; $R\Gamma_S(\cdot)$ and
$R\Gamma_{[S]}(\cdot)$ are the derived versions of $\Gamma_S(\cdot)$
and $\Gamma_{[S]}(\cdot)$, defining functors on the derived category
of complexes of $\Oc_X$-modules $D(\Oc_X) := D(\Mod (\Oc_X))$. The
local cohomology modules of $M^\bullet \in D(\Oc_X)$ are
$H^i_S(M^\bullet) = R^i\Gamma_S(M^\bullet)$ and $H^i_{[S]}(M^\bullet)=
R^i\Gamma_{[S]}(M^\bullet)$, and there exist distinguished triangles
triangles\footnote{See \cite{be-be-de} for an account of triangulated
  categories.} in $D(\Oc_X)$ 
\begin{eqnarray}\label{dist-triangle}
  R\Gamma_S(M^\bullet) \to M^\bullet \to Rj_*j^*(M^\bullet)\to,\\\nonumber 
R\Gamma_{[S]}(M^\bullet)\to M^\bullet \to M^\bullet(*S)  \to,
\end{eqnarray}
which coincide when $X$ is a scheme.

The work of Beilinson and Bernstein \cite{bei-ber:jantzen} contains a
succinct treatment of Lie algebroids and serves as an important source
for this section; this section can be regarded as a complement to some
of the most basic material in {\it loc.  \ cit.}.

\begin{definition}\label{Lie algebroid}
  A Lie algebroid is a triple $(\Oc_X, \g_X , \alpha)$, where:
  \begin{enumerate}
  \item $\g_X$ is a Lie algebra over $k$ on a space $(\Oc_X,X)$;
  \item $\g_X$ is a coherent  $\Oc_X$-module;
  \item $\alpha$ is a homomorphism of sheaves
  \begin{displaymath}
    \alpha : \g_X \to T_X
  \end{displaymath}
  such that $\alpha([\delta_1 , \delta_2]) = [\alpha(\delta_1), \alpha
  (\delta_2)]$ and $\alpha (f\delta)= f\alpha(\delta)$ for $f\in
  \Oc_X$, $\delta_1 , \delta_2 \in \g_X$. ($\alpha$ is a homomorphism
  for the structures in 1 and 2).  We require the compatibility
  $[\delta_1 ,f \delta_2] = \alpha (\delta_1)(f) \delta_2 + f
  [\delta_1 , \delta_2] $, $f \in \Oc_X$, $\delta_1 , \delta_2 \in
  \g_X$.
  \end{enumerate}
 Let $\Lie_X$ be the category of Lie
  algebroids on $X$.  We say that $(\Oc_X , \g_X , \alpha) \in \Lie_X$
  is {\it transitive} \/ if $\alpha$ is surjective; if $\alpha =0$
  then $\g_X$ is a sheaf of Lie algebras over $\Oc_X$.
\end{definition}
A homomorphism of Lie algebroids $\phi :(\g_X ,\alpha) \to (\af_X ,
\beta )$ (on the same space $X$) is a homomorphism of $\Oc_X$-modules
and $k$-Lie algebras $\phi :\g_X \to \af_X$ such that $\alpha = \beta
\circ \phi$. Of course, we often abbreviate a Lie algebroid $(\Oc_X,
\g_X , \alpha)$ to $\g_X$. It should be obvious what is meant by an
ideal of a Lie algebroid; for instance $\Ker (\phi)$ is an ideal if
$\phi$ is a homomorphism of Lie algebroids and the sub-sheaf $\bfr =
\Ker (\alpha)$ is an ideal which is moreover a coherent $\Oc_X$-Lie
algebra.

Let $\af_X$ be a subalgebroid of a Lie algebroid $\g_X$. The {\it
  normalizer}\/ of $\af_X$ in $\g_X$ is defined by the presheaf
$\nf_X(U)$ that is formed by the $\Oc_X(U)$-submodule of $\g_X(U)$
that is generated by the $k$-subspace
  \begin{displaymath}
\{\delta \in \g_X (U) : [\delta , \af_X(U)]\subset \af_X(U)\}
  \end{displaymath}
  where $U\subset X$ is an open subset. Clearly, the presheaf $\nf_X$
  is a sheaf satisfying all conditions in \Definition{Lie algebroid}
  except possibly $(2)$. For example, let $X= \Cb^1$ be the complex
  line with coordinate function $t\in \Oc_{\Cb^1}(\Cb^1)$, $\af_X =
  \Oc_{\Cb^1} t\partial_t$, and $\g_X= T_{\Cb^1}=
  \Oc_{\Cb^1}\partial_t$. Then $\nf_{\Cb^1} \cong
  j_!(T_{\Cb^1 \setminus 0})\notin \coh (\Oc_X)$, the extension by $0$ of $T_{\Cb^1
    \setminus 0}$ from $\Cb^1\setminus 0$ to $\Cb^1$.
  
  We put a condition on the $\af_X$-module $N := \g_X/\af_X$ that
  ensures that $\nf_X$ be coherent over $\Oc_X$.
\begin{lemma}\label{normalizer} 
  \begin{enumerate}
  \item If the action of $\af_X$ on $\g_X$ is $\Oc_X$-linear, i.e.
    $[a, \phi \delta]= \phi[a,\delta]$ when $a\in \af_X, \delta \in
    \g_X, \phi \in \Oc_X$, then $\nf_X$ is a sub-Lie algebroid.
  \item Assume that each point in $X$ has a neighbourhood $U\subseteq
    X$ such that $\af_U$ is generated by sections
    $\{a_i\}_{i=1}^r\subset \af_U(U)$ with the property: the
    $\Oc_U$-module $\Oc_U N^{a_i}$ that is generated by the subsheaves
    of $k$-vector spaces 
    \begin{displaymath}
      N^{a_i} := \{n\in N : [a_i,n]=0\}\subset    N_U
    \end{displaymath}
    is locally generated by the sections of $N^{a_i}$.  Then $\nf_X$
    is a sub-Lie algebroid of $\g_X$.
  \end{enumerate}
\end{lemma}
(1) implies, not surprisingly, that the normalizer is 
always coherent when $\g_X$ is a coherent $\Oc_X$-Lie algebra.
\begin{pf} 
  (2): We need to prove that $\nf_X$ is coherent and since this is a
  local problem we may assume that $U=X$, so the $a_i$ are
  global generators of $\af_X$ and $\Oc_XN^{a_i}$ is
  generated by $N^{a_i}(X)$. Since
  \begin{displaymath}
 \{\delta \in \g_X: [\delta, a_i]\in \af_X, i=1,\dots , k\}=
    \cap_{i=1}^k \{\delta\in \g_X: [\delta, a_i]\subset \af_X\}
  \end{displaymath}
  it follows that $ \nf_X = \cap_{i=1}^k \nf^{(i)}_X $ where
  $\nf^{(i)}_X $ is the $\Oc_X$-submodule of $\g_X$ that is
  generated by $\{\delta\in \g_X: [\delta, a_i]\subset
  \af_X\}$.  Hence, by \cite[0. \S 5, Cor.  5.3.6]{EGA1} it
  suffices to see that $\nf^{(i)}_X$ is coherent.  As
  $\af_X,\g_X\in \coh (\Oc_X)$ [loc.\  cit., Prop. 5.3.2] implies
  that $N = \g_X/\af_X\in \coh(\Oc_X)$ and if $\bar
  \nf^{(i)}_X$ is the image of $\nf^{(i)}_X$ in the canonical
  morphism $\g_X \to N$, since $\Ker (\g_X \to N)= \af_X \in
  \coh (\Oc_X)$, it suffices to see that $\bar \nf^{(i)}_X\in
  \coh (\Oc_X)$. But since $\bar \nf^{(i)}_X = \Oc_XN^{a_i}$ is
  generated by its sections $N^{a_i}(X)\subset \bar \nf^{(i)}(X)$
  and $\Oc_X$ is noetherian \Th{noetherianrings} this implies
  that $\bar \nf^{(i)}_X\in \coh (\Oc_X)$. 
  
  (1): Again it suffices to prove that $\bar \nf^{(i)}_X\in \coh
  (\Oc_X)$. Since the problem is local on $X$ and $\g_X$ is coherent
  we may assume that $\g_X$ has a finite set of global generators
  $\{\delta_i\}\subset \Gamma (X, \g_X)$.  Let $\bfr_X$ be the
  coherent submodule of $N$ that is generated by the image of
  $[\delta_i,a]$ in $N(X)$. We have
\begin{displaymath}
\nf^{(i)}_x = \{\sum_i\phi_i\delta_i \in \g_x : \sum \phi_i
[a,\delta_i] \in \af_x\}
\end{displaymath}
showing  that the image $\bar \nf^{(i)}_X = \Ann (\bfr_X)\in \coh (\Oc_X)$ \cite[0. 5.3.10]{EGA1}.
\end{pf}

A Lie algebroid $(\Oc_X, \g_X , \alpha)$ is {\it affine}\/ if $X$ is
an affine scheme over a field $k$; alternatively we denote it by $(A,
\g_A, \alpha)$ when $X= \Spec A$.

A {\it connection} \/ on $\g_X$ is an $\Oc_X$-homomorphism $\nabla :
T_X \to \g_X$ such that $\alpha \circ \nabla = \id_{T_X}$.  The {\it
  curvature} \/ of $\nabla$ is the $\Oc_X$-linear homomorphism
$R_{\nabla} : T_X\wedge T_X \to \bfr$, $\delta \wedge \eta \mapsto
[\nabla(\delta), \nabla (\eta)] - \nabla ([\delta , \eta])$.  The
connection is integrable if $R_\nabla =0 $, hence an integrable
connection is the same as a homomorphism $T_X \to \g_X$. The
connections on $\g_X$ form a sheaf in an obvious way, defining a
torsor over the group $Hom_{\Oc_X}(T_X , \bfr)$.

\begin{example}\label{examples-lie}
  Assume that a Lie algebra $\g$ over $k$ acts on a space $X$ by a
  Lie homomorphism
\begin{displaymath}
  \alpha_0 :\g \to T_X.
\end{displaymath}
Then $\g_X = \Oc_X \otimes_k \g$ forms a Lie algebroid $(\Oc_X, \g_X ,
\alpha)$ defining $\alpha : \g_X \to T_X$ by $\alpha =
\id_{\Oc_X}\otimes \alpha_0 $ and the bracket
\begin{displaymath}
  [f\otimes_k \delta_1, g\otimes_k \delta_2] = fg[\delta_1 ,
  \delta_2]+ f \alpha_0(\delta_1)(g)\delta_2 -
  g\alpha_0(\delta_2)(g)\delta_1.
\end{displaymath}
If $\alpha_0$ induces surjective mappings on the tangent spaces
$T_x/\mf_x T_x$, where $\mf_x$ is the maximal ideal of a point $x\in
X$, then $\g_X$ locally has a connection. If $T_X$ is
generated by its sections $T_X(X)$, then $\g_X$ has a globally defined
connection.
\end{example}

Suppose that $\g_X$ is a Lie algebroid which moreover is a right
$\Oc_X$-module such that $\delta \cdot \phi = \alpha (\delta)(\phi)
\delta + \phi \delta$, when $\phi \in \Oc_X$ and $\delta \in\g_X$;
next section contains an example.  One can then define the {\it
  opposite} \/ algebroid $\g^o_X$ of $\g_X$. As sheaf of $k$-algebras
$\g^o_X$ is the same as $\g_X$.  But if $[\cdot , \cdot]^\circ$ is the
Lie bracket on $\g^o_X$ and $\circ$ is the multiplication between
$\Oc_X$ and $\g^o_X$ one has
\begin{eqnarray*}
  [\delta, \eta]^\circ = [\eta, \delta],\\
  \phi \circ \delta = \delta \phi, \quad \delta \circ  \phi = \phi \delta,
\end{eqnarray*}
for $\delta ,\eta \in \g_X$ and $\phi \in \Oc_X$. This is again  a Lie 
algebroid $(\Oc_X, \g^o_X, \alpha^o)$ if one defines $\alpha^o : \g^o_X \to T_X$ by $\delta \mapsto - \alpha(\delta)$.

\subsection{Linear Lie algebroids} 
\label{linear} 
Let $M$ be an $\Oc_X$-module and $\gl_k(M) := End_k(M)$ be its sheaf
of $k$-Lie algebras.  We have an obvious $k$-linear mapping 
\begin{displaymath}
  i : \Oc_X\to \gl_{k}(M),
\end{displaymath}
defining a structure of $\Oc_X$-bi-module on
$\gl_{k}(M)$, by $(f\cdot \phi)(m) = f\phi(m) $ and $(\phi \cdot f)(m)
= \phi (f m)$, $f\in \Oc_X$, $m\in M$.  Let $\gl_{\Oc_X}(M) \subseteq
\gl_k(M)$ be the $\Oc_X$-submodule of $\Oc_X$-linear homomorphisms $M
\to M$.  Clearly $\gl_{\Oc_X}(M)$ is a Lie sub-algebra which commutes
with the commutative sub-algebra $i(\Oc_X)$. Define $\Dc_X^1(M)
\subseteq \gl_k(M)$ the first order matrix differential operators on
$M$ as the sub-sheaf of $k$-linear mappings $\delta$ such that
$[\delta , \gl_{\Oc_X}(M)] \subseteq \gl_{\Oc_X} (M)$. Evidently
$\Dc^1_X(M)$ is a Lie algebra over $k$.

We temporarily forget the condition $(2)$ in \Definition{Lie
algebroid}, so Lie algebroids are not necessarily coherent.
Now the Lie algebra $\Dc_X^1(M)$ is not provided with a natural
mapping $\alpha$ to $T_X$; its $\Oc_X$-submodule
\begin{equation}\label{linearalgebra}
  \cf_X(M) := \{\delta \in \Dc_X^1 (M) : [\delta , i(\Oc_X)]
  \subseteq i(\Oc_X)\}
\end{equation}
also has a Lie bracket, but again there is no mapping $\alpha$.
However, if $M$ is a {\it faithful} $\Oc_X$-module, i.e.\ the mapping
$i$ is injective, we identify $\Oc_X$ with a sub-ring of
$\gl_{\Oc_X}(M)$, and since $[\delta ,\phi] \cdot m = \delta \cdot
(\phi m) - (\phi \delta)\cdot m$, for $\delta\in \cf_X(M)$, $\phi \in
\Oc_X$, $m\in M$, we get a well-defined homomorphism $\alpha :
\cf_X(M) \to T_X$, $\delta \mapsto [\delta, \phi]$. Therefore $(\Oc_X
, \cf_X(M), \alpha )\in \Lie_X$. It is called the {\it linear Lie
  algebroid} \/ of the faithful $\Oc_X$-module $M$.  $\cf_X(M)$ can
also be described as the set of pairs $(\delta , \partial)$ where
$\partial \in T_X$ and $\delta$ is a lift of $\partial$ to an action
on $M$ such that $\delta (\phi m) = \partial (\phi)m + \phi\delta
(m)$. The homomorphism $\alpha : \cf_X(M)\to T_X$ is $(\delta ,
\partial) \mapsto \partial$, and the bracket is $[(\delta , \partial),
(\delta' , \partial')]= ([\delta, \delta'], [\partial, \partial'])$.
If $M$ is locally free, then $\cf_X(M)$ is locally isomorphic as Lie
algebroid to $\gl_{\Oc_X}(M)\oplus T_X$ where the action of $T_X$ is
defined using a local isomorphism $M_\Omega \cong \Oc_\Omega^n $ over
an open subset $\Omega \subset X$. If $X$ is non-singular and $M$ is
quasi-coherent with presentation
\begin{displaymath}
  \oplus_I \Oc_X \xrightarrow{\phi} \oplus_J \Oc_X
  \twoheadrightarrow M,
\end{displaymath}
then $\cf_X(M)$ is locally isomorphic to
\begin{displaymath}
\frac{\{\delta \in \cf_X(\oplus_J \Oc_X): \delta \phi (\oplus_I \Oc_X) \subseteq \phi (\oplus_I \Oc_X)\} }{ \{\delta \in \cf_X(\oplus_J \Oc_X): \delta  (\oplus_J \Oc_X) \subseteq \phi (\oplus_I \Oc_X)\}}.
\end{displaymath}

Let $I$ be an ideal of $\Oc_X$.  Then $\cf_X(I) = \Oc_X \oplus
T_{X}(I)$, where $T_X (I)$ is the sheaf of derivations that
preserve the ideal $I$, and if $\g_X$ is a Lie algebroid we
let $\g_X(I)$ be the sub-algebroid of vectors $\delta$ such
that $\alpha (\delta) \in T_X(I)$. 

The Lie algebroid $\cf_X(M)$ is also a right $\Oc_X$-module, where $\phi \in \Oc_X$ acts
on $\delta \in \cf_X(M)$ as $\delta \cdot \phi = [\delta , \phi_M] +
\phi_M \delta$, and $\phi_M \in \gl_k(M)$ is the multiplication by
$\phi$ on $M$.

\begin{remark}
  Let $M$ be locally free $\Oc_X$-module of rank $r<\infty$. One may
  construct `classical' Lie algebroids $\g_X$ from an alternating or
  symmetric bilinear form $<\cdot , \cdot> : M\otimes_{\Oc_X}M \to
  \Oc_X$ 
\begin{displaymath}\label{classical}
 \g_X =  \{\delta \in \cf_X(M) : \alpha (\delta)<m_1, m_2> = <\delta
   \cdot m_1 , m_2> + <m_1 , \delta\cdot m_2> \}.
\end{displaymath}
Letting $\det M = \wedge^r M$ be the determinant bundle of $M$ there exists a
homomorphism of Lie algebroids $ \tr : \cf_X(M) \to \cf_X(\det
M)$ defined by
 \begin{displaymath}
 \tr (\delta) \cdot m_1\wedge \cdots \wedge m_d = (\delta \cdot
 m_1)\wedge \cdots \wedge m_r + \cdots + m_1 \wedge \cdots
 \wedge (\delta\cdot m_r),
\end{displaymath}
where $m_i\in M$ and $\delta\in \cf_X(M)$, and one may define the
`special' Lie algebroid by $\Sl (M):=\Ker (\tr) \subseteq \Ker (\alpha
: \cf_X(M) \to T_X)$. Note that $[\cf_X(M) , \cf_X(M)] \nsubseteq \Sl
(M)$, which can be compared to $[\gl_{\Oc_X}(M), \gl_{\Oc_X}(M)]
\subset \Sl_{\Oc_X}(M)$.
\end{remark}
For applications of GAGA and checking coherence it is useful to
identify $\Dc_X^1(M)$ with an $\Oc_X$-module of $\Oc_X$-linear
mappings. We recall this identification (see \cite[\S16, Prop. 16.8.8]{EGA4}): Let
$\Delta: X \to X\times_k X$ be the diagonal morphism.  One
identifies sheaves on $X$ with sheaves on $X\times_k X$ whose
support belongs to $\Delta (X)$.  Let $I$ be the ideal of
$\Delta(X)$ and put $P^1_X= \Oc_{X\times_k X}/I^2$, which is
regarded as an $(\Oc_X, \Oc_X)$-bimodule in an obvious way;
similarly $\Omega_X:=I/I^2$ and $\Oc_X$ are regarded as
bimodules.  One has the multiplication homomorphism $P_X^1\to
\Oc_X$, defined by $\phi \otimes_k \phi'\mapsto \phi \cdot
\phi'$, $\phi, \phi' \in \Oc_X$, and a short exact sequence of
bimodules $0 \to \Omega_X \to P_X^1 \to \Oc_X \to 0$. Since
this sequence has a canonical left $\Oc_X$-linear split $d:
P_X^1 \to \Omega_X$, $\phi\otimes \phi' \mapsto \phi \cdot
\dop \phi'$, putting $\Omega_{X/k}(M)=
\Omega_{X/k}\otimes_{\Oc_X}M$ and $P_X^1(M)=
P^1_X\otimes_{\Oc_X}M$ one gets a short exact sequence of left
$\Oc_X$-modules

\begin{equation}\label{can-ex}
  0 \to \Omega_{X/k}(M) \to P_X^1(M) \to M \to 0.
\end{equation}
Then
\begin{displaymath}
  \Dc_X^1(M) = Hom_{\Oc_X}(P_X^1(M), M),
\end{displaymath}
where $\delta \in \Dc_X^1(M)$ is identified with the mapping
$P^1_X(M)\to M$ defined by $\phi \otimes \phi' \otimes m \mapsto \phi
\delta (\phi' \cdot m) $.
\begin{lemma}
  \begin{enumerate}
  \item ($X$ is a space) $P^1_{X/k},
    \Omega_{X/k}\in \coh (\Oc_X)$.
  \item If $M, \Omega_X\in \coh (\Oc_X)$, then $\cf_X(M)\in \coh (\Oc_X)$.
  \end{enumerate}
\end{lemma}
\begin{pf}
  (1) is well-known. (2): $\cf_X(M)$ is the normalizer of $i(\Oc_X)$
  in $\Dc^1_X(M)= Hom_{\Oc_X} (P^1_X, M)\in \coh (\Oc_X)$.  Since the
  action of $i(\Oc_X)$ on $\Dc^1_X(M)$ is $\Oc_X$-linear the assertion
  follows from (1) in \Lemma{normalizer}.
\end{pf}

A connection on $M$ is a $k$-linear mapping $\nabla : M \to
\Omega_{X/k}(M) $ such that $\nabla (\phi m) =
\dop (\phi) m + \phi \nabla(m)$. Let us denote the (sheaf of)
connections on $M$ by $\Der_k(M , \Omega_{X/k}(M))$.
 \begin{lemma}\label{lemgaga}
   $ \Der_k(M , \Omega_{X/k}(M)) = Hom_{\Oc_X}(T_X,
   \cf_X(M))$
\end{lemma}
\begin{proof}
  The mappings
\begin{eqnarray}\nonumber
&\rho :& \Der_k(M , \Omega(M)) \to Hom_{\Oc_X}(T_X,
\cf_X(M)), \\ \nonumber && \rho(\nabla)( \partial) \in
\cf_X(M): m \mapsto <\nabla(m),\partial> ; \\ \nonumber &\psi
:& Hom_{\Oc_X}(T_X, \cf_X(M)) \to \Der_k(M , \Omega(M)), \\
\nonumber & & \psi (\phi)(m) \in \Omega(M) : \partial \mapsto
\psi (\phi)(m)(\partial)= \phi(\partial)(m)
\end{eqnarray}
are mutually inverse isomorphisms.
\end{proof}
Thus if $M$ is a  faithful $\Oc_X$-module then a connection on
$\cf_X(M)$ is the same as a connection on $M$.

The short exact sequence (\ref{can-ex}) defines an element
\begin{displaymath}
c(M)\in \Ext_{\Oc_X}^1(M, \Omega_X(M)),
\end{displaymath}
called the Atiyah class of the $\Oc_X$-module $M$. Thus $c(M)=0$ when
there exists a globally defined connection $\nabla\in \Der_k(M,
\Omega_X(M))$ and if $M$ is locally free, so connections exist
locally, then $c(M)\in \Ext^1_{\Oc_X}(M, \Omega_X(M)) = H^1(X,
\Omega_X\otimes_{\Oc_X}\gl_{\Oc_X}(M))$.  Any vector in $T_X$ defines
locally a homomorphism $\Omega_X(M) \to M$ and thus a push-out of
(\ref{can-ex}), defining a homomorphism of $\Oc_X$-modules $\psi :T_X
\to Ext^1_{\Oc_X}(M,M)$.  Applying $Hom_{\Oc_X}(\cdot , M)$ to the
short exact sequence (\ref{can-ex}) we get long exact sequences
\begin{displaymath}
  \xymatrixnocompile{0 \ar[r]& \gl_{\Oc_X}(M)\ar[r] & \Dc^1(M) \ar[r]&
  Hom_{\Oc_X}(\Omega_{X/k}(M),M) \ar[r] & Ext^1_{\Oc_X}(M,M)
  \ar[r]& \cdots & \\ 0 \ar[r]& \gl_{\Oc_X}(M)\ar[r] \ar@{=}[u] &
  {\cf_X(M) \ar[r]^{\alpha}}\ar@{^{(}->}[u] & T_X \ar[u] \ar[r]^{\psi}&
  Ext^1_{\Oc_X}(M,M)\ar[r]\ar@{=}[u]& \cdots & }
\end{displaymath}
In particular:
\begin{displaymath}
  \Ker (\psi) = \alpha(\cf_X(M)).
\end{displaymath}

There is one obstruction class in $\Ext^1_{\Oc_X}(\alpha(\cf_X(M)),
\gl_{\Oc_X}(M))$ for the existence of a split of the short exact
sequence
\begin{displaymath}
  0 \to \gl_{\Oc_X}(M) \to \cf_X(M) \to \alpha(\cf_X(M)) \to 0
\end{displaymath}
and one obstruction class in
$\Ext^1_{\Oc_X}(\frac{T_X}{\alpha(\cf_X(M))},
\alpha(\cf_X(M)))$ that $\alpha(\cf_X(M))$ be a direct summand
of $T_X$; these two obstruction classes vanish if and only if
$\cf_X(M)$ has a connection.

\subsection{Differential operators}
\label{diffoperators} 
We define differential operators as follows.
\begin{definition}\label{diffop}
  A ring of differential operators $(\Oc_X , \Dc_X)$ is a $k$-algebra
  $\Dc_X$ on a space $(X,\Oc_X)$ that has:
  \begin{enumerate}
  \item a filtration, $0= \Dc^{-1} \subseteq\Dc_X^{0} \subseteq
\Dc_X^{1} \subseteq \cdots \subseteq \Dc_X$ such that
$\Dc_X=\cup_n \Dc_X^n$ and the associated graded ring $\gr
\Dc_X := \oplus_{n\geq0} \Dc_X^n/\Dc_X^{n-1}$ is commutative.
  \item an inclusion $i : \Oc_X \hookrightarrow \Dc^0_X$ such
that $[\Dc_X^n, i(\Oc_X)] \subseteq \Dc_X^{n-1}$.
  \end{enumerate}
\end{definition}
Let $\Diff_X$ be the category whose objects are rings of
differential operators on $X$ and the homomorphisms are required to be
compatible with the filtrations. Any ring of differential
operators defines a Lie algebroid 
\begin{equation}\nonumber
\g_X = \{\delta \in \Dc^1_X : [\delta , i(\Oc_X)] \subset i(\Oc_X)\}
\end{equation} 
(compare to (~\ref{linearalgebra}) in Subs.~\ref{linear})
defining a functor $Lie :\Diff_X \to \Lie_X$. The functor $Lie$
has a left adjoint $\Dc_X(\cdot): \Lie_X \to \Diff_X$, where
$\Dc_X:=\Dc_X(\g_X)$ is called the enveloping ring of
differential operators of $(\Oc_X,\g_X,
\alpha )\in \Lie_X$; it is  constructed similarly to the enveloping algebra of a Lie algebra.

Let $\So_X = \So_X(\g_X)$ be the symmetric algebra of the
$\Oc_X$-module $\g_X$.

In a ringed space $(X, R_X)$ we recall that $R_X$ is {\it
noetherian}\/ if the following holds: $(i)$ $R_X$ is coherent as
left $R_X$-module; $(ii)$ for each point $x$ in $X$ the stalk
$R_x$ is a noetherian ring; $(iii)$ for any open subset
$\Omega$ of $X$ any increasing family of coherent sub-modules
of a coherent $R_\Omega$-module is locally stationary.

\begin{thm}\label{noetherianrings}($X$ is a space)
  The rings $\Oc_X$, $\So_X$, $\gr (\Dc_X)$ and $\Dc_X$ are
  noetherian.
\end{thm}
\begin{pf}
  $\Oc_X$ is noetherian: When $X$ is a variety see \cite[Cor.
  1.5.3]{EGA1} (it suffices that $(X,\Oc_X)$ be locally noetherian).
  When $X$ is a complex analytic space $(i)$ is Oka's theorem, see
  \cite{gunning:2}; a proof of $(i)$ and $(ii)$ can be found in
  \cite{gunning:2}. $(iii)$ is due to Serre \cite{serre:prolongement}
  and Grauert \cite[Th 8]{grauert}.
  
  $\So_X$, $\gr (\Dc_X)$, and $\Dc_X$ are noetherian: Follow the
  proofs in B. Kaup's article in \cite{borel:Dmod}. Note that [loc.\ 
  cit., Prop.  3.5] should be changed to: the canonical homomorphism
  $\So_X \to \gr (\Dc_X)$ is surjective. See also \cite{bjork:analD}.
\end{pf}

The image of the canonical inclusion $\sigma_s : \g_X \to
\So_X$ generates the algebra $\So_X $, making it possible to
define a $k$-linear Lie product
\begin{displaymath}
  \{\cdot , \cdot \}: \So_X \otimes_k \So_X \to \So_X
\end{displaymath}
by $\{\sigma_s(\partial) , \sigma_s(\eta)\} =
\sigma_s([\partial , \eta]) $, $\{\phi , \psi\} = 0$,
$\{\sigma_s(\partial) , \phi \}= \alpha (\partial)(\phi)$, and
inductively by the rules
\begin{enumerate}
\item $\{a,b \} = - \{b,a\}$
\item $\{a, bc\} = \{a,b\}c + \{a, c\}b$
\end{enumerate}
for all $a,b,c \in \So_X$. This defines a Poisson algebra $(\So_X
,\cdot , \{\cdot , \cdot \})$. In the same way one defines a bracket
$\{\cdot , \cdot \}$ on $\gr \Dc_X$ resulting in a Poisson algebra
$(\gr \Dc_X , \cdot , \{\cdot , \cdot \})$.  The canonical mapping
$\g_X \to \gr (\Dc_X)$ lifts to a surjective homomorphism of Poisson
algebras $p:\So_X \to \gr (\Dc_X)$. We have the following
generalization of   the Poincar\'e-Birkhoff-Witt theorem:
\begin{thm}[\cite{rinehart:difforms}]\label{rinehart}
  If $\g_X$ is locally free over $\Oc_X$ then $p$ is an isomorphism.
  \end{thm}
Let $M$ be a coherent $\g_X$-module (see \Definition{gmodule}
below), $\Omega\subset X$ an open set and $M^0_\Omega$ a
coherent $\Oc_\Omega$-module generating $M_\Omega$ over
$\g_\Omega$. Set $M^n_\Omega = \g_\Omega M^{n-1}_\Omega =
\Dc_\Omega^n \cdot M^0_\Omega$ for $n\geq 1$, and
\begin{displaymath}
G (M_\Omega) = \oplus M^n_\Omega/M^{n-1}_\Omega,
\end{displaymath}
defining a coherent $\So_\Omega$-module.  When $X$ is a variety we let
$\Specb \So_X$ be the spectrum of $\So_X$, and if $X$ is a complex
analytic space we use the same notation $\Specb \So_X$ to denote the
analytic space that can be associated to $\So_X$; note that the
$\Oc_X$-algebra $\So_X$ is locally finitely generated. There is a
morphism of spaces $p: \Specb \So_X \to X$ such that $p^*(\So_X)=
\Oc_{\Specb \So_X}$.  The support $\SSo (M_\Omega)\subseteq \Specb \So_\Omega$ of
$p^*(G(M_\Omega))$ is independent of the choice of $M^0_\Omega$ (see
e.g. \cite{bjork:analD}), and therefore the construction globalises
defining a conic set $\SSo (M_X)\subset \Specb \So_X $, which is a
subspace since $p^*(G(M_\Omega))$ is coherent over $\Oc_{\Specb
  \So_X}$.  The set $\SSo (M_X)$ is the {\it singular support} \/ of
$M$, and the ideal $\Jo (M)$ of $\SSo (M_X)$ is the characteristic
ideal of $M$.

\begin{thm}[{\cite[Th. I]{gabber:integrability}}]
\label{gabber}
  \begin{displaymath}
    \{\Jo (M), \Jo (M)\} \subseteq \Jo (M).
  \end{displaymath}
\end{thm}

Let $p : \Specb \So_X \to X$ be the canonical projection and
$i: X \to \Specb \So_X $ the canonical section. It is easy to
see that $M_x$ is of finite type over $\Oc_x$ if and only if
$\SSo (M_x) \subset i(x)$, i.e.\  the ideal $\Jo(M_x)=\Jo(M)_x$ of $\SSo (M_x)$ is 
$(\sigma (\g_x)) \subseteq \So (\g_x)$.  Set $\dSSo (M) = \SSo (M) \setminus i(X)$.

\begin{definition}\label{strong-supp}
  The {\it strong support} \/ of a $\g_X$-module $M$ is the set
    \begin{displaymath}
      \ssupp M = \{x \in X : M_x \text{ is not of finite type
      over } \Oc_x \} = p(\dSSo (M)).
    \end{displaymath}
\end{definition}

We shall need the following lemma, whose proof is well-known.
\begin{lemma}\label{strong-var}($X$ is a  space)
  Let $M$ be a coherent $\g_X$-module. Then $\ssupp M$ is a subspace
  of $X$.
\end{lemma}
\begin{pf} We prove this only when  $X$ is a complex analytic space.
  The set $\SSo (M) = \supp \Go (M)$ is a conic analytic set, hence
  $\dSSo (M)$ is a conic analytic set. Let $\Projb (\So_X)$ be the
  projective space of $\Specb (\So_X)$ (the analytic space that is
  associated with the graded $\Oc_X$-algebra $\So_X $), with canonical
  projection $p_1: \Projb (\So_X ) \to X$. So we have an open
  inclusion mapping $\phi : \Specb (\So_X ) \setminus i(X) \to \Projb
  (\So_X )$, satisfying $p_1 \circ\phi =p' $, where $p'$ is the
  restriction of $p$ to $\Specb (\So_X )\setminus i(X)$. That $\dSSo$
  is conic means that $\dSSo (M) = \phi^{-1} (V)$ where $V$ is a
  complex analytic space in $\Projb (\So_X )$. Hence
  \begin{displaymath}
    \ssupp (M) = p'(\dSSo (M))= p_1 (V),
  \end{displaymath}
  so the result follows from Remmert's proper mapping theorem
  \cite[Ch. N Th. 1]{gunning:2}.
\end{pf}

\subsection{Homomorphisms}
\label{morphisms}
\subsubsection{Remark on the tangent mapping}\label{tangentmap}
Let $\pi : Y \to X$ be a morphism of spaces.  There is then the
canonical exact sequence $\pi^*(\Omega_{X/k}) \to \Omega_Y \to
\Omega_{Y/X} \to 0$, inducing a homomorphism
\begin{displaymath}
  T_Y=Hom_{\Oc_Y}(\Omega_{Y/k}, \Oc_Y) \to Hom_{\Oc_Y}(\pi^*(\Omega_{X/k}), \Oc_Y).
\end{displaymath}
Consider also the canonical homomorphism
\begin{displaymath}
\psi : \pi^*(T_X) = \pi^*(Hom_{\Oc_X} (\Omega_{X/k}, \Oc_X))
\to Hom_{\Oc_Y}(\pi^*(\Omega_{X/k}), \Oc_Y).
\end{displaymath}
If the induced homomorphism between stalks at any point $x\in X$ is an
isomorphism, then $\psi$ is an isomorphism. As $\Omega_{X/k}\in \coh
(\Oc_X)$ we have 
\begin{displaymath}
  Hom_{\Oc_X}(\Omega_{X/k},\Oc_X)_x = Hom_{\Oc_x}(\Omega_{\Oc_x/k}, \Oc_x),
\end{displaymath}
so the condition is that the
canonical morphism
\begin{displaymath}
  \Oc_y \otimes _{\Oc_x} \Hom_{\Oc_x}(\Omega_{\Oc_x/k}, \Oc_x)
  \to Hom_{\Oc_y}(\Oc_y \otimes_{\Oc_x} \Omega_{\Oc_x/k},
  \Oc_y)
\end{displaymath}
is an isomorphism. Here $\Oc_x \to \Oc_y$ is the homomorphism of local
rings defined by $\pi$. Thus $\psi$ is an isomorphism if either $X$ is
non-singular, so $\Omega_{\Oc_x/k}$ is free of finite rank over
$\Oc_x$, or if $\pi$ is flat (to see this use also the fact that
$\Omega_{\Oc_x/k}$ is finitely presented over $\Oc_x$).  When $\psi$ is
an isomorphism  there exists a canonical homomorphism
\begin{displaymath}
   d\pi : T_Y \to \pi^*(T_X).
\end{displaymath}
{\it Also when $\psi$ is not an isomorphism we shall abuse the
notation in the following by writing $\pi^*(T_X)$ when we
actually mean $Hom_{\Oc_Y}(\pi^*(\Omega_{X/k}), \Oc_Y)$}.  Thus
the homomorphism $d\pi$ is always defined.

\subsubsection{Homomorphisms of Lie algebroids} 
The definition of homomorphisms $\g_Y \to \g_X$ between Lie algebroids
on different spaces is slightly involved, because $\pi^*(\g_X)= \Oc_Y
\otimes_{\pi^{-1}(\Oc_X)}\pi^{-1}(\g_X) $ does not in general have a
natural structure of Lie algebroid.

\begin{remark}\label{seidenberg}
  The $\Oc_Y$-module $\pi^*(\g_X)$ has a structure of Lie algebroid
  when there exists a homomorphism of Lie algebroids over
  $\pi^{-1}(\Oc_X)$, $\psi : \pi^{-1}(T_X) \to T_Y$, defining the Lie
  bracket by
  \begin{eqnarray}\nonumber
    [\phi_1\otimes \partial_1, \phi_2\otimes \partial_2] & := &
    \phi_1\phi_2\otimes [\partial_1, \partial_2]+ \phi_1\psi
    (\alpha(\partial_1)(\phi_2))\otimes\partial_2
    -{}\\\nonumber & & {}- \phi_2\psi
    (\alpha(\partial_2)(\phi_1))\otimes\partial_1.
  \end{eqnarray}
  Here are two instances of this situation:
  \begin{enumerate}
  \item $\pi$ is a projection, for then $\pi^*(T_X)\subset T_Y$.
  \item Let $\pi :X' \to X$ be the normalization morphism of an
    integral space.  By \cite{seidenberg} every section of $T_X$ has a
    unique lift to a section of $T_{X'}$, i.e. \ if $U \subset X'$ is
    an open subset and $\delta \in \pi^{-1}(T_X)(U)$, then there
    exists a unique section $\delta' \in T_{X'}(U)$ such that its
    restriction to $\pi^{-1}(\Oc_X)(U) \subset \Oc_{X'}(U)$ coincides
    with $\delta$; see {\it loc.\ cit.}  \/ for a more general and
    precise condition.
\end{enumerate}
\end{remark}

\begin{definition}\label{lie-morphism} Let $(Y, \g_Y, \alpha) \in \Lie_Y, (X,
  \g_X , \beta)\in \Lie_X$.  A homomorphism $(\pi',\pi) :(Y,
  \g_Y, \alpha) \to (X, \g_X , \beta)$ is:
  \begin{enumerate}
  \item a morphism $\pi : Y \to X$ of spaces;
  \item a homomorphism of $\Oc_Y$-modules
  \begin{displaymath}
    \pi' :\g_Y \to \pi^*(\g_X)
  \end{displaymath}
  such that the following diagram commutes
  \begin{equation}\nonumber
\begin{CD}
  \g_Y @>\pi'>> \pi^*(\g_X) \\ @V\alpha VV @VV \beta V \\ T_Y
  @>d\pi >> \pi^*(T_X)
\end{CD}
\end{equation}
   \end{enumerate}
   and $\pi'$ satisfies:
\begin{displaymath}
  \pi' ([\delta_1,\delta_2])= \sum a_ib_j \otimes [\eta_i,
  \eta_j] + \sum (\alpha(\delta_1)(b_i) -\alpha(\delta_2)(a_i))
  \otimes \eta_i
\end{displaymath}
if $\pi' (\delta_1) = \sum a_i \otimes \eta_i$, $\pi'(\delta_2)
= \sum b_j \otimes \eta_j$, where $\delta_1 , \delta_2 \in
\g_Y$, $a_i, b_i\in \Oc_Y$, $\eta_i \in \g_X$.  (Recall our
abuse of notation in defining $d\pi$.)
\end{definition}
If $\pi = \id $ then $\pi'$ is a homomorphism of Lie algebroids
on $X$. If $\alpha = 0$ then $\g_Y$ is a sheaf of Lie algebras
over $\Oc_Y$ and a homomorphism $\g_Y \to \g_X$ is a
homomorphism of the $\Oc_Y$-Lie algebra $\g_Y$ onto a Lie
algebra in the $\Oc_Y$-module $\pi^*(\g_X)$.

\subsubsection{Pull-back} 
\label{pullbacklie}
Let $(\Oc_X , \g_X , \alpha) \in \Lie_X$. The fibre product
with respect to the $\Oc_Y$-linear homomorphisms
$\pi^*(\alpha): \pi^*(\g_X)\to \pi^*(T_X) $ and $d \pi : T_Y
\to \pi^*(T_X)$ is the $\Oc_Y$-module
\begin{displaymath}
  \pi^*(\g_X)\times_{\pi^*(T_X)} T_Y = \{(\delta ,\partial) \in
  \pi^*(\g_X)\times T_Y : \pi^*(\alpha)(\delta) = d\pi
  (\partial)\}.
\end{displaymath}
We define a bracket on $ \pi^*(\g_X)\times_{\pi^*(T_X)} T_Y$ by
\begin{eqnarray}\nonumber
\lefteqn{ [(\phi_1 \otimes \delta_1, \partial_1) , (\phi_2
\otimes \delta_2, \partial_2)] } \\ \nonumber & = &
(\partial_1(\phi_2)\otimes\delta_2 - \partial_2(\phi_1)\otimes
\delta_1+ (\phi_1\phi_2)\otimes [\delta_1, \delta_2],
[\partial_1, \partial_2]),
\end{eqnarray}
for $(\phi_i\otimes\delta_i, \partial_i) \in
\pi^*(\g_X)\times_{\pi^*(T_X)} T_Y$, and letting $\beta :
\pi^*(\g_X)\times_{\pi^*(T_X)} T_Y \to T_Y $ be the projection on the
second factor we get a structure of Lie algebroid on
$\pi^*(\g_X)\times_{\pi^*(T_X)} T_Y$; in that capacity we denote it by
$\pi^+(\g_X)$. Clearly, we get a functor $(\Oc_X , \g_X , \alpha)
\mapsto (\Oc_Y, \pi^+(\g_X), \beta ) $ from $\Lie_X$ to $\Lie_Y$;
being a composition of a right exact and left exact functor, $\pi^+$ is in
general neither left nor right exact. Let $g : \pi^+(\g_X) \to
\pi^*(\g_X)$ be the projection on the first factor.  Then if $\phi :
\g_Y \to \g_X$ is a homomorphism of Lie algebroids there exists a
unique homomorphism $f:\g_Y \to \pi^+(\g_X)$ such that $\phi = g\circ
f$; hence the functor $\g_X \mapsto \pi^+(\g_X)$ represents the
functor $\Lie_Y \to \Set$, $\g_Y \mapsto Hom_{\Lie}(\g_Y , \g_X)$. We
call $\pi^+(\g_X)$ the pull-back of $\g_X$ by $\pi$. 

If $Y \xrightarrow{f\ }Z \xrightarrow{g\ } X$ are morphisms of spaces,
we have a canonical homomorphism
\begin{eqnarray}\label{pullbackmorphism}
&c_{fg}: f^+(g^+(\g_X)) \to   (g\circ f)^+(\g_X),\\ \nonumber
&(h(y)\otimes_{f^{-1}(\Oc_Z)}(\phi(z)\otimes_{g^{-1}(\Oc_X)}\delta (x),
\partial_z ), \partial_y) \mapsto (h(y)\bar \phi(y)\otimes_{(g\circ f)^{-1}(\Oc_X)} \delta (x), \partial_y),
\end{eqnarray}
where $h(y)\in \Oc_Y$, $\phi (z)\in \Oc_Z$, $\delta (x)\in \g_X$,
$\partial_z\in T_Z$, $\partial_y\in T_Y$, and $\bar \phi (y)$ is the
image of $\phi (z)$ in $\Oc_Y$.  Unfortunately, $c_{fg}$ is not always
an isomorphism, so the projection functor $\Lie \to \Sch$ makes $\Lie$
only into a prefibred category over the category of schemes/k (or over
the category of complex analytic spaces); see \cite[Exp. VI]{SGA1} for
this terminology. The following lemma is slightly more general than
\cite[Lemma 1.5.1]{bei-ber:jantzen}.

\begin{lemma}\label{prefib-iso}
If either  of the conditions 
\begin{enumerate}
\item $f$ is flat, or
\item $g$ is smooth and $Z,X$ are non-singular
\end{enumerate}
are satisfied, then $c_{fg}$ is an isomorphism.
\end{lemma}
\begin{pf} The problem is local on $Y,Z, X$ so we may
  instead consider a composition of local rings $A \xrightarrow{ g \ }
  B \xrightarrow {f \ } C$, and an affine Lie algebroid $(A, \g_A)$.
  Put $\g_B = g^+(\g_A) = B\otimes_A \g_A \times_{B\otimes_A T_A} T_B$
  and $\g_C = (g\circ f)^+(\g_A)= C\otimes_A\g_A
  \times_{C\otimes_AT_A} T_C$. Then $f^+(g^+(\g_A))= f^+(\g_B) =
  C\otimes_B \g_B \times_{C\otimes_B T_B} T_C$, and we have a
  canonical homomorphism
  \begin{displaymath}
    C\otimes_B \g_B \times_{C\otimes_B T_B} T_C \to C\otimes_A\g_A  \times_{C\otimes_AT_A} T_C.
  \end{displaymath}
When necessary we here abuse the notation as in (\ref{tangentmap}).

Assume $(1)$. Then the exact sequence $0 \to \g_B
\to B\otimes_A\g_A \times T_B$ gives the exact sequence $0 \to
C\otimes_B \g_B \to C \otimes_A \g_A \times C\otimes_B T_B$, implying 
that  $C\otimes_B \g_B \cong C\otimes_A \g_A \times_{C\otimes_A T_A}
C\otimes_B T_B$. Therefore
\begin{eqnarray*}
  f^+(\g_B) &=& C\otimes_B \g_B \times_{C\otimes_B T_B} T_C \\
&\cong&  (C\otimes_A \g_A \times_{C\otimes_A T_A} C\otimes_B T_B)  \times_{C\otimes_B T_B} T_C \\
  &\cong& C\otimes_A \g_A \times_{C\otimes_A T_A} T_C \\ &=& (g\circ
  f)^+(\g_A).
\end{eqnarray*}

Assume $(2)$. Then we have a split short exact sequence
\begin{displaymath}
  0 \to T_{B/A}\to T_B \to  B\otimes_A T_A \to 0.
\end{displaymath}
Choose a splitting, so $T_B = T_{B/A} \oplus B\otimes_A T_A$. Then,
using the fact that $N\times_L(M+L) = N\times M$  when $N,L,M$ are
$B$-modules with morphisms $N \to L$, $M +L \to L$ and $M$ is the
kernel of the last morphism, one gets
\begin{displaymath}
\g_B= B\otimes_A \g_A \times_{B\otimes_A T_A} (T_{B/A}\oplus B\otimes AT_A) \cong T_{B/A} \oplus B\otimes_A \g_A
\end{displaymath}
where the morphism $T_{B/A} \oplus B\otimes_A \g_A \to T_B$ is $\delta 
+ b\otimes \eta \mapsto \delta + b \otimes \alpha(\eta) \in T_{B/A}\oplus B\otimes_A T_A = T_B$.
Therefore
\begin{eqnarray*}
    f^+(\g_B) &=&  C\otimes_B (T_{B/A} \oplus B\otimes_A
    \g_A)\times_{C\otimes_B T_B} T_C \\
&=& (C\otimes_B T_{B/A} \oplus    C\otimes_A \g_A) \times_{C\otimes_B
  T_{B/A} \oplus C\otimes_A T_A}T_C \\
&=& C\otimes_A \g_A \times_{C\otimes_A T_A} T_C\\
&=& (g\circ f)^+(\g_A).
\end{eqnarray*}
\end{pf}

\begin{examples}\label{ex-pullback}
  \begin{enumerate}
  \item $ \pi^+(T_X)= T_Y$. When $\alpha =0$, then $\pi^+(\g_X)=
    \pi^*(\g_X)$. 
  \item If $\pi : Y \to X$ is an open embedding, then $\pi^+(\g_X)=
    \pi^*(\g_X)$.
  \item Let $\pi : Y \to X$ be a closed embedding to a non-singular space
  $X$, and  $I_Y$ be the ideal of $\pi(Y)$. Then
  \begin{displaymath}
    \pi^+(\g_X)= \pi^+(\g_X(I_Y)) = \pi^{-1}(\frac
    {\g_X(I_Y)}{I_Y\g_X}),
  \end{displaymath}
  since $T_Y = \pi^{-1}(T_X(I_Y))/I_YT_X)$. Note that
  $\pi^+(\g_X)$ need not be locally free when $\g_X$ is locally 
  free even if $Y$ is non-singular.
\item If $\alpha(\g_X)$ is flat over $\Oc_X$, and $i_x: \{x\} \to X$
  is an inclusion of a point in $X$, then $i_x^+(\g_X)= i^*(\bfr_X)=
  \bfr_x/\mf_x \bfr_x$.
\item If $X$ is non-singular and $\pi$ is a submersion, i.e. $d\pi :
  T_Y \to \Oc_Y\otimes_{\pi^{-1}(\Oc_X)} \pi^{-1}(T_X)$ is surjective,
  then there exists a canonical injection $\pi^*(\g_X) \hookrightarrow
  \pi^+(\g_X)$.
  \end{enumerate}
\end{examples}

Clearly, if $\g_X$ is transitive then $\pi^+(\g_X)$ again is
transitive, and if $\alpha =0$ then $\pi^+(\g_X) =
\pi^*(\g_X)$ is locally free. We will later need the following
lemma.
\begin{lemma}\label{pullfree} Let $\pi : Y\to X$ be morphism of non-singular
  spaces and $\g_X$ a locally free and transitive Lie
  algebroid. Then $\pi^+(\g_X)$ is locally free and transitive.
\end{lemma}
\begin{pf}
  Factorize $\pi$ as $\pi = p\circ i$ where $i: Y \to Y\times
  X$ is the graph morphism and $p : Z=Y\times X \to X$ is the
  projection on the second factor. Since $p$ is smooth and $Z$
  and $X$ are non-singular $\pi^+(\g_X) = i^+(p^+(\g_X))$
  \Lem{prefib-iso}.  
  
  One has $\g_Z:=p^+(\g_X)= p^*(\g_X)\times_{p^*(T_X)}T_Z =
  p^*(\g_X)\oplus q^*(T_Y)$ where $q: Z \to Y$ is the projection
  to the first factor, so this is locally free and transitive.
  
  Letting $i: Y \to Z$ be a closed embedding and $\g_Z$ a locally free
  and transitive Lie algebroid it remains to see that $i^+(\g_Z)$ is
  locally free. There exists a short exact sequence
\begin{displaymath}
  0 \to i^+(\g_Z) \to i^*(\g_Z) \to \frac
  {i^*(\g_Z)}{i^+(\g_Z)} \to 0
\end{displaymath}
where $i^*(\g_Z)$ is locally free. Thus it suffices to prove that $\frac
{i^*(\g_Z)}{i^+(\g_Z)}$ is locally free, for  then $i^+(\g_Z)$ is locally a
direct summand of a free module; hence $i^+(\g_Z)$ is locally free.  $Y$
is non-singular  so there exists a split  short exact sequence
\begin{displaymath}
  0 \to T_Y  \to i^*(T_Z) \to \frac {i^*(T_Z)}{T_Y }\to 0;
\end{displaymath}
hence $\frac {i^*(T_Z)}{T_Y}$ is locally free (the normal bundle of
$Y$ in $Z$). Now by transitivity $\frac {i^*(\g_Z)}{i^+(\g_Z)}\cong
\frac {i^*(T_Z)}{T_Y}$, implying the assertion.
\end{pf}

\begin{lemma}($X$ is a non-singular space)
  Let $M$ be a locally free $\Oc_X$-module. Then $ \pi^+(\cf_X(M))
  \cong \cf_Y(\pi^*(M))$.
\end{lemma}
\begin{pf}  In \Section{pullbackmod} we give a  canonical homomorphism $\rho :\pi^+(\cf_X(M))\to
  \cf_Y(\pi^*(M))$, which we now contend is an isomorphism. This is a
  local question on $X$ so we may assume that $M$ is free, and hence
  that $\cf_X(M)= Hom_{\Oc_X}(M,M)\oplus T_X$. Then
  \begin{eqnarray*}
    \pi^+(\cf_X(M)) &=&  \pi^+(Hom_{\Oc_X}(M,M)\oplus T_X)
\\ &=&  \pi^*(Hom_{\Oc_X}(M,M)\oplus T_X)\times_{\pi^*(T_X)}T_Y 
\\ &=& \pi^*(Hom_{\Oc_X}(M,M))\oplus T_Y  
\\ &=& Hom_{\Oc_Y}(\pi^*(M),  \pi^*(M))\oplus T_Y  \quad \quad \text{($M$ is  locally free)}
\\ &=& \cf_Y(\pi^*(M)).
  \end{eqnarray*}
  It remains to the check that these isomorphisms compose to $\rho$,
  but we omit the details.
\end{pf}

\subsection{Picard Lie algebroids}
\label{picardliealg}
Let $(\Oc_X, \g'_X, \alpha^1)\in \Lie_X$, put $\bfr= \Ker (\alpha^1)$
and let $\chi: \bfr \to \Oc_X$ be a homomorphism of Lie algebroids (a
character of $\bfr$). The {\it push-out}\/ of $\g'_X$ by $\chi$ is the
Lie algebroid $\g^{\chi}_X = \{(\phi,\delta)\in \Oc_X \oplus
\g'_X\}/J$, where $J=\{(\chi(b), -b):\in \Oc_X\oplus \bfr \}$.  A
$\g'_X$-module $M$ such that the action of $\bfr$ is determined by $b
\cdot m = \chi (b)m$, $b\in \bfr$, $m\in M$, corresponds to a
$\g^\chi_X$-module.

\begin{definition}\label{defpic}
  A Lie algebroid $(\Oc_X, \g_X , \alpha)$ is a Picard Lie algebroid
  if $\alpha$ is surjective and $\Ker (\alpha) \cong \Oc_X$. We {\it
    identify} $\Ker (\alpha) = \Oc_X$ by choosing a global central
  section $1_{\g_X} \in \g_X$.  Let $\LPic_X$ be the category of
  Picard Lie algebroids $\g_X$ which are locally isomorphic to
  $T_X\oplus \Oc_X$ in the category of $\Oc_X$-modules (this is of
  course automatic when $X$ is non-singular).
\end{definition}

Thus if $(\Oc_X , \g_X , \alpha)\in \Lie_X$ is transitive and
$\chi$ is a character of $\Ker (\alpha)$ then the push-out
$\g^\chi_X \in \LPic_X$.

Homomorphisms of Picard Lie algebroids $\phi : \g_1 \to \g_2$ are
required to satisfy 
\begin{displaymath}
  \phi(1_{\g_1}) = 1_{\g_2} = 1_{\Oc_X},
\end{displaymath}
so they are isomorphisms, i.e.  $\LPic_X$ is a groupoid. Note that in
the category $\Lie_X$ we do not require $\phi (1_{\g_1})= 1_{g_2}$, so
$\LPic_X$ is not a full sub-category of $\Lie_X$.

We shall determine the cohomology group classifying $\LPic_X$.
In general we have a short exact sequence
\begin{displaymath}
  0 \to \Oc_X \to \g_X \xrightarrow{\alpha} T_X \to 0
\end{displaymath}
which is locally split as $\Oc_X$-modules, locally identifying $\g_X$
with $\Oc_X \oplus T_X$. 

{\it Local structure}\/ : Suppose that $X_0\subset X$ is an affine subset over which the exact
sequence splits; thus putting $\g= \Gamma (X_0, \g_{X_0})$ one has a
short split exact sequence in the category of $A$-modules
\begin{displaymath}
  0 \to A \to \g \to \Der_k(A) \to 0.
\end{displaymath}
A connection $\nabla : \Der_k(A) \to \g$,
determines an isomorphism $\g \cong A \oplus \Der_k(A)$ (as
$A$-modules). Its curvature is the $A$-linear mapping $R_
\nabla : \Der_k(A)\wedge \Der_k(A) \to A$, $\partial_1 \wedge
\partial_2 \mapsto \eta(\partial_1 \wedge \partial_2) = [\nabla
(\partial_1),\nabla (\partial_2)]- \nabla([\partial_1 , \partial_2])$,
 defining an element $\eta$ in the vector space of closed 2-forms
$\Omega^{2,cl}(A)$; the condition that $\eta$ be closed follows from
the Jacobi identity.  Hence, identifying the $A$-module $\g$ with
$A\oplus \Der_k(A)$, the Lie bracket is expressed by
\begin{displaymath}
  [(a_1, \partial_1),(a_2, \partial_2)] = (\partial_1(a_2) -
  \partial_2(a_1) + \eta (\partial_1 \wedge \partial_2),
  [\partial_1 , \partial_2]).
\end{displaymath}
Conversely, any closed 2-form $\eta$ defines a structure
of Lie algebroid on $A\oplus \Der_k(A)$ whose curvature 2-form
is $\eta$. Denote by $\g^\eta$ $(\cong \g)$ the object of
$\LPic$ which corresponds to $\eta$.

Automorphisms of $\g $ are determined by $A$-linear homomorphisms $(a,
\partial)\mapsto (a + \omega(\partial), \partial)$ where $\omega \in
\Omega^{1,cl}(A)$ is a closed $1$-form, hence $\Aut (\g) \cong
\Omega^{1,cl}(A)$. For $\omega\in \Omega^1(A)$, the mapping
$(a,\partial)\mapsto (a+ \omega(\partial), \partial) $ defines an
isomorphism of Lie algebroids $\g^\eta \to \g^{\eta + d\omega } $
(identifying both sides with $A\oplus \Der_k(A)$ as $A$-modules).
There are no non-trivial automorphisms $\phi : \g \to \g$ such that
$\phi \circ \nabla= \nabla$, so therefore the group
$\Omega^{2,cl}(A)/d\Omega^1(A)$ classifies $A$-split affine Lie
algebroids.  Note that in the Zariski topology a closed 2-form need
not be locally the differential of a 1-form; here is the typical
example:
\begin{example} Let $X = \{(x,y) \in \Cb^2 : xy \neq 0 \}$, and
  $A=\Cb[x,y ,x^{-1},y^{-1}]$ be its affine $\Cb$-algebra, $\Oc_{X_h}$
  the sheaf of holomorphic functions in the Hausdorff topology on $X$
  and $\Oc_X$ the sheaf of regular rational functions in the Zariski
  topology.  The 2-form $\eta = 1/(xy) dx\wedge dy$ in $X$ is
  evidently closed.  Put $\g^\eta = A\oplus A\partial_x \oplus A\partial_y$
  and define the bracket $[\partial_x , \partial_y] = \eta(\partial_x
  \wedge \partial_y) = 1/(xy)$. The canonical projection defines a Lie
  homomorphism $\alpha : \g \to \Der_\Cb (A) = A\partial_x +
  B\partial_y$. Thus $\g$ is an affine Lie algebroid over $A$.  Making
  a cut in the $x$-plane in $\Cb^2$ we have $\eta^{an} = d (log (x)/y
  dy)$, so $\eta$ is locally exact in the Hausdorff topology. On the
  other hand, if $\phi $ is a rational 1-form the residue
  $\res_{y=0}(d\phi)=0$ (for fixed $x$), so $d \phi \neq \frac 1{xy}$.
  Therefore $\Oc_{X_h} \otimes_A \g$ is locally trivial whereas $\Oc_X
  \otimes_A\g$ is not locally trivial.
\end{example}

Define the truncated complex $\Omega^{\geq 1}(A) := \Omega^1(A)\to
\Omega^{2,cl}(A)$.  An $\Omega^{\geq 1}(A)$-torsor is a pair $(C,
\eta)$ where $C$ is an $\Omega^1(A)$-torsor and $\eta : C \to
\Omega^{2,cl}(A)$ is a mapping such that $\eta (c+ \phi) = \eta(c)+
d\phi$.
  
The set of connections $\Con(\g)$ on a Lie algebroid $(A, \g ,
\alpha)$ is an $\Omega^1(A)$-torsor, where the $\Omega^1(A)$-action is
defined by $\nabla \mapsto \nabla + \omega$, $\omega\in \Omega^1(A)$.
The mapping $\eta: \Con (\g) \to \Omega^{2,cl}(A)$, $\nabla \mapsto
R_\nabla$, satisfies $\eta(\nabla + \phi)= \eta(\nabla)+ d\phi$.
Therefore the pair $(\Con (\g), \eta)$ forms an $\Omega^{\geq
  1}(A)$-torsor, and we have a functor $\Cc : \LPic_A \to \Omega^{\geq
  1}(A)$-tors, $\g \mapsto (\Con (\g), \eta)$. The next lemma is taken
from \cite{bei-ber:jantzen}.
\begin{lemma}
  The functor $\Cc :\LPic_A \to \Omega^{\geq 1}$-tors is an
  equivalence of categories.
\end{lemma}

{\it Global structure}\/ : Set $\Omega^\bullet_X = \Omega^\bullet
(T_X, \Oc_X)$ (see~\ref{dr-complex}).  When $X$ is non-singular this
is just the ordinary de~Rham complex. Let $\Omega^{\geq 1}_X =
\Omega^1_X \to \Omega^{2,cl}_X$ be its truncation at degree $1$.  If
$\Omega \subset X$ is an open subset, the automorphism group $\Aut
(\g_\Omega) \cong \Gamma (\Omega , \Omega^{1,cl}_X)$, and if $\g_X,
\g'_X \in \LPic_X$ have connections $\nabla, \nabla'$ defined over
$\Omega$, then $\g_U \cong \g'_\Omega$ precisely when $R_\nabla -
R_{\nabla'} \in d \Omega_X^1(\Omega)$. It follows that the category
$\LPic_X$ is a sheaf forming a torsor over the sheaf of abelian groups
$\Omega_X^{\geq 1} = \Omega_X^1 \to \Omega_X^{2,cl}$, and that the
isomorphism classes of $\LPic_X$ are classified by the group $H^2(X,
\Omega_X^{\geq 1 })$.  \v{C}ech-representatives for vectors in $H^2(X,
\Omega_X^{\geq 1 })$ are, in evident notation, pairs $(\eta_i,
\phi_{ij})$ where $\eta_i$ are closed 2-forms defined in open sets
$U_i \subset X$ and $\phi_{ij}$ are 1-forms defined in $U_i\cap U_j$
such that $\eta_i-\eta_j = d\phi_{ij}$.

The connections on a locally $\Oc_X$-split Picard Lie algebroid $\g_X$
form a torsor over $\Omega^1_X$, so the sheaf of connections on $\g_X$
gives rise to a class in $H^1(X, \Omega^1_X)$ (the Atiyah class of
$\g_X$). A Picard Lie algebroid is {\it locally integrable} \/ if it is
locally isomorphic to $\Oc_X \oplus T_X$, i.e. \ if locally the
curvature 2-form of its connections are exact. These are the same as
$\Omega_X^{1, cl}$-torsors, hence they are classified by $H^1(X,
\Omega_X^{1, cl}) \subseteq H^2(X, \Omega_X^{\geq 1 })$.  If $X$ is a
complex manifold, then by Poincar\'e's lemma $\Omega^{\geq 1}_X$ is
quasi-isomorphic to $\Omega^{1,cl}$, so any Picard Lie algebroid is
locally integrable.  If $X$ is a compact complex algebraic manifold
then $H^2(X, \Omega_X^{\geq 1 }) = F^1H^2_{DR}$, the Hodge filtration
subspace, and the locally integrable Picard Lie algebroids correspond
to those classes that vanish on some Zariski open subset of $X$, i.e.,
precisely to $\Cb$-linear combinations of the algebraic cycle classes.

The category $\LPic_X$ is a $k$-vector space in categories, where the
$k$-vector space structure of $H^2(X, \Omega_X^{\geq 1 })$ is defined
directly on $\LPic_X$ using the Baer sum construction. Namely, if
$(\g^1_X, \alpha_1), (\g^2_X, \alpha_2) \in \LPic (X)$, $a , b \in k$,
then $ a\g^1_X = \Oc_X\oplus\g_X^1/\{(a\phi, - \phi); \phi \in \Oc_X=\Ker (\alpha_1)\}$ and the linear
combination $ a \g^1_X +b \g^2_X$ is a Picard algebroid $\g_X$ that is
equipped with a homomorphism of Lie algebroids $s_{a ,b} : \g^1_X
\times_{ T_X}\g^2_X \to \g_X$ such that $s_{ a,b} (\delta_1,\delta_2)= (a \delta_1 +b \delta_2)$.

{\it Relation with the Picard group}\/ : Letting $\Pic_X$ be the groupoid
of invertible sheaves on $X$ there exists a functor
\begin{eqnarray*}
  g :\Pic_X &\to& \LPic_X\\
\lambda &\mapsto& \cf_X(\lambda) \quad \quad (\text{see (\ref{linear})}).
\end{eqnarray*}
If $\phi : \lambda_1 \to \lambda_2$ is an isomorphism, then the
isomorphism $g(\phi) : \cf_X(\lambda_1)\to \cf_X(\lambda_2)$ is
defined by transport of structure, i.e. $g(\phi)(\delta_1)= \phi
\delta_1 \phi^{-1} \in \cf_X(\lambda_2)$.  Also, $\Pic_X$ is a group
in categories with composition $\otimes$ and $g$ is a homomorphism of
abelian groups in categories, expressed by
\begin{displaymath}
  \cf_X(\lambda_1 \otimes_{\Oc_X} \lambda_2)\cong \cf_X(\lambda_1) +  \cf_X(\lambda_2)\cong \lambda_1 \otimes_{\Oc_X}\cf_X(\lambda_2)\otimes_{\Oc_X}\lambda_1^{-1}.
\end{displaymath}
Here all isomorphisms are unique to unique automorphism of either
$\cf_X(\lambda_1 \otimes_{\Oc_X} \lambda_2)$, $\cf_X(\lambda_1) +
\cf_X(\lambda_2)$, or $ \lambda_1
\otimes_{\Oc_X}\cf_X(\lambda_2)\otimes_{\Oc_X}\lambda_1^{-1}$.  On the
level of isomorphism classes $[\lambda]$ and $[\g_X]$, parametrized by
\v{C}ech cohomology classes for an acyclic covering $X= \cup
U_\alpha$, $g$ defines a homomorphism of groups $dlog : H^{1}(X,
\Oc^*_X) \to H^{1}(X, \Omega^{1,cl})$, $[(\xi_{\alpha , \beta})]
\mapsto [(dlog(\xi_{\alpha , \beta }))]$, so that $dlog [\lambda] =
[\cf_X(\lambda)]$ (using obvious notation).

\begin{remark}\label{isopic} For any integer $n$ and $\g_X\in \LPic_X$ and $\lambda \in
  \Pic_X$ we get $\lambda^n\otimes_{\Oc_X}\g_X
  \otimes_{\Oc_X}\lambda^{-n} = \g_X + \cf_X(\lambda^n)\in \LPic_X$
  and there exists a push-out diagram
\begin{displaymath}
  \xymatrixnocompile{0 \ar[r]& \Oc_X\ar[r] \ar[d]^n & \lambda
    \otimes_{\Oc_X}\g_X\otimes_{\Oc_X} \lambda^{-1} \ar[r] \ar[d]^\phi&
  T_X \ar[r]\ar@{=}[d] & 0 & \\ 0 \ar[r]& \Oc_X\ar[r]  & \lambda^n\otimes_{\Oc_X}\g_X\otimes\lambda^{-n}\ar[r] & T_X\ar[r]&  0&  }
\end{displaymath}
where $\phi$ is an isomorphism in $\Lie_X$ (not in $\LPic_X$); see
above \Definition{defpic}.
\end{remark}

Assume that $X$ is non-singular of pure dimension $n$ and put
$\omega_X = \Omega^n_X \in \Pic_X$.  Since $\g_X$ is also a right
$\Oc_X$-module (\ref{linear}) we may consider the opposite \/ Lie
algebroid $\g^o_X$ of $\g_X\in \LPic_X$ (\ref{def-section-def}).  This
is again a Picard Lie algebroid and we have
\begin{displaymath}
  \g^o_X \cong  \cf_X(\omega_X) -\g_X  
\cong   \omega_X\otimes_{\Oc_X}  (-\g_X) \otimes_{\Oc_X} \omega^{-1}_X \cong - \omega^{-1}_X\otimes_{\Oc_X}\g_X\otimes_{\Oc_X}\omega_X.
\end{displaymath}
An isomorphism is given by
\begin{eqnarray*}
   \mu : \g^o_X &\to& -\omega^{-1}_X\otimes_{\Oc_X} \g_X
   \otimes_{\Oc_X} \omega_X \\
&=& \frac {\Oc_X \oplus (\omega^{-1}_X\otimes_{\Oc_X} \g_X
  \otimes_{\Oc_X} \omega_X)}{\{(\phi,\phi): \phi \in \Oc_X = \Ker
(\omega^{-1}_X\otimes_{\Oc_X} \g_X \otimes_{\Oc_X} \omega_X \to T_X)\}} \\
  \delta &\mapsto&  (0, -h \otimes \delta\otimes h^{-1} -  h\otimes
  \frac {\alpha(\delta)\cdot  h}{h}\otimes h^{-1}) \\
&=& (-  h\otimes  \frac {\alpha(\delta)\cdot  h}{h}\otimes h^{-1}, h \otimes  \delta\otimes h^{-1}),
\end{eqnarray*}
where $h \in \omega_X$, $h\neq 0$, and $\alpha (\delta)\in T_X$ acts
on $h$ by the negative of the Lie derivative, see (\ref{interchange}).
In particular, $\cf_X(\lambda)^o \cong
\omega_X\otimes_{\Oc_X}\cf_X(\lambda^{-1})\otimes_{\Oc_X}\omega^{-1}_X
\cong \cf_X(\omega_X\otimes_{\Oc_X}\lambda^{-1})$. Here one may either
regard $\omega_X \otimes_{\Oc_X}\lambda^{-1}$ as a left module over
$\cf_X(\lambda)^o$ or a right module over $\cf_X(\lambda)$; compare to
(\ref{interchange}). Let us check that the isomorphisms are consistent
with the fact that the operation $o$ is an involution:
\begin{eqnarray*}
  \g_X^{oo} &\cong& (\omega_X\otimes_{\Oc_X} -\g_X \otimes_{\Oc_X}
  \omega^{-1}_X)^o\\
& \cong& (-\omega^{-1}\otimes_{\Oc_X} \g_X \otimes_{\Oc_X}
\omega_X)^o\\ 
&\cong& \omega_X\otimes_{\Oc_X}(\omega^{-1}\otimes_{\Oc_X}
\g_X\otimes_{\Oc_X} \omega_X)\otimes_{\Oc_X}\omega_X^{-1} \cong \g_X,
\end{eqnarray*}
and more particularly $\cf(\lambda)_X^{oo} = \cf_X(\lambda^{-1}\otimes \omega_X)^o = \cf_X (\lambda 
\otimes \omega_X^{-1}\otimes \omega_X) = \cf_X(\lambda).$

Now let $X_h$ be a complex manifold. We then have 3 horizontal short exact
sequences, forming a commutative diagram where the vertical sequences
also form short exact sequences:
\begin{equation}\label{diagram}
  \xymatrix{
 & 0 \ar[d] &  2\pi i \Zb \ar @{^{(}->} [d]                & \Cb^*\ar @{^{(}->}[d]  & \\
0 \ar[r]    & 2\pi i\Zb \ar[d] \ar[r] &  \Oc_{X_h} \ar[r]^{\exp}\ar@{=}[d]     &
\Oc_{X_h}^* \ar[d]^{dlog} \ar[r]& 0  \\
0 \ar[r]& \Cb \ar[r]\ar[d]^{\exp} & \Oc_{X_h}  \ar[r]^{d}  \ar[d]^{c_1} &
\Omega^{1,cl}_{X_h} \ar@{=}[d] \ar[r] & 0 \\
0 \ar[r] & \Cb^* \ar[d]\ar[r] & \Oc_{X_h}^* \ar[d]\ar[r]^{dlog} & \Omega^{1,cl}_{X_h}
\ar[r]\ar[d] & 0  \\
& 0 & 0 & 0 & 
}
\end{equation}
inducing a commutative diagram where the  vertical sequences are exact
\begin{equation}\nonumber
  \xymatrix{
 &  H^0(X_h, \Oc^*_{X_h}) \ar[d]\ar[r]^{dlog} &  H^0(X_h, \Omega^{1,cl}_{X_h}) \ar[d] & \\
   &  H^1(X_h,2\pi i \Zb) \ar[d] \ar[r]^1 & H^1(X_h, \Cb^*)  \ar@{=}[r]\ar[d]     &
 H^1(X_h, \Cb^*)   \ar[d] &   \\
& H^1(X_h, \Oc_{X_h}) \ar[r]^{\exp}\ar[d] & H^1(X_h,
\Oc_{X_h}^*)  \ar[r]^{c_1}  \ar[d]^{dlog} & H^2({X_h}, 2\pi i \Zb) \ar[d]^f &  \\
 & H^1(X_h, \Oc^*_{X_h}) \ar[d]^{c_1}\ar[r]^{dlog} & H^1(X_h,\Omega^{1,cl}_{X_h} ) \ar[d]\ar[r]^{c} & H^2({X_h},\Cb)\ar[d] &   \\
& H^2(X_h, 2\pi i\Zb) \ar[r]^1 & H^2(X_h,\Cb^* )\ar@{=}[r] & H^2(X_h,\Cb^* ) & 
}
\end{equation}

The following proposition is due to A.  Weil \cite[Ch. V, no, 4, Lemme
2]{weil:kahler}.  Put $\bar H^2(X_h, 2\pi i \Zb) = \Imo (f)$.
\begin{proposition}\label{weil}($X_h$ is a complex manifold)
A Picard Lie algebroid  $\g_{X_h}$ is isomorphic to  $\cf_{X_h}(\lambda)$ for
some invertible sheaf  $\lambda$ if and only if  $c([\g_{X_h}])\in \bar H^2(X_h, 2\pi i \Zb)$.
\end{proposition}
In other words, $\Imo (dlog) = c^{-1}(\bar H^2(X_h,2\pi i \Zb))$.
\begin{pf} If $\{\xi_{\alpha
    \beta}\}$ is a \v{C}ech representative of $[\g_{X_h}] \in H^1(X_h,
  \Omega^{1,cl})$, then, since $\{\partial \xi_{\alpha \beta} \} \in
  \bar H^2(X_h, 2\pi i \Zb)$,  using the holomorphic Poincar\'e lemma
  (assuming the \v{C}ech covering is sufficiently fine), the chain
  $\{\phi_{\alpha \beta}\} =\{ \exp^{\int \xi_{\alpha \beta}}\}$ is a
  cocycle representing a class in $H^1(X_h, \Oc_{X_h}^*)$ such that $dlog
  (\phi_{\alpha \beta}) = \xi_{\alpha \beta}$. Choosing $\lambda\in
  \Pic_{X_h}$ such that $[\lambda]= \{\phi_{\alpha \beta}\}$ we have
  $\cf_{X_h}(\lambda)\cong \g_{X_h}$.
\end{pf}
\begin{remark}
  Poincar\'e's lemma holds more generally when $X_h$ is a locally
  holomorphically contractible complex analytic space \cite{reiffen},
  so \Proposition{weil} holds for such spaces.
\end{remark}

Assume that $X$ is a non-singular complex projective variety and $X_h$
its associated compact complex manifold. Put $\LPic^0_{X_h}= \Ker
(c)$; see the horizontal sequence (non-exact) above containing the
mapping $c$.  Consider the long exact sequence in homology of the
second horizontal short exact sequence in (\ref{diagram}) and let
$\bar H^1(X_h, \Cb)$ be the image of $H^1(X_h, \Cb)$ in $H^1(X_h,
\Oc_{X_h})\cong H^1(X, \Oc_X)$ (GAGA).  In analogy with the ordinary
Picard variety one may call
\begin{displaymath}
 \LPic^0_X \cong \frac {H^1(X, \Oc_X)}{\bar H^1(X_h, \Cb)} 
\end{displaymath}
the {\it Lie Picard variety} \/ of $X$. The functor $g$ above induces a
proper surjective morphism from the Picard variety of $X$ to its Lie
Picard variety, $\Pic^0_X \to \LPic_X^0$, which is simply the canonical mapping
\begin{displaymath}
 \frac {H^1(X, \Oc_X)}{H^1(X_h, 2\pi i \Zb)}\to \frac{H^1(X,
    \Oc_X)}{\bar H^1(X_h, \Cb)}.
\end{displaymath}
When $X_h$ is a curve one has ${\LPic}^0_{X_h}=0$, hence
$\LPic_{X_h}\subseteq H^2(X_h, \Cb)$; in particular,
$\cf_{X_h}(\lambda_1) \cong \cf_{X_h}(\lambda_2)$ precisely when
$\lambda_1, \lambda_2\in \Pic_{X_h}$ have the same degree. See also
(\ref{modpic}).

If $X_h$ is a Stein manifold it follows from the first horizontal
short exact sequence and Cartan's theorem B that $H^1(X_h,
\Oc_{X_h}^*)$ $=$ $ H^2(X_h, 2\pi i \Zb)$ and $H^2(X, \Oc_{X_h}^*)=
H^3(X_h, 2\pi i\Zb)$, so using the long exact sequence in cohomology
of the third horizontal short exact sequence we get an exact sequence
\begin{displaymath}
  H^1(X_h, \Cb^*) \to H^2(X_h, 2\pi i \Zb)\to H^1(X_h,
  \Omega_{X_h}^{1,cl} ) \to H^2(X_h,\Cb^*) \to H^3(X_h, 2\pi i \Zb).
\end{displaymath}

\section{Modules}\label{Modules}
\subsection{Basic facts}\label{basic-facts}
Let $(\Oc_X,\g_X, \alpha)$ be a Lie algebroid.
\begin{definition}\label{gmodule}
  A {\it left} $\g_X$-module $M$ is an $\Oc_X$-module $M$ and a
  homomorphism of $k$-Lie algebras
  \begin{displaymath}
    \rho : \g_X \to \cf_X(M) \quad (\text{see
    (\ref{linearalgebra})})
  \end{displaymath}
  satisfying $\rho (f \delta) \cdot m = f \rho (\delta)\cdot m$, and
  $\rho (\delta)\cdot fm = \alpha (\delta)(f)m + f \rho (\delta)\cdot
  m$, for all $f\in \Oc_X$, $\delta \in \g_X$, $m\in M$. When $M$ is
  faithful over $\Oc_X$ these conditions express that $\rho$ is a
  homomorphism of Lie algebroids.
  
  A {\it right} $\g_X$-module is defined by a $k$-linear homomorphism 
  \begin{displaymath}
    \nu :  \g_X\to \gl_k(M)
  \end{displaymath}
satisfying $\nu([\delta,  \eta])=-[\nu(\delta), \nu(\eta)]$, and $\nu (f \delta) \cdot m = f \nu
  (\delta)\cdot m - \alpha(\delta)(f)m = \nu(\delta)\cdot(fm)$, where  $\delta, \eta \in \g_X$.
  \end{definition}
  
  When we say module without specifying left or right we will always
  mean left module. A $\g_X$-module $M$ is the same as a
  $\Dc(\g_X)$-module and it is {\it coherent} \/ ({\it
    quasi-coherent}\/) if it is coherent (quasi-coherent) over the
  ring $\Dc(\g_X)$.  $M$ is {\it torsion free}\/ if for every non-zero coherent
  $\Oc_X$-submodule $N\subset M$ one has $\dim \supp N = \dim X$; this
  coincides with the usual notion of torsion free modules when $X$ is
  integral.  Note that if $M=\Dc(\g_X)M^0$ and $M^0$ is torsion free
  it does not follow that $M$ is torsion free, but if $M\in \coh
  (\g_X)$ there exists locally $M^0\in \coh (\Oc_X)$ such that $M=
  \Dc(\g_X)M^0$ so $\supp M = \supp M^0$, implying that $\supp M$ is
  an algebraic (analytic) subspace of $X$.
  
  We let $\Mod (\g_X)= \Mod (\Dc(\g_X))$ be the category of
  $\g_X$-modules and $\coh(\g_X) \subset \qcoh(\g_X)$ its
  sub-categories of coherent and quasi-coherent $\g_X$-modules.  The
  triangulated derived category of complexes of $\g_X$-modules
  $M^\bullet$ with coherent homology $H^\bullet(M^\bullet)\in \coh
  (\g_X)$ is denoted by $D_{\coh}(\g_X)$ and the subcategories of
  complexes that are bounded above and below are 
  $D^-_{\coh}(\g_X)$ and $D^+_{\coh}(\g_X)$. We shall mostly deal with
  the subcategory of bounded complexes $D^b_{\coh}(\g_X)=
  D^+_{\coh}(\g_X) \cap D^-_{\coh}(\g_X)$ mainly because it is
  generated by its homology objects, i.e., every object $M^\bullet \in
  D^b_{\coh}(\g_X)$ is isomorphic to a complex of coherent
  $\g_X$-modules; see the proof of \cite[I, Prop.  12.8]{borel:Dmod}.
  
  For completeness we quote results of Deligne and Bernstein [\S VI,
  Th. 2.10, Prop.  2.11, loc.\ cit] about complexes on a noetherian
  variety $X$. Let $D^b(\coh (\g_X))$ ($D^b(\qcoh (\g_X))$) be the
  category of bounded complexes of coherent (quasi-coherent)
  $\g_X$-modules, and $D^b_{\qcoh}(\g_X)$ the category of bounded
  complexes of $\g_X$-modules $M^\bullet$ such that
  $H^\bullet(M^\bullet)\in \qcoh (\g_X)$. Then the natural inclusion
  functors
\begin{displaymath}
D^b (\coh(\g_X)) \to D^b_{\coh}(\qcoh (\g_X))
\end{displaymath}
and
\begin{displaymath}
D^b (\qcoh (\g_X)) \to D^b_{\qcoh}(\Mod (\g_X))
\end{displaymath}
are equivalences of categories. Hence we have also a third equivalence
\begin{displaymath}
  D^b(\coh (\g_X))  \to D^b_{\coh}(\Mod (\g_X)).
\end{displaymath}

A {\it $\g_X$-connection} \/ on $M$ is a homomorphism of
$\Oc_X$-modules $\rho : \g_X \to \cf_X(M)$. When $\g_X=T_X$ this is
just a connection on the Lie algebroid $\cf_X(M)$ (assuming that $M$
is faithful over $\Oc_X$).  The $k$-linear homomorphism $ \nabla : M
\to Hom_{\Oc_X}(\g_X,M)$, $\nabla(m)(\delta)= \rho(\delta)(m)$, is the
first step in a sequence of $k$-linear homomorphisms of
$\Oc_X$-modules. Put $l= \rk \g_X$ and let $\{\delta_1, \ldots,
\delta_l\}$ be generators of the $\Oc_X$-module $\g_X$ ($l$ and the
generators are defined locally).  Put $\Omega^{p}(\g_X, M)=
Hom_{\Oc_X}(\wedge^p \g_X , M)$ and consider
  \begin{equation}\label{dr-complex}
     M \xrightarrow{\nabla} \Omega^{1}(\g_X, M)
     \xrightarrow{\nabla}  \cdots \xrightarrow{\nabla}
     \Omega^{l}(\g_X, M)
  \end{equation}
  where the mapping $\nabla: \Omega^p(\g_X,M) \to
  \Omega^{p+1}(\g_X,M)$ is defined in the usual way
  \begin{eqnarray}\nonumber
    &&(\nabla\omega)(\delta_{i_1}\wedge \dots \wedge
    \delta_{i_{p+1}}) =\\ \nonumber{} &&\sum_{1\leq s <t \leq
    p+1}(-1)^{s+t}\omega([\delta_{i_s}, \delta_{i_t}] \wedge
    \delta_{i_1} \wedge \cdots \hat {\delta}_{i_s} \cdots \hat
    {\delta}_{i_t} \cdots \delta_{i_{p+1}}) \\ \nonumber{} &&+
    \sum_{1\leq s\leq p+1} (-1)^{s+1}\delta_{i_s} \cdot \omega
    (\delta_{i_1} \wedge \cdots \hat{\delta}_{i_s} \cdots
    \wedge \delta_{i_{p+1}}).
  \end{eqnarray}
  If $\rho$ gives $M$ the structure of $\g_X$-module this is a
  complex, i.e.  $\nabla^2=0$, and we call $(\Omega^{\bullet}(\g_X,
  M), \nabla)$ the {\it de~Rham}\/ (-Chevalley-Hochschild) complex of
  the $\g_X$-module $M$; see also
  \cite[Lem. 2.3]{kallstrom:2smooth}.\typeout{Referens till II!!}
  
  The sheaf $Hom_{\Oc_X}(\g_X,\cf_X(M))$ forms a torsor over the sheaf
  of abelian groups $N_X:=Hom_{\Oc_X}(\g_X, \gl_{\Oc_X}(M))$. Hence
  supposing that local $\g_X$-con\-nections exist we get an
  obstruction class in $H^1(X, N_X)$ for the existence of a global
  $\g_X$-connection.

Some basic facts about the Lie algebroid $T_X =  Hom_{\Oc_X}(\Omega_{X/k}, \Oc_X)$ are summarized as follows:

\begin{proposition}\label{bfacts}  
  ($X$ is a space which in the algebraic case is assumed to be locally
  a sub-scheme of a regular scheme of finite type over $k$) The following are equivalent:
\begin{enumerate}
\item $X/k$ is non-singular;
\item $\Omega_{X/k}$ is locally free; hence $T_X$ is locally
free;
\item $T_X$ is a simple Lie algebroid, i.e.\ $T_X$ contains no proper
  coherent ideal;
\item $\Oc_X$ is a simple $T_X$-module, i.e. \  $\Oc_X$ contains no
proper coherent $T_X$-submodule.
  \end{enumerate}
\end{proposition}

\begin{pf} 
  $(1)\Leftrightarrow (2)$: See \cite{matsumura} ($k$ is
  algebraically closed of characteristic $0$).
   
  $(3)\Rightarrow (1)$, $(4)\Rightarrow (1)$: The assertion is local
  in $X$ so one may assume that $X \subset Y$ where, in the complex
  analytic case $Y= \Cb^n$, and in the algebraic case $Y= \Ab^n$. Let
  $j: X \to Y$ be the associated closed immersion of spaces and $I$ be
  the kernel of the surjection $\Oc_Y \to j_*(\Oc_X)$. Put $t= \dim X$
  (the maximal value of $\dim \Oc_x$ when $x\in X$). Then if $D_1,
  \dots , D_t \in T_Y$ and $a_1, \dots , a_t \in I$ one can form the
  elements $\det (D_i(a_j))\in \Oc_Y$ and we let $J$ be the ideal they
  generate.  We assert that the ideal $I^s:=j^*(J)\subset \Oc_X$ is
  stable under the action of $T_X$. Since $T_X= j^{-1}(T_Y(I)/IT_Y)$
  it suffices to see that $\delta \cdot (\det (D_i (a_j))) \in J$ when
  $\delta \in T_Y(I)$. This holds more generally, namely for any
  integer $l\geq 1$ letting $J_l$ be the ideal generated by
  elements$\det (D_i(a_j))$ where $D_1, \dots , D_l \in T_Y$ and $a_1,
  \dots , a_l \in I$, we have $T_Y(I)\cdot J_l \subset J_l$.  To see
  this consider a term of $\det (D_i(a_j))$, which has the form
   \begin{displaymath}
    \pm D_{1}(a_{i_1})D_2(a_{i_2})\cdots D_l(a_{i_l}),
   \end{displaymath}
   where $(i_1, \dots, i_l)$ is a permutation of $(1, \dots , l)$.
   Applying $\delta$ to such a term one gets $2l$ terms occuring
   either by replacing one of the $D_i$ by $D'_i=[\delta, D_i]$ or one
   of the $a_{i_s}$ by $a'_{i_s}:= \delta(a_{i_s})\in I$. This implies
   that $\delta \cdot \det (D_i(a_j))$ is a sum of $2l$ terms each
   belonging to $J_l$.
   
   By the Jacobian criterion of regularity for spaces of this type the
   singular locus of $X$ is the zero locus of $I^s$.  Assuming that
   $X$ is singular it follows that $I^s \subset \Oc_X$ is a proper
   ideal which is moreover a $T_X$-submodule. Also, $I^s T_X\subset
   T_X$ is an ideal of the Lie algebroid $T_X$, which is proper by
   Nakayama's lemma.  Hence $\Oc_X$ is not a simple $T_X$-module and
   $T_X$ is not a simple Lie algebroid.
  
  $(4)\Rightarrow (3)$: Let $\af \subset T_X$ be a proper coherent
  ideal of $T_X$.  It suffices to prove that $1\notin \af \cdot
  \Oc_X$, for then $\af \cdot \Oc_X \subset \Oc_X$ is a proper ideal,
  and since clearly $\af \cdot \Oc_X$ is a coherent $T_X$-module it
  follows that $\Oc_X$ is not simple.  Thus suppose that $\partial \in
  \af, \eta \in T_X$ and $\partial (a)= 1$, then $\eta = \partial
  (a)\eta= [\partial, a\eta]- a[\partial, \eta] \in \af$, implying
  that $T_X \subseteq \af$,  contradicting the assumption that
  $\af$ is a proper ideal.
  
  $(2) \Rightarrow (4)$: This is of course well-known.  Let $I \subset
  \Oc_X$ be a non-zero coherent $T_X$-submodule. It suffices to see
  that $I_x=\Oc_x$ when $x$ is a closed point in $X$.  By
  $(2)\Rightarrow (1)$ there exists a regular system of parameters
  $(x_1, \dots , x_n)$ of the regular local ring $\Oc_x$, so letting
  $\hat \Oc_x $ denote the completion of $\Oc_x$ by Artin-Rees lemma
  and Cohen's structure theorem one has $\Oc_x \subset \hat \Oc_x
  \cong k[[x_1, \dots , x_n]]$. (2) implies that $(\Omega_{X/k})_x$ is
  free over $dx_i$, hence there exist derivations $\partial_{x_i}\in
  \Der_k(\Oc_x)$ such that $\partial_{x_i}(x_j)= \delta_{ij}$.  Then
  considering a non-zero $f\in I_x \subset k[[x_1, \dots , x_n]]$, one checks
  that there exists a multi-index $\alpha$ such that $\partial^\alpha
  (f)\in I_x$ is invertible in $\Oc_x$; therefore $I_x = \Oc_x$.
\end{pf}

  \begin{remarks}
    \begin{enumerate}
      
    \item When $\alpha \neq 0$ the mapping $\ad : \g_X \to
      \cf_X(\g_X)$, $\ad (\partial)(\delta)= [\partial,\delta]$ is not
      $\Oc_X$-linear, and so does not define a module structure on
      $\g_X$.
      \typeout{Referens till Zariski-Lipman}
    \item Recall Zariski-Lipman's conjecture: If $T_X$ is locally
      free, then the equivalent conditions in \Proposition{bfacts}
      hold, and $\Omega_X\cong \Omega^1_X(T_X, \Oc_X )$
      (\ref{dr-complex}).  \cite{saito-kyoji:log} contains conditions
      ensuring that $T_X(I)$ be free when $I$ is the ideal of a
      divisor on a complex manifold.
      
    \item Let $\pi: X' \to X$ be the normalization morphism of an integral
      space. Let $J\subset \Oc_X$ be the conductor ideal, that is $J =
      \{\phi \in \Oc_X: \phi \pi_*(\Oc_{X'})\subset \Oc_X\}$.  By
      Seidenberg's result in \Remark{seidenberg} it easily follows that
      $J$ is $T_X$-submodule of $\Oc_X$. Therefore, by
      \Proposition{bfacts} we get a curious proof of the well-known
      fact that if $X$ is non-singular, then $X$ is normal.
  
\item A noetherian $k$-algebra $R$ of characteristic $p>0$ is a simple
  $\Der_k(R)$-module if and only if it has the form $R=k[t_1,\dots,
  t_n]/(t_1^p, \dots , t_n^p)$, where $k$ is a field of characteristic
  $p$ \cite{harper}.
\end{enumerate}
  \end{remarks}

\subsection{Interchanging left and right modules} \label{interchange}
We make a slight extension of the discussion in \cite[VI]{borel:Dmod}.
Let $M_1,M_2$ be left $\g_X$-modules and $N_1,N_2$ be right $\g_X$
modules \Defn{gmodule} defined by homomorphisms
\begin{eqnarray*}
  &l_i& : \g_X \to \cf_X(M_i)\\
&r_i& : \g_X \to \gl_k(N_i).
\end{eqnarray*}
Then $M_1\otimes_{\Oc_X} M_2$, $Hom_{\Oc_X}(M_1,M_2)$ and
$Hom_{\Oc_X}(N_1,N_2)$ are left $\g_X$-modules defined by the
homomorphisms $\rho_1, \rho_2, \rho_2$
\begin{eqnarray}\label{tensorprod}
 \rho_1 &:& \g_X \to \cf_X(M_1 \otimes_{\Oc_X} M_2)\\ \nonumber
\delta &\mapsto& l_1(\delta)\otimes \id_{M_1} + \id_{M_1} \otimes
l_2(\delta)\\
 \rho_2&:& \g_X \to \cf_X(Hom_{\Oc_X}(M_1, M_2))\\ \nonumber
\delta &\mapsto& \rho_2(\delta): f \mapsto l_i(\delta)\cdot  f - f \cdot l_2(\delta)\\
\label{r-l} \rho_3 &:& \g_X \to \cf_X(Hom_{\Oc_X}(N_1,N_2))\\ \nonumber
\delta &\mapsto& \rho_3(\delta): f \mapsto f\cdot l_1(\delta) - l_2(\delta)\cdot f.
\end{eqnarray}
The $\Oc_X$-modules $M_1\otimes N_1$ and $Hom_{\Oc_X}(M_1, N_1)$ are
right $\g_X$-modules defined by the homomorphisms $\nu_1 , \nu_2$
\begin{eqnarray}\label{lefttoright}
  \nu_1 &:& \g_X\to \gl_k(M_1\otimes N_1)\\\nonumber
\delta &\mapsto& -l_1(\delta)\otimes \id_{N_1} + \id_{M_1}\otimes
r_1(\delta),\\
\label{l-r} \nu_2 &:& \g_X \mapsto \gl_k(Hom(M_1, N_1))\\ \nonumber
\delta &\mapsto& \nu_2(\delta): f \mapsto f \cdot l_1(\delta) + r_1(\delta)\cdot f
\end{eqnarray}
\begin{remark}
  The $\Oc_X$-modules $Hom_{\Oc_X}(N_1, M_1)$ and
  $N_1\otimes_{\Oc_X}N_2$ are in general neither right nor left
  $\g_X$-modules \cite{oda}.
\end{remark}

Assume now that $X$ is non-singular of pure dimension $n$, so that
$T_X$ is a locally free $\Oc_X$-module of rank $n$. Then
\begin{displaymath}
  \omega_X := \Omega^n(T_X,\Oc_X)
\end{displaymath}
is an invertible module which moreover is a {\it right}\footnote{See
  \Section{picardliealg} for the case of Picard Lie algebroids.}
$\g_X$-module. The homomorphism $\nu : \g_X\to \omega_X$ is defined by
\begin{equation*}
\begin{aligned}
  (\phi\cdot \delta)(d) &= \phi (\delta \cdot d)- \alpha
  (\delta)(\phi(d)), \text{ where }\\ 
 \delta \cdot d &= [\alpha (\delta),\partial_1]\wedge \partial_2 \wedge \cdots \partial_n + \cdots +
\partial_1\wedge \cdots \wedge [\alpha (\delta) , \partial_n] \in \bigwedge^n T_X ,
\end{aligned}
\end{equation*}
$\partial_1, \dots , \partial_r \in \g_X$ form a local basis of
$\g_X$, and $d= \partial_1 \wedge \cdots \wedge \partial_n$.  Note
that $\Ker (\g_X \to T_X)$ acts trivially on $\omega_X$.\footnote{Of
  course, in the same way the $\Oc_X$-module $\Omega^r(\g_X, \Oc_X)=
  Hom_{\Oc_X}(\bigwedge^r\g{X}, \Oc_X)$, if $\rk \g_X =r $, is a right
  $\g_X$-module for which $\Ker (\g_X \to T_X)$ acts non-trivially.}
Then if $M$ is a left $\g_X$-module it follows that
\begin{displaymath}
  M_r=\omega_X\otimes _{\Oc_X}M
\end{displaymath}
is a right $\g_X$-module (\ref{l-r}), and starting with
a right $\g_X$-module $N$, 
\begin{displaymath}
  N_l= Hom_{\Oc_X}(\omega_X, N)= \omega^{-1}_X\otimes_{\Oc_X}N
\end{displaymath}
is a left $\g_X$-module (\ref{r-l}). Since $\omega_X$ is invertible
the pair $(\cdot \otimes_{\Oc_X} \omega_X , \omega^{-1}_X
\otimes_{\Oc_X}\cdot )$ defines an equivalence between the categories
of left and right $\g_X$-modules.

\subsection{Modules over Picard Lie algebroids}\label{modpic}
When $\g_X$ is a locally $\Oc_X$-split Picard Lie algebroid we {\it
  identify} $\Ker (\alpha)$ with $\Oc_X$, so in this case a
$\g_X$-module structure is locally defined by the (non-Lie) action of
$T_X$ with respect to a connection. But if $\rho : T_X \to \cf_X(M)$
is a $T_X$-module, then the composition $\rho\circ\alpha : \g_X \to
\cf_X(M)$ does not satisfy $\rho\circ\alpha (\phi)= i(\phi)$, where
$\phi \in \Oc_X \subseteq \g_X$ and $i: \Oc_X \to \gl_k(M)$. Hence
$T_X$-modules are considered as $\g_X$-modules only when $\g_X = \Oc_X
\oplus T_X$.

Picard Lie algebroids do not always have modules that are
coherent over $\Oc_X$. If  $r\in k$  we  let $\g^r_X$ denote
the push-out of $\g_X$ by the character $\Oc_X \to \Oc_X$,
$\phi \mapsto r\cdot \phi$ (see above~\ref{defpic}).
\begin{prop}\label{cohpic} ($X$ is a non-singular space)
  Let $\g_X$ be a Picard algebroid. Any $\g_X$-module which is
  coherent over $\Oc_X$ is locally free. The following are
  equivalent:
  \begin{enumerate}
  \item $\g_X$ has a non-zero module $M$ which is coherent over
$\Oc_X$, of some rank $r$.
  \item $\g^r_X = \cf_X (\lambda)$ for some invertible sheaf
$\lambda$, where $\lambda = \det M$, the determinant bundle of $M$.
  \end{enumerate}
 We have $ \cf_X (\det M) = \cf^{\tr}_X
  (M)$, the push-out of the linear algebroid $\cf_X (M)$ by the
  trace homomorphism $\tr : \gl_{\Oc_X}(M) \to \Oc_X$ (see
  \Remark{classical}).
\end{prop}
$(1)\Rightarrow (2)$ is proven in \cite[2.3.1]{bei-ber:jantzen};
the remaining parts are also well-known.
\begin{pf}  
  The first assertion follows since $\g_X$ is transitive and $X$ is
  non-singular; see, e.g.  \cite{borel:Dmod} for a corresponding
  statement for $\Dc$-modules. Thus if $M\in \coh (\Oc_X)$ is a
  $\g_X$-module, then $M$ is locally free, so $\lambda := \det M =
  \wedge^r_{\Oc_X} M$ is an invertible sheaf. Define a homomorphism of
  Lie algebroids $\g_X \to \cf (\lambda)$, by $\delta \cdot (m_1\wedge
  m_2 \wedge \cdots \wedge m_r) = (\delta\cdot m_1)\wedge m_2\wedge
  \cdots \wedge m_d + \cdots + m_1 \wedge \cdots \wedge (\delta \cdot
  m_r)$ (Leibniz' rule), for $\delta \in \g_X$ and $m_i\in M$; in
  particular $1_{\g_X}$ is mapped to $r 1_{\cf(\lambda)}$. This is an
  isomorphism. The last statement should now also be evident.
\end{pf}
\begin{remark}
  Let $\g_X\in \LPic_X$. Then $\omega_X$ is not a right $\g_X$-module
  (requiering that $1_{\g_X}= 1_{\Oc_X}$) unless $\g_X= \Oc_X \oplus
  T_X$. Moreover, the formula (\ref{l-r}) does not give
  $\omega_X\otimes M$ a structure of right module in this sense, but
  we will  later (\ref{directim}) regard it as a right
  $\g_X$-module  such that $1_{\g_X}$ acts by $-1$.
\end{remark}

\subsubsection{Twisted sheaves}\label{twisted} 
Let $X$ be a complex analytic manifold. The third horizontal short
exact sequence in the diagram before \Proposition{weil} gives a long
exact sequence in homology
\begin{displaymath}
\cdots \to H^1(X,\Cb^*) \to   H^1(X,\Oc_X^*) \xrightarrow{dlog} H^1(X, \Omega^{1,cl}_X)
\xrightarrow{t\ } H^2(X, \Cb_X^*)
  \to \cdots
\end{displaymath}
So if $\g_X\in \LPic_X$ and $ \phi := t([\g_X])\neq 1 \in H^2(X,
\Cb_X^*)$, then $\Mod(\g_X)$ contains no coherent $\Oc_X$-modules
\Prop{cohpic}.  One may then instead of the category of sheaves
consider the category of $\phi$-twisted sheaves and its sub-category
of twisted $\g_X$-modules $\Mod(\g_X, \phi)$. We will make some
remarks about the possibility to use twisted sheaves; see also
\cite[\S 3]{kashiwara:rep}. Objects $M\in \Mod(\g_X, \phi)$ are given
by a descent datum with respect to a covering of $X$, so that the
co-cycle condition for sheaves $[\partial \psi_{ij}] =
[\psi_{ik}\psi^{-1}_{jk}\psi^{-1}_{ij}] = 1 \in H^2(X, \Oc_X^*)$ (in
\v {C}ech-notation) is broken to $[\partial \psi_{ij}] = \phi \id $.
Morphisms of $\phi$-twisted sheaves are still global sections of
sheaves of morphisms, and $\Mod(\g_X, \phi)$ is an abelian category.

For example, let $c\in \Cb^*$ and $\lambda \in \Pic_X$, and suppose
that open sets $U_i$ are chosen so that $\lambda_{U_i}$ is trivial,
giving isomorphisms $\psi_{ij}: (\lambda_{U_j})_{U_i}\to
(\lambda_{U_i})_{U_j}$ where $\psi_{ij}\in \Oc^*_{U_i\cap U_j}$.
Choosing for each pair $(i,j)$ a $c$th root $\psi^c_{ij}$ of
$\psi_{ij}$, we get a 2-cocycle $\phi_{ijk}=
\psi^c_{ik}\psi^{-c}_{jk}\psi^{-c}_{ij}$ defining a class
$[\phi_{ijk}]\in H^2(X, \Cb^*_X)$.  The descent datum $\{f_{ij},
\Oc_{U_i}\}$ then defines an object $\lambda^c \in
\Mod(\cf_X(\lambda^c), \phi)$. Note that the linear algebroid
$\cf_X(\lambda^c) = c\cdot \cf_X(\lambda)$ is on ordinary sheaf.

Let $\{\psi_{ij}\}$ be a \v {C}ech representative of $[\g_X]\in
H^1(X,\Omega^{1,cl})$, choose $\int \psi_{ij}\in \Oc_{U_{ij}}$, using
the holomorphic Poincar\'e lemma, such that $\psi_{ij}= d\int
\psi_{ij}$, and put $f_{ij}=\exp^{\int \psi_{ij}}$, defining a 1-chain
with values in $\Oc_{X}^*$. So $\{f_{ij}\}$ is determined up to a
1-cocycle with values in $\Cb^*$, but let us make a {\it choice}\/ of
this 1-chain. The isomorphisms $f_{ij}:\Oc_{U_{ij}}\to \Oc_{U_{ji}}$
then gives a descent datum $\{f_{ij, \Oc_{U_{i}}}\}$, defining a
$\phi$-twisted sheaf of $\g_{X}$-modules when $\phi = t([\g_X])$. We
denote this $\phi$-twisted sheaf by $\lambda = \lambda
(\g_X)$\footnote{We make one choice of $\lambda$ for each $\g_X\in
  \LPic_X$.}.

Let $\g^1_X, \g^2_X\in \LPic_X$ and $M\in \Mod (\g^1_X, \phi_1), N\in \Mod(\g_X^2, \phi_2)$.

\medskip
{\it Tensor products.}  \/ We have $M\otimes_{\Oc_X}N \in \Mod(\g^1_X +
\g^2_X , \phi_1 \cdot \phi_2)$ where the sum $\g^1_X + \g^2_X$ is
taken in $\LPic_X$. When $N$ locally is invertible as
$\Oc_X$-module, then
\begin{eqnarray*}
 &M\to M\otimes_{\Oc_X}N \\
   &\Mod(\g_X^1,\phi_1)\to \Mod (\g_X^1+\g_X^2, \phi_1\cdot \phi_2)
\end{eqnarray*}
defines an equivalence of categories. This equivalence is not unique,
for one may change $N$ to $N\otimes_{\Cb}L$ where $L$ is any
$\phi_2$-twisted invertible $\Cb^*$-module.

In particular, if $\g^1_{X}= \g_{X}$ one can choose $\phi=
t([-\g_X])=t([\g_X])^{-1}$, $\g^2_X = -\g_X$, and $N= \lambda(\phi)$,
resulting in an equivalence
\begin{displaymath}
  \Mod(\g_X) \cong \Mod(T_X, \phi),
\end{displaymath}
between a category of sheaves and a category of $\phi$-twisted
sheaves.

\medskip
{\it Homomorphisms.} \/ We have $Hom_{\Oc_X}(M,N) \in \Mod (\g^2_X-
\g^1_X, \phi_2 \cdot \phi^{-1}_1)$.  Thus if $H^2(X, \Cb^*)$ has
torsion, there may exist ``untwisted'' morphisms between $M$ and $N$
even if $\phi_1 \neq \phi_2$. When $\phi_1=\phi_2=\phi$ then
$Hom_{\Oc_X}(M,N)$ is an ordinary sheaf of $T_X$-modules. Suppose that
$X= \cup X_i$ where $X_i \subset X$ are open subsets, and $M_i \in
\Mod (\g_{X_i}, \phi_i) $ are $\phi_i$-twisted sheaves of
$\g^i_{X_i}$-modules. Restricting $M_i $, $\g^i_{X_i}$, and $\phi_i$
to $X_{ij}=X_i\cap X_j$ results in a $\phi_{ij}$-twisted sheaf of
$\g^1_{X_{ij}}$-modules, where possibly $[\phi_{ij}]=1 \in
H^2(X,\Cb^*)$ so the $M_{ij}:=({M_i})_{X_j}$ may be ordinary sheaves.
Let $\{\psi_{ij}\}$ be a \v {C}ech-cocycle of isomorphisms of Lie
algebroids $\psi_{ij}: \g^i_{X_{ij}}\to \g^j_{X_{ji}}$. Thus the Lie
algebroids $\g^i_{X_i}$ on $X_i$ can be glued to a global Lie
algebroid $\g_X\in \LPic_X$ such that $\g^i_{X_i}\cong \g_{X_i}$. Now
consider $M_{ji}$ as a $\g^i_{X_{ij}}$-module using $\psi_{ij}$. Let
$f_{ij}: M_{ij}\to M_{ji}$ be isomorphisms of $\phi_{ij}$-twisted
sheaves of $\g^i_{X_{ij}}$-modules.  Then if $\{f_{ij}\}$ satisfies
the condition $[\partial f_{ij}]= \phi \in H^2(X, \Cb^*)$ there exists
$M\in \Mod (\g_X, \phi)$, unique to unique isomorphism, such that
$M_{X_i}\cong M_i$.

The above discussion shows that the category of $\phi$-twisted
sheaves on $X$ form a stack over the category of schemes/$k$ (or over
the category of complex analytic spaces)\cite{giraud}.

\subsection{Inverse images} 
\label{pullbackmod}
Let $\rho : \g_X \to \cf_X(M)$ be a $\g_X$-module and set $\g_Y =
\pi^+(\g_X)$ (\ref{pullbacklie}).  The pull-back functor $\pi^! :
\Mod (\g_X) \to \Mod (\g_Y)$ is defined by $\pi^!(M)= \Oc_Y
\otimes_{\pi^{-1}(\Oc_X)} \pi^{-1}(M)$ and the homomorphism
\begin{eqnarray}\nonumber
\pi^!(\rho) : \g_Y &\to& \cf_Y(\pi^*(M)),\\ \nonumber
\pi^!(\rho)(\phi\otimes \delta(x) , \partial_y) &=& \partial \otimes \id + \phi\otimes\rho (\delta),
\end{eqnarray}  
where $(\phi \otimes \delta(x), \partial_y)\in \g_Y=
\pi^*(\g_X)\times_{\pi^*(T_X)}T_Y$. We  check that $\pi^!(\rho)$ is
well-defined. Let $\psi \in \pi^{-1}(\Oc_X)$, $m\in \pi^{-1}(M)$, and  $\bar \psi$ be
the image of $\psi$ in $\Oc_Y$. Since $\partial_y (\bar \psi)= \phi
\bar {\alpha(\delta(x))(\psi)}$ one gets
\begin{eqnarray*}
  (\phi\otimes\delta (x),\partial_y)(\bar \psi \otimes m) &=&  \partial_y(\bar \psi )\otimes m + \bar \psi \phi\otimes \delta(x)\cdot m \\
&=&  \phi \bar {\alpha(\delta(x))(\psi)} \otimes m  + \bar \psi \phi\otimes \delta(x)\cdot m\\
&=& \phi \otimes \alpha (\delta(x))(\psi)m +  \bar \psi \phi\otimes
\delta(x)\cdot m \\ &=& \phi \otimes \delta(x)\cdot (\psi m) \\ 
&=& (\phi\otimes \delta (x), \partial_y)(1\otimes \psi m).
\end{eqnarray*}
Set 
\begin{displaymath}
  \Dc_{Y\to X}= \pi^!(\Dc(\g_X))= \Oc_Y\otimes_{\pi^{-1}(\Oc_X)}\pi^{-1}(\Dc(\g_X)).
\end{displaymath}
This is a
$(\g_Y, \pi^{-1}(\g_X))$-bimodule in an obvious way, and one
easily sees that $\pi^! ( \cdot )$ is isomorphic to $ \Dc_{Y\to
X}\otimes_{\pi^{-1}(\Dc(\g_X))}\pi^{-1}(\cdot)$. 

When $M\in \coh (\g_X)$, then $\pi^!(M)$ is a countable union of its
coherent submodules, and if
\begin{displaymath}
  \pi^*(\Imo (\alpha : \g_X \to T_X))\subseteq \Imo (d\pi : T_Y
  \to \pi^*(T_X)),
\end{displaymath}
then $\pi^! : \coh (\g_X) \to \coh (\g_Y)$. This is easily
checked by showing that the stalks $\pi^!(M)_y$, $y\in Y$, are
$\g_y$-modules of finite type. 
\begin{examples}\label{pull-examples}
  \begin{enumerate}
  \item Let $\pi : Y=\Ab^1 \to X=\Ab^2$ be the embedding of
affine spaces $ y \mapsto (0,y)$ and $M= \Dc (T_X) \in \coh
(T_X)$. Then $\pi^+(T_X)= T_Y$ and $\pi^!(M)= \Oc_Y
\otimes_{\pi^{-1}(\Oc_X)}\pi^{-1}(M)= \Dc_Y[ \partial_x]$,
(polynomials in $\partial_x$ with coefficients in the ring of
differential operators on $Y$).  Clearly $\pi^!(M)\notin \coh
(T_Y)$.
  \item Let $\pi : Y=\Ab^2 \to X= \Ab^1$ be the projection on
the first coordinate $(x,y) \mapsto x$. Then $\pi^+(T_X)= T_Y$
and $\pi^!(M)= \Oc_Y\otimes_{\pi^{-1}(\Oc_X)}\pi^{-1}(M) \in
\coh (T_Y)$.
\end{enumerate}
\end{examples}

Let $Y \xrightarrow{f \ } Z \xrightarrow{g \ } X$ be a composition of
morphisms. Choosing pull-backs $\g_Z = g^+(\g_X)$, $\g_Y=(g \circ
f)^+(\g_X)$ and $f^+(\g_Z)$, there is a canonical homomorphism
$c_{fg}: f^+(\g_Z)\to \g_Y$ which in general is not an isomorphism
(see \Lemma{prefib-iso}). Hence if $N$ is a $\g_Z$-module, its inverse
image $f^!(N)$ is in general not a $\g_Y$-module.  This will
nevertheless not cause problems for composed inverse images.
\begin{lemma} Let $M$ be a $\g_X$-module. There exists a
  canonical isomorphism of $\g_Y$-modules
\begin{displaymath}
  (g \circ f)^!(M) \cong f^! (g^!(M)).
\end{displaymath}

\end{lemma}
\begin{pf}
  Considered as $\Oc_Y$-modules the isomorphism is the canonical one
  \begin{displaymath}
    \Oc_Y\otimes_{(g\circ f)^{-1}(\Oc_X)}(g\circ f)^{-1}(M) \cong \Oc_Y
  \otimes_{f^{-1}(\Oc_Z)}f^{-1}(\Oc_Z \otimes_{g^{-1}(\Oc_X)}  g^{-1}(M)),
  \end{displaymath}
and this also defines the action of $\g_Y$ on the right
  side. One should verify that the $f^+(\g_Z)$-action on the right
  side is compatible with this $\g_Y$-action and the morphism $c_{fg}:
  f^+(\g_Z)\to \g_Y$, but we omit these straightforward details.
\end{pf}

We thus can and will always consider $f^!(g^!(M))$ as a $\g_Y$-module.

The derived version of $\pi^!$ is the functor
\begin{displaymath}
L\pi^! : D^-_{\coh}(\g_X) \to D^-_{\qcoh}(\g_Y), \quad M^\bullet \mapsto
  \Dc_{Y\to X }\otimes^{L}_{\pi^{-1}(\Dc(\g_X))}\pi^{-1}(M^\bullet)
\end{displaymath}
and if $X$ is non-singular $L\pi^!$ defines a functor
$D^b_{\coh}(\g_X)\to D^b_{\qcoh}(\g_Y)$.  One can compute
$L\pi^!(M^\bullet)$ either by taking a resolution of the $\g_X$-module
$M$ which is flat as $\Oc_X$-module or a resolution of $\Dc_{Y\to X}$
in the category of $(\g_Y, \pi^{-1}(\g_X))$-bimodules consisting of
modules which are flat as right $\pi^{-1}(\g_X)$-module.  Abusing the
notation slightly we again use $\pi^!$ for $L\pi^!$, so to define
$\pi^!$ completely one thus has to choose one resolution for each
$M^\bullet\in D^-(\g_X)$.  If $\g_X$ is locally free, then by the
Poincar\'e-Birkhoff-Witt theorem for Lie algebroids \Th{rinehart}
$\Dc(\g_X)$ is locally free over $\Oc_X$, so a resolution $F^\bullet
\to M^\bullet$ which is
flat as $\g_X$-module is also flat as $\Oc_X$-module; hence 
\begin{eqnarray*}
   \Dc_{Y\to X}\otimes^L_{\pi^{-1}(\Dc(\g_X))}\pi^{-1}(M^\bullet) &=&
   \Dc_{Y\to X}   \otimes_{\pi^{-1}(\Dc(\g_X))}\pi^{-1}(F^\bullet) \\
&=&   \Oc_Y\otimes_{\pi^{-1}(\Oc_X)}\pi^{-1}(F^\bullet)\\ &=&
 \Oc_Y   \otimes^L_{\pi^{-1}(\Oc_X)} M^\bullet,
\end{eqnarray*}
and $H^{-i}(\pi^!(M^\bullet))= Tor_{i}^{\pi^{-1}(\Oc_X)}(\Oc_Y,\pi^{-1}(M^\bullet))$.

\begin{prop}\label{composed} 
  ($\g_X$ is locally free over  $\Oc_X$) Let $Y \xrightarrow {f\ } Z \xrightarrow {g\ } X$ be a
  composition of morphisms. Then there exists a canonical isomorphism
  of functors $D^-(\g_X)\to D^-(\g_Y)$
  \begin{displaymath}
 (g\circ f)^!(\cdot) \to f^! \circ g^! (\cdot).
\end{displaymath}
\end{prop}
See  \cite[Ch. II,   Th. 2.3.21]{bjork:analD}  and \cite[VI, Prop
4.3]{borel:Dmod} for proofs when $\g_X=T_X$ and all spaces are
non-singular. The proof below is very similar.
\begin{pf}  By the above discussion, since $\g_X$ is
  locally free, there exists a resolution $F^\bullet$ of $M^\bullet\in
  D^{-}(\g_X)$ such that the terms $F^i$ are flat as $\Oc_X$-modules.
  Since flatness is preserved under base change and $f^!(\cdot)\cong
  \Oc_Y\otimes^L_{f^{-1}(\Oc_X)}(\cdot)$ (restricted to
  $D^{-}(\g_Z)\subset D^-(\Oc_Z)$) this gives:
\begin{eqnarray*}
  f^!(g^!(M^\bullet)) &=& \Dc_{Y\to  Z}\otimes^L_{f^{-1}(\Dc(\g_Z))}f^{-1}(\Dc_{Z\to
    X}\otimes^L_{g^{-1}(\Dc(\g_X))} g^{-1}(M^\bullet))\\ &=&
  \Oc_Y\otimes^L_{f^{-1}(\Oc_Z)}f^{-1}(\Oc_Z\otimes_{g^{-1}(\Oc_X)} g^{-1}(F^\bullet)) \\
&=& \Oc_Y\otimes_{f^{-1}(\Oc_Z)} f^{-1}(\Oc_Z\otimes_{g^{-1}(\Oc_X)}g^{-1}(F^\bullet)) \\
&=& \Oc_Y\otimes_{(g\circ f)^{-1}(\Oc_X)}(g\circ f)^{-1}(F^\bullet)\\
&=& \Dc_{Y \to X}\otimes^L_{(g\circ f)^{-1}(\Dc(\g_X))} (g\circ
f)^{-1}(F^\bullet)\\
&=& (g\circ f)^!(M^\bullet).
\end{eqnarray*}
\end{pf}

\subsection{Direct images}
\label{directim} 
{\it In this section all spaces are non-singular, and all morphism
  $\pi :Y\to X$ of spaces are such that the functor $\pi_*$ on the
  category of sheaves on $Y$ has finite cohomological dimension.}
\medskip

Let $\pi : Y\to X$ be a morphism of spaces and set $\g_Y =
\pi^+(\g_X)$. We will alway assume that $\g_Y$ is {\it locally free}.  Following Kashiwara's procedure for $T_X$-modules we
will define direct images of $\g_Y$-modules, using the $(\g_Y,
\pi^{-1}(\g_X))$-bimodule $\Dc_{Y\to X}$. The direct image of a
bounded complex of right $\g_Y$-modules $N^\bullet$ is defined by
\begin{displaymath}
  \pi_+^r(N)= R\pi_*(N^\bullet\otimes^L_{\Dc (\g_Y)}\Dc_{Y\to X}).
\end{displaymath}
To define the direct image of a left module one first notices  that by (\ref{r-l}) and (\ref{l-r}) 
\begin{displaymath}
  \Dc_{X \gets Y}= \omega_Y\otimes_{\Oc_Y}\pi^*(\omega^{-1}_X)\otimes_{\Oc_Y}\Dc_{Y\to X}.
\end{displaymath}
is a $(\pi^{-1}(\g_X), \g_Y)$-bimodule.  Then define $\pi_+ :
D^b(\g_Y)\to D^b(\g_X)$ by
\begin{displaymath}
  \pi_+(M^\bullet)=   R\pi_*(\Dc_{X\gets Y}\otimes^L_{\Dc (\g_Y)}M^\bullet)\in D^b(\g_X)
\end{displaymath}
when $M^\bullet \in D^b(\g_X)$. Since $\g_Y$ is locally free and $Y$
is non-singular it follows that the functor $\Dc_{X\leftarrow Y}\otimes^L_{\Dc(\g_Y)}$ has finite homological dimension (see for
instance \cite{kallstrom:gldim}), so $\pi_+(M^\bullet)$ does indeed
belong to  $D^b(\g_X)$.
Note that by the projection formula:
\begin{equation}
  \label{eq:proj-formula}
  \pi_+(M^\bullet) = \omega_X^{-1}\otimes_{\Oc_X} \pi_+^r(\omega_Y \otimes_{\Oc_Y}M^\bullet).
\end{equation}

We now discuss the basic functorial property of the direct image
functor. Consider a composition of morphims of spaces
\begin{displaymath}
  Y \xrightarrow{f\ }Z \xrightarrow{g\ } X
\end{displaymath}
and put $\g_Z = g^+(\g_X)$, and $\g_Y = f^+(\g_Z)$.  Then in general
\begin{displaymath}
 (g\circ f)_*(\Dc_{X\leftarrow Y}\otimes_{\Dc(\g_Y)}M) \neq
  g_*(\Dc_{X\leftarrow Z}\otimes_{\Dc(\g_Z)} f_*(\Dc_{Z\leftarrow Y}\otimes_{\Dc(\g_Y)}M)).
\end{displaymath}
So the non-derived direct image does not compose well, but we will see
that the situation is more satisfying when working with derived
categories. 
 \begin{prop}\label{direct-compose}
   ($\g_X$ and $\g_Z$ are locally free) As functors $D^b(\g_Y)\to
   D^b(\g_X)$ we have
  \begin{displaymath}
    (g\circ f)_+ (\cdot)  \cong  g_+ f_+(\cdot).
  \end{displaymath}
\end{prop}
The assumption that $\g_X$ and $\g_Z$ are locally free are satisfied
for instance if $\g_X$ is locally free and $g$ is submersive, or in
the situation of \Lemma{pullfree}.
\begin{pf}
  The proof is similar to  \cite[Ch. II Th.
  2.3.21]{bjork:analD}.

a) We first  prove 
  \begin{eqnarray}\label{trans-iso}
    \Dc_{Y\to X} &\cong& \Dc_{Y\to Z}\otimes_{f^{-1}(\Dc(\g_Z))}
    f^{-1}(\Dc_{Z\to X}) \\ \nonumber
&\cong& \Dc_{Y\to Z}\otimes^L_{f^{-1}(\Dc(\g_Z))} f^{-1}(\Dc_{Z\to X})
  \end{eqnarray}
  where the isomorphism holds in the category  of (complexes of) $(\g_Y,
  (g\circ f)^{-1}(\g_X))$-bimodules.  
  
  Begin with the first isomorphism. Note that the right side need not
  be a $\g_Y$-module when $c_{fg}: f^+(\g_Z)\to \g_Y$
  (\ref{pullbackmorphism}) is not an isomorphism, but both sides are
  $f^+(\g_Z)$-modules, so if we prove the isomorphisms in the category
  of $(f^+(\g_Z), (g\circ f)^{-1}(\g_X))$-modules, then since the
  $f^+(\g_Z)$-action in the left side is determined by the
  $\g_Y$-action (using the homomorphism $c_{fg}$), we may regard the
  isomorphisms in the category of $(\g_Y, (g\circ
  f)^{-1}(\g_X))$-modules. The morphism
  \begin{eqnarray*}
    & G : \Dc_{Y\to X} \to  \Dc_{Y\to Z}\otimes_{f^{-1}(\Dc(\g_Z))}
    f^{-1}(\Dc_{Z\to X})
\end{eqnarray*} 
defined by
\begin{multline*}
k(y) \otimes_{(g\circ f)^{-1}(\Oc)}(g\circ f)^{-1}(P(x)) \\
 \mapsto k(y)  \otimes_{f^{-1}(\Oc_Z)}f^{-1}(1_Z) \otimes_{f^{-1}(\Dc(\g_Z))}
f^{-1}(1_Z\otimes_{g^{-1}(\Oc_X)} g^{-1}(P(x)))
\end{multline*}
  where $k(y)\in\Oc_Y$ and $P(x)\in \Dc(\g_X)$, obviously is an
  isomorphism in the category of $(\Oc_Y, (g\circ
  f)^{-1}(\g_X))$-bimodules. It remains to see that it is
  $f^+(\g_Z)$-linear (see (\ref{pullbackmod})).  In the notation of (\ref{pullbackmorphism}),
  letting
  $r(y)=(h(y)\otimes_{f^{-1}(\Oc_Z)}f^{-1}(\phi(z)\otimes_{g^{-1}(\Oc_X)}g^{-1}(\delta
  (x)), \partial_z ), \partial_y)\in f^+(\g_Z)$, we have
  \begin{displaymath}
    c_{fg}(r(y))=  (h(y)\bar \phi(y)\otimes_{(g\circ f)^{-1}(\Oc_X)}(g\circ f)^{-1}(\delta(x)),\partial_y),
  \end{displaymath}
where $\bar \phi (y)$ is the image of $\phi(z)\in \Oc_Z$ in $\Oc_Y$.
Therefore
\begin{eqnarray*}
&&  G (r(y) \cdot  (k(y)\otimes_{(g\circ f)^{-1}(\Dc(\g_X))}(g\circ f)^{-1}(P)))\\
&=&  G (\partial_y(k)\otimes_{(g\circ f)^{-1}(\Dc(\g_X))}\otimes (g\circ f)^{-1}(P(x))\\
&& +  k(y)h(y)\bar \phi(y)\otimes_{(g\circ f)^{-1}(\Dc(\g_X))}\cdot(g\circ f)^{-1}(\delta(x)\cdot P) )\\
&=& \partial_y(k)\otimes_{f^{-1}(\Oc_Z)}f^{-1}(1_Z) \otimes_{f^{-1}(\Dc(\g_Z))}
f^{-1}(1_Z\otimes_{g^{-1}(\Oc_X)} g^{-1}(P(x))) \\
&&  + k(y)h(y) \bar \phi (y)\otimes_{f^{-1}(\Oc_Z)}f^{-1}(1_Z) \otimes_{f^{-1}(\Dc(\g_Z))}
f^{-1}(1_Z\otimes_{g^{-1}(\Oc_X)} g^{-1}(\delta(x)P(x)))\\
&=&  \partial_y(k)\otimes_{f^{-1}(\Oc_Z)}f^{-1}(1_Z) \otimes_{f^{-1}(\Dc(\g_Z))}
f^{-1}(1_Z\otimes_{g^{-1}(\Oc_X)} g^{-1}(P(x))) \\
&&  + k(y)h(y)\otimes_{f^{-1}(\Oc_Z)}f^{-1}(1_Z) \otimes_{f^{-1}(\Dc(\g_Z))}
f^{-1}(\phi(z)\otimes_{g^{-1}(\Oc_X)} g^{-1}(\delta(x)P(x)))\\
&=& r(y)\cdot k(y)\otimes_{f^{-1}(\Oc_Z)}f^{-1}(1_Z) \otimes_{f^{-1}(\Dc(\g_Z))}
f^{-1}(1_Z\otimes_{g^{-1}(\Oc_X)} g^{-1}(P(x)))\\
&=& r(y)\cdot G (k(y)\otimes_{(g\circ f)^{-1}(\Dc(\g_X))}(g\circ f)^{-1}(P)).
\end{eqnarray*}
This proves the first isomorphism in (\ref{trans-iso}). To see the
second isomorphism, first note that the $\Oc_Z$-module $\Dc(\g_Z)$ is
locally free since $\g_Z$ is locally free \Th{rinehart}, so flat $\g_Z$-modules are
flat as $\Oc_Z$-modules, and then since $\Dc(\g_X)$ is also locally
free it follows that $\Dc_{Z\to X}= g^*(\Dc(\g_X))$ is flat as
$\Oc_Z$-module. This gives the following isomorphisms in $D^b(\g_Y)$,
where $F^\bullet$ is a flat resolution of the $\g_Z$-module $\Dc_{Z\to
  X}$:
\begin{eqnarray*}
  \Dc_{Y\to Z}\otimes^L_{f^{-1}(\Dc(\g_Z))}f^{-1}(\Dc_{Z\to X})&\cong& 
  \Dc_{Y\to Z}\otimes_{f^{-1}(\Dc(\g_Z))} f^{-1}(F^\bullet) \\
&\cong& \Oc_Y\otimes_{f^{-1}(\Oc_Z)}f^{-1}(F^\bullet) \\
&\cong& \Oc_Y\otimes_{f^{-1}(\Oc_Z)}f^{-1}(\Dc_{Z\to X})\\
&\cong& \Dc_{Y\to Z}\otimes_{f^{-1}(\Dc(\g_Z))}f^{-1}(\Dc_{Z\to X}).
\end{eqnarray*}

b) Let first $N$ be a right $\g_Y$-module. Applying the projection
formula and (\ref{trans-iso}) one gets
\begin{eqnarray*}
  g^r_+ (f^r_+(N)) &=& Rg_*(Rf_*(N\otimes^L_{\Dc(\g_Y)} \Dc_{Y\to
    Z})\otimes^L_{\Dc(\g_Z)}\Dc_{Z\to X}) \\
&=& Rg_*(Rf_*(N\otimes^L_{\Dc(\g_Y)} \Dc_{Y\to
  Z}\otimes^L_{f^{-1}(\Dc(\g_Z))} f^{-1}(\Dc_{Z\to X}))) \\
&=& Rg_* \circ Rf_*(N\otimes^L_{\Dc(\g_Y)}\Dc_{Y\to X})\\
&=& R(g\circ f)_*(N\otimes^L_{\Dc(\g_Y)}\Dc_{Y\to X})= (g\circ f)^r_+(N),
\end{eqnarray*}
proving  the composition formula in the proposition for right
$\g_Y$-modules. Let now $M$ be a left $\g_Y$-module. By (\ref{eq:proj-formula}) 
\begin{eqnarray*}
  g_+(f_+(M)) &=& \omega_X^{-1}\otimes_{\Oc_X}
  g^r_+(\omega_Z \otimes_{\Oc_Z}\omega_Z^{-1}\otimes_{\Oc_Z}f^r_+(\omega_Y\otimes_{\Oc_Y}M)) \\
&=&
\omega_X^{-1}\otimes_{\Oc_X}g^r_+(f^r_+(\omega_Y\otimes_{\Oc_Y}M)))\\
&=& \omega_X^{-1}\otimes_{\Oc_X}(g\circ
f)^r_+(\omega_Y\otimes_{\Oc_Y}M))\\
&=& (g\circ f)_+(M).
\end{eqnarray*}
\end{pf}

Let $D^b_{\coh,l}(\g_Y) \subset D^b_{\coh}(\g_Y)$ be the subcategory
of bounded complexes $M^\bullet$ whose homology
$H^\bullet(M^\bullet)$, locally in $X$, is generated by an
$\Oc_Y$-coherent submodule. When $X$ is a noetherian space, then
$D^b_{\coh,l}(\g_X)= D^b_{\coh}(\g_X)$.
\begin{prop} \label{direct-coh}($\g_Y$ and $\g_X$ are locally free) Let $\pi : Y \to X$ be a proper
  morphism of spaces.
  We then have
  \begin{displaymath}
    \pi_+(D^b_{\coh,l}(\g_Y)) \subset D^b_{\coh,l}(\g_X).
  \end{displaymath}
\end{prop}
\begin{pf}
  The proof in \cite[Th. 2.8.1]{bjork:analD} works fine so we make
  only a few points for comparison.
  
1) The theorem holds  for right $\g_Y$-modules of
the form $K\otimes_{\Oc_Y} \Dc(\g_Y)$ when $K\in \coh (\Oc_Y)$: If
$F^\bullet \to K$ is a flat resolution of $K$, then $F^\bullet
\otimes_{\Oc_Y}\Dc(\g_Y)$ is a flat resolution of $K \otimes_{\Oc_Y}
\Dc(\g_Y)$, so
\begin{eqnarray*}
  \pi_+^r(K\otimes_{\Oc_Y}\Dc(\g_Y)) &=& R\pi_*(K\otimes_{\Oc_Y}\Dc(\g_Y)\otimes^L_{\Dc(\g_Y)} \Dc_{Y\to
    X} )\\
&\cong & R\pi_*(F^\bullet \otimes_{\pi^{-1}(\Oc_X)}\pi^{-1}(\Dc(\g_X))) \\
&\cong&  R\pi_*(F^\bullet)\otimes_{\Oc_X}\Dc(\g_X) \\
&\cong & R\pi_*(K)\otimes _{\Oc_X} \Dc(\g_X) ,
\end{eqnarray*}
where the third line follows from the projection formula since
$\Dc(\g_X)$ is locally free over $\Oc_X$. By properness, Grauert's or
Grothendieck's theorem implies that $R\pi_*(K)\in D^b_{\coh}(\Oc_X)$,
hence $\pi_+^r(K\otimes_{\Oc_Y}\Dc(\g_Y))\in D^b_{r,\coh}(\g_X)$
(complexes of right $\g_X$-modules).

2) Let $M$ be a left $\g_Y$-module and put $l=\rk \g_Y$.  Consider the
induced $(\g_Y, \g_Y$) -bimodule $M_{\ind}=M\otimes_{\Oc_Y} \Dc(\g_Y)$
(the left action is (\ref{tensorprod}) and the right action is
multiplication in $\Dc(\g_Y)$). Put $\Omega^i_Y(\g_Y)=
Hom_{\Oc_Y}(\bigwedge^i \g_Y, \Oc_Y)$; since $\g_Y$ is locally free
the de~Rham complex $\Omega^\bullet (M_{\ind}) =
\Omega^i(\g_Y)\otimes_{\Oc_Y}M \otimes_{\Oc_Y}\Dc(\g_Y)$
(\ref{dr-complex}). Then $\Omega^\bullet_Y(\g_Y)$ is an acyclic left
resolution of the right $\g_Y$-module $\Omega^l(\g_X, M) =
\Omega^l(\g_Y)\otimes_{\Oc_Y} M$ (\ref{l-r}). This is proven as
follows: First filter the complex $\Omega^\bullet (M_{\ind})$ by the
subcomplexes $F^\bullet_k$, where $F^i_k =
\Omega^i(\g_Y)\otimes_{\Oc_Y} M \otimes_{\Oc_Y}\Dc^{i+k}(\g_Y)$. Then
the associated graded complex $\gr F^\bullet \cong M\otimes_{\Oc_Y}
K^\bullet$, where $K^\bullet$ is the Koszul complex of the symmetric
algebra $\So(\g_Y)$ of the locally free $\Oc_Y$-module $\g_Y$ ($K^i=
\So(\g_Y)\otimes_{\Oc_Y}\bigwedge^{-i} \g_Y$), so it is
acyclic in degree $\leq - 1$. But this implies that $F^\bullet$ is
acyclic in degrees $\leq - 1$, and one sees that its homology in degree
$0$ is $\Omega^l(\g_X, M)$.
  
3) The module $M= \cup_{k\geq 0} M_k$ is filtered by $\Oc_Y$-coherent
submodules $M_k$ so that $\g_Y\cdot M_k \subset M_k$. Then for each
integer $k \geq 0$ the de~Rham complex $\Omega^\bullet (M_{\ind})$
contains a subcomplex $S_k^\bullet(M)$ (the Spencer complex of the
filtered module $M$) where
\begin{displaymath}
  S^i_k(M)= \Omega^i(\g_Y)\otimes_{\Oc_Y} M_{k+i}\otimes_{\Oc_Y}\Dc(\g_Y).
\end{displaymath}
One proves that locally in $Y$ when $k \gg 1$ then $S^\bullet_k(M)$ is
quasi-isomorphic to $\Omega^\bullet (M_{\ind})$; but as $\pi$ is
proper this quasi-isomorphism  holds locally in $X$ when $k\gg 1$, and since the assertion to
be proven is local in $X$ we may assume that  $k$ is choosen  so that, by 2),
\begin{displaymath}
  S^\bullet_k(M)\cong \Omega^l(\g_Y)\otimes_{\Oc_Y}M
\end{displaymath}
in $D^b_r(\g_Y)$.

4) If a right $\g_Y$-module $N$ has a bounded acyclic resolution
\begin{displaymath}
  F^l \to F^{l-1}\to \cdots \to F^0 \to N
\end{displaymath}
such that $\pi_+^r(F^i)\in D^b_{r,\coh}(\g_X)$, then $\pi_+^r(N)\in
D^b_{r,\coh}(\g_X)$. This is proven by a straighforward  induction in $l$.  From this
it follows by 1) and 3) that
$\pi_+^r(\Omega^l(\g_Y)\otimes_{\Oc_Y}M)\in D^b_{r,\coh}(\g_X)$. Then
$M_1 := \omega_Y\otimes_{\Oc_Y}\Omega^{l}(\g_Y)^{-1}\otimes_{\Oc_Y}M$
is a coherent left $\g_Y$-module and (see (\ref{interchange}))
\begin{eqnarray*}
  \pi_+(M) &=& \omega_X^{-1}\otimes_{\Oc_X} \pi_+^r(\omega_Y  \otimes_{\Oc_Y}M) \\
&=&  \omega_X^{-1}\otimes_{\Oc_X} \pi_+^r(\Omega^l(\g_Y)\otimes_{\Oc_Y}M_1) \in D^b_{\coh}(\g_X).
\end{eqnarray*}
\end{pf}
Let $M^\bullet[d]$ be the translation of a complex to the left $d$
steps, and $\dop_{Y,X}= \dim Y - \dim X$ be the relative dimension of
a morphism of pure dimensional spaces $\pi : Y \to X$.

\begin{thm}\label{adj-prop} ($\g_Y$ and $\g_X$ are locally free)
  Let $\pi : Y \to X$ be a proper morphism of non-singular spaces.
  Then 
\begin{displaymath}
  R\pi_*  RHom_{\g_Y}(M^\bullet, \pi^!(N^\bullet)[\dop_{Y,X}]) =
  RHom_{\g_X}(\pi_+(M^\bullet), N^\bullet)
\end{displaymath}
when $M^\bullet \in D^b_{\coh}(\g_Y)$ and $N^\bullet \in D^b_{\coh}(\g_X)$.
\end{thm}
This adjointness property is proven in a similar way as for
$\Dc$-modules \cite[VIII, Prop. 9.10]{borel:Dmod}, so we omit the
quite long proof. Let us only remark that it suffices to construct the
trace morphism $\tr : \pi_+ \pi^!(M^\bullet[d_{Y,X}]) \to M^\bullet $
and for this it suffices to consider closed embeddings and projections
separately. To see this, factorize $\pi = p\circ i$, where $i: Y \to
Y\times X$ is the graph embedding and $p: Y\times X \to X$ is the
projection on $X$. Then if $\tr_i :i_+ i^!(N^\bullet)[-\dim X]\to
N^\bullet$ and $\tr_p :p_+p^!(M^\bullet)[\dim Y] \to M^\bullet $ are
defined for any $N^\bullet \in D^b_{\coh}(\g_{Y\times X})$, $M^\bullet
\in D^b_{\coh}(\g_X)$, we first get $\tr_i:i_+i^! p^!(M^\bullet) \to
p^!(M^\bullet)$, and applying $p_+$
\begin{displaymath}
  \pi_+\pi^!(M^\bullet) = p_+i_+ i^!p^!(M^\bullet)  \xrightarrow{p_+(\tr_i)} p_+p^!(M^\bullet)\xrightarrow{\tr_p}M^\bullet,
\end{displaymath}
where the equality is because $\pi^! = i^!p^! $ and $ \pi_+= p_+i_+ $
\Props{direct-compose}{composed}.  Note that in the situation of
\Theorem{kashiwara} below $\tr_i$ is an isomorphism.

Letting $M^\bullet \in D^b(\g_X)$ and $Y \subset X$ be a closed algebraic
(analytic) subset one has a distinguished triangle in $D^b(\Oc_X)$
(\ref{dist-triangle})
\begin{equation}\label{localcoh}
   R\Gamma_{[Y]}(M^\bullet) \to M^\bullet \to M^\bullet(*S)\to.
\end{equation}
Since $\Dc(\g_X)$ is locally free over $\Oc_X$, injective
$\g_X$-modules are injective as $\Oc_X$-modules, so an injective
resolution of $M^\bullet$ is acyclic for the functor
$\Gamma_{[Y]}(\cdot)$, hence (\ref{localcoh}) is a distinguished
triangle in $D^b (\g_X)$.

Let now $\pi : Y \to X$ be a closed immersion of non-singular spaces
and $D^b_Y(\g_X)\subset D^b(\g_X)$ be the subcategory of complexes
$M^\bullet $ such that $\supp H^\bullet (M^\bullet)\subseteq i(Y)$.
Next result is basically due to Kashiwara.
\begin{thm}\label{kashiwara} 
  Let $Y\to X$ be a closed immersion of non-singular spaces and $\g_X$
  a locally free and transitive Lie algebroid.  Then we have a
  distinguished triangle
  \begin{equation}\label{eq:kash} 
    \pi_+\pi^!(M^\bullet)[\dop_{Y,X}] \to M^\bullet  \to  M^\bullet (*Y) \to 
  \end{equation}
  The functor $\pi_+$ defines an equivalence of categories $D^b(\g_Y)
  \to D^b_Y(\g_X)$ with quasi-inverse $\pi^!(\cdot)[\dop_{Y,X}]$;
  $\pi_+$ preserves coherence and the restriction of
  $\pi^!(\cdot)[\dop_{Y,X}]$ to $D^b_Y(\g_X)$ also preserves
  coherence.
\end{thm}
Because of \Theorems{adj-prop}{kashiwara} one may find it
convenient to include the translation in the definition of the inverse
image functor, redefining it to $M\mapsto
\pi^!(M^\bullet)[\dop_{Y,X}]$.\footnote{Which  is the convention in
  [loc. \ cit.].}

\begin{pf}
  When $\g_X=T_X$ the proof can be found in \cite{bjork:analD} or
  \cite[VI, Th.  7.13]{borel:Dmod} so we only make a sketch to see how 
  to get (\ref{eq:kash}), following the
  proof in [loc.\  cit.]. 
  
  a) Put $d=-d_{Y,X}$ and let $\{f_1, \ldots ,f_d\}$ be a regular
  sequence locally generating the ideal $I$ of $Y$. Let $\partial_1,
  \ldots \partial_d$ be a dual basis of the $k$-vector space
  $\sum_{i=1}k^d df_i$. Locally we have
  $\Oc_Y  \otimes_{\pi^{-1}(\Oc_X)}\pi^{-1}(\Oc_X df_1 \wedge \cdots \wedge
  df_d)= \omega^{-1}_Y\otimes_{\pi^{-1}(\Oc_X)}\omega_X$.
  
  b) Using the Koszul complex based on the regular sequence $\{f_1,
  \ldots ,f_d\}$ to resolve the $\pi^{-1}(\Oc_X)$-module $\Oc_Y$ one
  gets, by a),
\begin{displaymath}
  H^{0}(\pi^!(M)[d_{Y,X}])= M^{I}\otimes_{\Oc_Y} \omega^{-1}_Y
  \otimes_{\pi^{-1}(\Oc_X)} \pi^{-1}(\omega_{X}),
\end{displaymath}
where $M^I=\{m\in M : I\cdot m=0\}\subset M$.

c) There  are isomorphisms of $\g_Y$-modules 
\begin{displaymath}
  \Dc_{Y\to X} = \pi^*(\Dc(\g_X)) \cong \Dc(\g_X)/I\Dc(\g_X) \cong
  \Dc(\g_Y)\otimes_k\Cb[\partial_1, \ldots , \partial_d].
\end{displaymath} 
In particular  $\Dc_{Y\to X}$ is locally free over $\g_Y$, and so
$\pi_+$ is exact.\  

d) Since $\pi_*$ is exact,  c) implies
\begin{eqnarray*}
  && \pi_+H^0(\pi^!(M)[d_{Y,X}]) \\
    &=&    \pi_*(\pi^{-1}(\omega^{-1}_X)\otimes_{\pi^{-1}(\Oc_X)}\Dc_{Y\to X}
    \otimes_{\Oc_Y}\omega_Y \otimes_{\Dc(\g_Y)} H^{0}(\pi^!(M)[d_{Y,X}]))
\\
&=& \omega_X^{-1}\otimes_{\Oc_X}\pi_* (\frac {\Dc(\g_X)}{I\Dc(\g_X)}\otimes_{\pi^{-1}(\Oc_X)}\omega_Y \otimes_{\Dc(\g_Y)} (M^I \otimes_{\Oc_Y}\omega^{-1}_Y\otimes_{\pi^{-1}(\Oc_X)}\pi^{-1}(\omega_X))).
\end{eqnarray*}
Then the mapping $\Dc(\g_X)\times M^I \to M$, $(P, m)\mapsto P\cdot M$ 
induces a   homomorphism of $\g_X$-modules
\begin{displaymath}
  \mu : \pi_+H^0(\pi^!(M)[d_{Y,X}]) \to \Gamma_{[Y]}(M).
\end{displaymath}
Using the transitivity of $\g_X$ one proves that $\mu$ is an
isomorphism in the same way as  [Prop. 7.10, loc.\  cit.].

e) By d) there exists an isomorphism of functors 
\begin{displaymath}
  \Gamma_{[Y]}(\cdot) \cong \pi_+ H^0(M^!(\cdot)[d_{Y,X}])
\end{displaymath}
on the category $\Mod(\g_X)$.  Since   $\Dc(\g_X)$ is locally free over
$\Oc_X$,  injective $\g_X$-modules are injective as
$\Oc_X$-modules; hence, as  $\pi_+$ is exact, injective resolutions
are acyclic both for $\Gamma_{[Y]}(\cdot)$  and $\pi_+
H^0(\pi^!(\cdot)) = \pi_+(Hom_{\Oc_X}(\Oc_X/I, \cdot))$. Therefore
$\mu$ induces an isomorphism of functors
\begin{displaymath}
  R\Gamma_{[Y]}(\cdot )\cong \pi_+\pi^!(\cdot)[d_{Y,X}]
\end{displaymath}
on $D^b(\g_X)$.  Hence the distinguished triangle (\ref{localcoh}) can
be identified with (\ref{eq:kash}).
\end{pf}

\section{Smooth modules}
\label{sec-smooth}
\subsection{Elementary properties}
Let $(\Oc_X , \g_X, \alpha)$ be a Lie algebroid on a space $X$. We are
interested in $\g_X$-modules which are locally of finite type with
respect to some distinguished sub-algebroid $\g'_X$ of $\g_X$.  One
example to keep in mind is $ \g_X(I) \subseteq \g_X$ for some ideal
$I$ of $\Oc_X$ (see~\ref{linear}).

\begin{definition}\label{smooth-def} 
  Let $M\in \coh(\g_X)$. Then $M$ is {\it smooth} \/ along $\g'_X$
  at a point $x\in X$ if there exists a $\g'_{x}$-submodule
  $M^0_x \subseteq M_x$ of finite type over $\Oc_{x}$, which
  generates the $\g_{x}$-module $M_x$.  Let $\Omega\subset X$
  be an open subset.  $M$ is point-wise smooth in $\Omega$ if
  it is smooth at all points $x$ in $\Omega$, and globally
  smooth along $\g'_X$ in $\Omega$ if $M_\Omega$ contains a
  coherent $\Oc_\Omega$-submodule $M_\Omega^0$ which generates
  the $\g_\Omega$-module $M_\Omega$ and satisfies $\g'_\Omega
  \cdot M^0_\Omega\subseteq M^0_\Omega$.
\end{definition}
It often follows from the context which sub-algebroid $\g_X'\subseteq
\g_X$ is intended, and we then abbreviate by simply saying that $M$ is
(point-wise) smooth. We then also let  $\coh_s(\g_X)\subset \coh
(\g_X)$   denote  the sub-category of smooth modules.

Obviously, a globally smooth module is point-wise smooth.  When $X$ is
a quasi-compact variety, so $X$ is a noetherian scheme, the converse
also holds \Prop{globalsmooth}, so there is no distinction between
globally smooth and point-wise smooth modules, and we say only smooth.
This is not true when $X$ is a quasi-compact complex manifold;
see~\ref{gaga-section}.

Notice that $\g'_x M^0_x\subseteq M^0_x$ implies that $M$
satisfies the following weaker condition, which occasionally is
useful.
\begin{definition}\label{along}  Let $\delta $ be a vector in
  $\g_X(U)$, where $U$ is an open subset of $X$. We say that $M$ is
  point-wise smooth along $\delta$ in $\g_U'$ when for any point $x
  \in U$ there exists an integer $n (=n(x))$ such that
  \begin{displaymath}
    \delta_x^n M^0_x \subseteq \Dc^{n-1}_X(\g'_x)M^0_x,
  \end{displaymath}
  where $M^0_x \subset M_x$ is a submodule of finite type over
  $\Oc_x$ such that $\Dc (\g_x)M^0_x = M_x$.
\end{definition}

The first part of the following lemma implies that the above
definitions do not depend on the choice of generating $\Oc_X$-module
$M^0 \subseteq M$.

\begin{lemma}\label{all-coh} Let $(A, \g , \alpha)$ be a Lie algebroid
  over a $k$-algebra $A$, $\bfr \subseteq \afr \subseteq \g$ be
  sub-algebroids, and $M\in \coh (\g)$. Then:
  \begin{enumerate}
  \item $M$ is smooth along $\afr$ if and only if every
$A$-submodule of finite type $M^1 \subseteq M$ generates an
$\afr$-module which remains of finite type over $A$;
  \item If $M$ is smooth along $\afr$, then $M$ is smooth along
$\bfr$.
  \end{enumerate}
\end{lemma}
\begin{pf}
  Trivial.
\end{pf}

\begin{prop}\label{extension}  Let
  \begin{displaymath}
    0\to M_1 \to M \to M_2 \to 0
  \end{displaymath}
  be a short exact sequence of coherent $\g_X$-modules. Then
  $M\in \coh_s(\g_X)$  if and only if $M_1, M_2\in \coh_s(\g_X)$.
\end{prop}
\begin{pf}
  It is evident that $M\in \coh_s(\g_X)$ implies $M_1, M_2 \in
  \coh_s(\g_X)$. To prove the converse, let $x$ be a point in $X$ and
  $M_x^0$ be an $\Oc_x$-submodule of $M_x$ of finite type, and
  put $N= \Dc(\g'_x)M_x^0$. Since $M_2 \in \coh_s(\g_X)$  the image $\bar
  N \subset (M_2)_x$ is of finite type over $\Oc_x$, so $N$ is of
  finite type if  $N_1:=N \cap (M_1)_x$ is of finite type. We have
  $N_1\subset N$; $N$   is of finite type over $\g'_x$ and
  $\Dc(\g'_x)$ is noetherian; hence $N_1$ is of finite type over
  $\g'_x$. Since $N_1 \subset (M_1)_x$ and $M_1\in \coh_s(\g_X) $ it
  follows that $N_1$ is of finite type over $\Oc_x$ \Lem{all-coh}.
\end{pf}

\begin{prop}\label{point} Let $\g'_X$  be a Lie sub-algebroid of
  $\g_X$ and $M\in \coh (\g_X)$.
  \begin{enumerate}
  \item ($X$ is a space) The subset
$\Omega=\{x\in X: M_x \text{ is smooth}\}$ is open.  If $M$ is
smooth at a point $x$, then there exists a neighbourhood
$\Omega$ of $x$ such that $M$ is globally smooth over $\Omega$.
\item ($X$ is a quasi-compact space) Let $S \subset X$ be a divisor
  whose sheaf of ideals is $I$, and assume that $M(*S)\in
  \coh(\Oc_X(*S))$. Then there exists an integer $n$ such that $M$ is
  point-wise smooth along $I^n\g_X$.
  \end{enumerate}
\end{prop}
By $(1)$ it follows that if $M\in \coh_s(\g_X)$ and $M_0\subset M$ is
a coherent $\Oc_X$-submodule, then $\Dc(\g'_X)M_0 \in \coh (\Oc_X)$.
It also implies that $M$ is point-wise smooth on a variety if it is
smooth at all closed points.

\begin{remark} We leave out the
  proofs of the following assertions: (i) Let $x\in S$ and $e(M,x)$ be
  the smallest integer $n$ such that $M_x$ is smooth along
  $I^{n+1}_x\g_x$.  Then for a short exact sequence of $\g_X$-modules
\begin{displaymath}
  0 \to M_1 \to M \to M_2 \to 0
\end{displaymath}
one has
\begin{displaymath}
  e(M,x) = \max \{e(M_1, x), e(M_2, x)\}.
\end{displaymath}
(ii) Let $k$ be a field of characteristic $0$ and $R$ a $k$-algebra
which is a complete discrete valuation ring with residue field $k=
R/\mf$, and $K$ be the fraction field of $R$. Let $M$ be an
$l$-dimensional $K$-vector space which is a module over the Lie
algebroid $T_K:=\Der_k(K)$; $T_K$ contains the subalgebroids
$T_R=\Der_k(R)$ and $\g:= \{\partial \in T_R: \partial (\mf)\subset
\mf\}$. N. Katz \cite{katz:galois} has defined certain non-negative
rational numbers $\lambda_1, \lambda_2, \dots , \lambda_l$ (the slopes
of $M$) and an irregularity index $\Irr (M)= \sum_{i=1}^n \lambda_i$;
put $\lambda_{max}= \max \{\lambda_1, \dots , \lambda_l\}$.  Define
$e(M)=e(M,0)$ as in $(i)$ with $\g_x= \g$.  Then
\begin{displaymath}
  e(M) = r \lambda_{max},
\end{displaymath}
where $r$ is the multiplicity of $\lambda_{\max}$. Hence $\Irr (M)=0$
if and only if $e(M)=0$.
\end{remark}

\begin{pf}
  (1): Assume that $M_x$ is smooth and let $M^0_x= \sum_{i=1}^l
  \Oc_x \bar m_i$ be an $\Oc_x$-module of finite type which
  generates the $\g_x$-module $M_x$, satisfying $\g'_x \cdot
  M^0_x\subseteq M^0_x$.  The vectors $\bar m_i$ are germs of
  sections $m_i$ which are defined in some open set $\Omega$;
  set $M_\Omega^0 = \sum \Oc_\Omega m_i$, defining a coherent
  submodule of $M_\Omega$.  The $\Oc_X$-module $\g'_X$ is
  locally of finite type, so choose sections $\delta_1 , \dots
  , \delta_n$ defined in an open subset $\Omega' \subset
  \Omega$, $x\in \Omega'$, such that $\g'_{\Omega'}= \sum
  \Oc_{\Omega'} \delta_j$. After further shrinking $\Omega'$ we
  can assume that $\delta_j m_i \subset M^0_{\Omega'}$, for
  every $i,j$; hence $M^0_{\Omega'}$ is a
  $\g'_{\Omega'}$-module which is coherent over $\Oc_{\Omega'}$
  and $M^0_x = (M^0_{\Omega'})_x$ generates the $\g_x$-module
  $M_x$.  Since $M$ is coherent over $\g_X$ there exists an
  open subset $\Omega_1 \subset \Omega'$, $x\in \Omega_1$, such
  that $(M^0_{\Omega'})_{\Omega_1}= M^0_{\Omega_1}$ generates
  the $\g_{\Omega_1}$-module $M_{\Omega_1}$, and therefore
  $M_y^0$ generates $M_y$ for all $y \in \Omega_1$.
  
  (2): By quasi-compactness the assertion follows from
  $(1)$ if we show that there exists an integer $n$ such that
  $M_x$ is smooth along $I_x^n \g$ for a point $x\in X$.  Let
  $j: M_x \to M(*S)_x$ be the natural mapping and put $M_1 =
  \Ker (j)$, $M_2= \Imo (j)$. It suffices to prove that $M_1$
  and $M_2$ are smooth \Prop{extension}.  Let $m_1, \ldots ,
  m_r $ be generators of the $\g_x$-module $M_1$ and $\delta_1,
  \cdots , \delta_s$ generators of the $\Oc_x$-module $\g_x$.
  Now there exists an integer $n$ such that $I^n_x\delta_j
  m_i=0$, for all $i,j$, which implies that $M_1$ is smooth
  along $I^n_x\g_x$.
  
  Let instead the $m_i$ denote vectors that generate the
  $\Oc(*T)_x$-module $M(*S)_x$, such that $M_2 \subseteq M_3:=
  \sum \Dc(\g_x) m_i$; by \Proposition{extension} it suffices
  to see that $M_3$ is smooth.  But since $M_3 \subseteq
  M(*S)_x$, and $M(*S)_x$ is of finite type over $\Oc(*T)_x$,
  it is evident that there exists an integer $n$ such that
  $I_x^n \delta_i m_j \subseteq \sum \Oc_x m_l \subseteq M_3$,
  for all $i,j$, implying that $M_3$ is smooth along $I_x^n
  \g_x$.
\end{pf}

\begin{prop}\label{globalsmooth} ($X$ is a quasi-compact variety)
  Let $M$ be a coherent $\g_X$-module.  If $M$ is point-wise
  smooth along a sub-algebroid $\g'_X \subseteq \g_X$, then it
  is globally smooth along $\g'_X$.
\end{prop}
\begin{pf}  Since $X$ is quasi-compact and $M$ is coherent over $\g_X$
  it follows from $(1)$ in \Proposition{point} that there exist
  a finite number of open sets $\Omega_i$ such that $\cup
  \Omega_i = X$ and coherent $\Oc_{\Omega_i}$-submodules
  $M^0_{\Omega_i} \subseteq M_{\Omega}$, satisfying
  $\g'_{\Omega_i}\cdot M^0_{\Omega_i} \subseteq M^0_{\Omega_i}$
  and $\Dc (\g'_{\Omega_i})\cdot M^0_{\Omega_i} =
  M_{\Omega_i}$.  By \cite[Prop.  2]{borel-serre:riemann-roch}
  there exist coherent $\Oc_X$-sub-modules $M^i \subseteq M$
  such that $M^i_{\Omega_i} = M^0_{\Omega_i}$. Now put $M^{00}=
  \sum \Dc (\g'_X)\cdot M^i$. Since each $\Oc_X$-module $\Dc
  (\g'_X)\cdot M^i$ is coherent \Lem{all-coh} and the sum is
  finite it follows that $M^{00}\subseteq M$ is a globally
  defined coherent $\Oc_X$-submodule such that $\g'_X \cdot
  M^{00} \subseteq M^{00}$ and $\Dc(\g_X)\cdot M^{00} = M$.
\end{pf}

\begin{prop}\label{maximal}
  Let $M\in \coh (\g_X)$. There exists a maximal Lie
  sub-algebroid $\g^M_X \subset \g_X $ such that $M$ is smooth
  along $\g^M_X$.
\end{prop}
\begin{pf}
  Let $\g^M_X\subseteq \g_X$ be the Lie subalgebroid that is generated
  by the sum of all sub-algebroids of $\g_X$ along which $M$ is
  smooth; it is coherent since it is locally generated by its sections
  and $\Oc_X$ is noetherian. By \Theorem{gabber} $M$ is smooth along
  $\g_X^M$.
\end{pf}

Let $M\in \coh (\g_X)$ and  $I$ be the coherent ideal whose zero locus is $\supp M$.
\begin{proposition}($X$ is a space) \label{propsing}
  \begin{displaymath}
I \g_X \subseteq \g_X^M \subseteq \g_X(I) \quad \text{(see~\ref{maximal})}
\end{displaymath}
\end{proposition}
\begin{pf} $I\g_X \subseteq \g_X^M$: If $x\notin Z$ there is nothing to prove  so let $x\in
  Z$, $\phi\in I_x$, $\nabla = \phi\delta \in I_x\g_x$, and
  $m_x \in M_x$.  One proves by induction that  for any positive integer $n$
  one can rewrite $\nabla^n\in \Dc (I_x\g_x)$ as
  \begin{displaymath}
    \nabla^n= \delta^n \phi^n + a_{n-1}\nabla^{n-1}+ \cdots +
    a_0,\quad {\text{ where} }\quad a_i \in \Oc_x.
  \end{displaymath}
  Therefore, if $\phi^n m_x=0 $, then $\nabla^n \cdot m=
  (a_{n-1}\nabla^{n-1}+ \cdots + a_0)m$, whence the assertion.
  
$\g_X^M \subseteq \g_X(I)$:  The assertion being local in $X$ we can assume that
  all sections below are globally defined. 

If $\delta \in \g_X^M$ it suffices to prove that there exists a
Zariski open dense subset $Z_0 \subseteq Z$ such that $\delta_x \in
\g_x(I)$ when $x\in Z_0$, for then if $\phi \in I$,  $\alpha
(\delta)(\phi)\in I( Z_0\cup (X\setminus Z)) = I(X)$ \Lem{hartog}.

Supposing that $\delta \notin \g(I)$ we prove that there exists a
Zariski open dense subset $Z_0 \subset Z$ such that $\delta_x \notin
\g^M_x$ when $x\in Z_0$: There exists a function $\phi \in I$ such
that $\alpha (\delta)(\phi)\notin I$.  Set $Z_0= \{x\in Z :
\alpha(\delta)(\phi)_x \notin \mf_x\}$; this is a Zariski open dense
subset of $Z$ since $\alpha (\delta)(\phi)_x \notin I_x$.  Let $x\in
Z_0$ and let $M_x^{\phi_x}$ the $\phi_x$-invariants of
$M_x$. Clearly $M_x^{\phi_x}\neq 0$ so we get a non-zero $\g_x$-module
  \begin{displaymath}
     M_1 = \Dc (\g_x)M_x^{\phi_x}.
  \end{displaymath}
  Let $\g^1_x$ be the Lie algebroid that is generated by
  $\g^M_x$ and $\delta_x$.  By \Proposition{extension} it
  suffices to prove that $M_1$ is not smooth along the
  sub-algebroid $\g^1_x\subset \g_x$.  As $\Dc (\g_x)$ is
  noetherian it follows that $M_1\subset M$ is of finite type,
  hence there exists a sub-module $M^0_1 \subseteq
  M_x^{\phi_x}$ of finite type over $\Oc_x$ such that
  $\Dc(\g_x)\cdot M^0_1 = M_1$.  Clearly it suffices to prove
  that $\delta^n_x \cdot M^0_1 \nsubseteq
  \Dc^{n-1}(\g^1_x)\cdot M^0_1$. Since $ \phi_x \delta \cdot
  M^0_1 = - \alpha (\delta)(\phi_x)\cdot M^0_1 \neq 0$,
  ($\alpha (\phi_x)\notin \mf_x$), so $\delta \cdot M^0_1
  \nsubseteq M_1^0$, this takes care of the $n=1$ case. Assume
  that $n\geq2$, and suppose the contrary that
\begin{displaymath}
    \delta^n\cdot M^0_1 \subseteq \Dc^{n-1}(\g^1_x)\cdot M^0_1.
  \end{displaymath}
  Then
\begin{displaymath}
  \phi \cdot \delta^n \cdot M^0_1 \subseteq \phi \cdot
  \Dc^{n-1}(\g^1_x)\cdot M^0_1 \subseteq \Dc^{n-2}(\g^1_x)\cdot
  M^0_1,
\end{displaymath}
and
\begin{displaymath}
\phi \cdot \delta^n\cdot M^0_1 = [\phi, \delta^n]\cdot M^0_1 =
\delta \cdot [\phi, \delta^{n-1}] \cdot M^0_1 - \alpha
(\delta)(\phi)\cdot \delta^{n-1}\cdot M^0_1;
\end{displaymath}
hence
\begin{displaymath}
  \alpha (\delta)(\phi) \cdot \delta^{n-1}\cdot M^0_1\subseteq
  \delta \cdot [\phi , \delta^{n-1}]\cdot M^0_1 +
  \Dc^{n-2}(\g^1_x)\cdot M^0_1 = \Dc ^{n-2}(\g^1_x)\cdot M^0_1.
\end{displaymath}
Since $\alpha (\delta)(\phi)$ is invertible an induction gives a contradiction.
\end{pf}

\subsection{Smooth complexes}
\label{preserve}
Throughout this subsection we let $\pi : Y \to X$ be a morphism of
non-singular spaces, and if  $\g'_X \subset \g_X$  is an inclusion of Lie
algebroids we let $\g'_Y$ be the image of $\pi^+(\g'_X)$ in $\g_Y:=
\pi^+(\g_X)$.

Let $D^b_s(\g_X)$ be the sub-category of bounded complexes $M^\bullet$
whose homology $H^\bullet(M^\bullet)\in \coh_s(\g_X)$.
\begin{lemma}
  Let 
  \begin{displaymath}
    M^\bullet_1 \to M^\bullet \to M^\bullet_2 \to
  \end{displaymath}
be a distinguished triangle in $D^b(\g_X)$. If two vertices of this
triangle belong to
$D^b_s(\g_X)$, then the third vertex also belongs to $D^b_s(\g_X)$.
\end{lemma}
\begin{pf}
  This follows immediately from \Proposition{extension}, using the long 
  exact sequence in homology.
\end{pf}

\subsubsection{Inverse images}
When $M^\bullet \in D^b_{\coh}(\g_X)$ then any submodule
of finite type of the homology module
$H^\bullet(\pi^!(M^\bullet))$ is coherent. Let us therefore
agree that a $\g_Y$-module is smooth if all its coherent
submodules are smooth.

The category $\coh_s(\g_x)$,  $x\in X$, does in
general not contain sufficiently many projectives (or flats), so the
next result is not quite immediate. 
\begin{prop}\label{preservation}
  ($\g_X$ is locally free  and transitive) If $M^\bullet \in D^b_s(\g_X)$, then the homology of
  $\pi^!(M^\bullet)$ is smooth along $\g'_Y$.
\end{prop}

\begin{pf} We will be brief since standard arguments for
  $\Dc$-modules are used. It suffices to prove the assertion
  for a generating class in $D^b_s(\g_X)$, so if $M$ is a
  smooth $\g_X$-module and we prove that the complex $\pi^!(M)$
  has smooth homology we are done.
   
  Factorize $\pi$ as $\pi = p \circ i$ where $i : Y \to
  Z=Y\times X$ is the graph morphism and $p: Z\to X$ is the
  projection on the second factor. As $\g_X$ is locally free
  and $p$ is  a smooth morphism of non-singular spaces we
  have $\pi^!(M) \cong i^!  p^!(M)$ \Prop{composed} and $\g'_Y=
  i^+(p^+(\g'_X))$ \Lem{prefib-iso}, so it suffices to prove
  the assertion for $i$ and $p$ separately.
   
  $p$: Since $p$ is flat and $\g_X$ is locally free, $p^!(M)$
  is the single degree complex 
  \begin{displaymath}
    \Dc_{Y\to    X}\otimes_{p^{-1}(\Dc(\g_X))}p^{-1}(M)=\Oc_Y\otimes_{p^{-1}(\Oc_X)}p^{-1}(M);
  \end{displaymath}
  see the proof of \Proposition{composed}.  This module clearly
  is smooth along $\g'_Z= p^+(\g'_X)$ (see (\ref{pullbackmod})).
   
  $i$: Letting $N\in \coh_s(\g_Z)$ we have to prove that
  $i^!(N)\in D^b_{s}(\g_Y)$. First note that $\g_Z= p^!(\g_X)$
  is locally free since $p$ is submersive and $\g_X$ is locally
  free. Now the closed embedding $i$ can be factorized into a sequence 
  of   closed embeddings of the type $\pi: Y\to Z$ where $\codim_Z
  \pi(Y) =1$, and the pull-back of a locally free transitive Lie
  algebroid on a non-singular space to a non-singular subspace
  again is locally free and transitive \Lem{pullfree}, so by
  \Proposition{composed} we may assume that $\codim_Z i(Y) =1$.
  Let $I_Y$ be the ideal of $Y$.  Then since $\g_Z = p^!(\g_X)$ is
  locally free and $Y$ is non-singular one can use the Koszul
  complex to compute $i^!(M)$; see above
  \Proposition{composed}.  The only non-zero cohomology
  $\g_Y$-modules are
\begin{displaymath}
  M/I_YM  \quad \text {and } \quad M^{I_Y}.
\end{displaymath}
That any coherent $\Oc_Y$-submodule of these modules generates a
$\g'_Y$-module $(i_*(\g'_Y) = \g'_X(I_Y)/I_Y\g'_X)$ which
remains coherent over $\Oc_Y$ follows immediately since any coherent
$\Oc_X$-submodule of $M$ generates a $\g'_X$-module that remains
coherent over $\Oc_X$.
\end{pf}

We will give a partial converse of \Proposition{preservation}.
Let $y \in Y$ and $x= \pi(y)\in X$. A morphism of Lie
algebroids $\pi ' : \g_Y \to \g_X$ is submersive at $y$ if the
induced morphism $\g_Y \to k_y \otimes_{\Oc_y} \pi^*(\g)_y$ is
surjective; this is, by Nakayama's lemma, the same as to say
that the morphism $\g_y \to \pi^*(\g)_y$ is surjective.
\begin{prop}\label{inv-equiv} ($\g_X$ 
  is locally free) Let $M^\bullet \in D^b_{\coh}(\g_X)$.  Assume that
  $\pi$ is flat on $\supp M^\bullet \subset X$, and that for each
  point $x\in \supp M^\bullet$ there exists a point $y\in \pi^{-1}(x)$
  at which $\pi' :\g'_Y \to \g'_X$ is submersive. Then the following
  are equivalent:
\begin{enumerate} 
\item $M^\bullet \in  D^b_s(\g_X)$;
\item $\pi^!(M^\bullet )\in D^b_s(\g_Y)$.
  \end{enumerate}
\end{prop}

\begin{remark} 
  When considering completely regular complexes one can improve
  \Proposition{inv-equiv}, not demanding that $\pi$ be flat
  \Prop{inv-equiv-complete}.
\end{remark}

\begin{pf} 
  $(1)\Rightarrow (2)$: Apply \Proposition{preservation}.  $(2)
  \Rightarrow (1)$: It suffices to prove this for a generating class
  of $D^b(\g_X)$, so letting $M\in \coh (\g_X)$ satisfy $\pi^!(M)\in
  D^b_s(\g_Y)$, we need to prove that $M\in \coh_s(\g_X)$.
  
  Let $x\in \supp M$ and put $A= \Oc_x$. If $M^0_x\subset M_x$ an
  $A$-submodule of finite type, we have to check that
  $\Dc(\g'_x)M^0_x$ remains of finite type over $A$. There exists a 
  point $y\in
  Y$ such  that $x= \pi (y)$ and the canonical mapping $ \g'_y \to
  B\otimes_A\g'_x$ is surjective, where we have put $B=\Oc_y$. Then
\begin{eqnarray*}
 \pi^!(M)_y &\cong& \pi^*(\Dc_x)\otimes^L_{\Dc(\g_x)} M_x \\
&\cong& B \otimes^L_{A} M_x \quad \quad  \text{($\g_x$ is    free over $\Oc_x$)}\\
&\cong&  B\otimes_{A}M_x  \quad \quad \quad \text{($A\to B$ is   flat)}.
\end{eqnarray*}
As the restriction $\g'_y \to \g'_x$ is surjective, so is its
extension to a mapping $\Dc(\g'_y) \to B \otimes_A\Dc (\g'_x)$; hence,
by flatness,
\begin{displaymath}
B\otimes_A \Dc (\g'_x)M^0_x \subset \Dc(\g'_y)B\otimes_AM^0_x.
\end{displaymath}
By $(2)$ the right side is of finite type over $B$, and since $B$ is
noetherian the left side is also of finite type. Then by faithful
flatness $\Dc (\g'_x)M^0_x$ is of finite type over $A$.
\end{pf}

\subsubsection{Direct images}
\label{directsmooth} As in \Proposition{direct-coh} we consider the
category $D^b_{sl}(\g_Y)$ of complexes $M^\bullet$ whose homology
$H^\bullet(M^\bullet) \in \coh_{s}(\g_Y)$ locally in $X$ is generated
by an $\Oc_Y$-coherent submodule. When $X$ is a quasi-compact variety
we have $D^b_{sl}(\g_Y)= D^b_s(\g_Y)$.

Put $T_{Y/X}= \Ker (T_Y \to \pi^*(T_X))$ and $\g_{Y/X}=
\alpha^{-1}(T_{Y/X})$. For simplicity we shall {\it assume} \/ that
$\bfr_Y := \Ker (\g_Y \to T_Y)\subset \g'_Y$.
\begin{thm}\label{direct-smooth} 
($\g_X$ is locally free and transitive) 
Let  $\g'_Y \subset \g_Y$ and $\g'_X \subset \g_X $ be sub-algebroids
and  $\afr_Y$  be the Lie algebroid that is generated by $\g'_Y +
\g_{Y/X}\subseteq \g_Y$.  Assume that $\pi$ is proper and 
that 
\begin{equation}\label{log-cond}
 \Imo (\pi^*( \g'_X) \to \pi^*(\g_X)) \subseteq  \Imo (\afr_Y \to
 \pi^*(\g_X))
\end{equation}
(canonical mappings). Then $\pi_+$ defines a functor
  \begin{displaymath}
D^b_{sl}(\g_Y) \to D^b_{s}(\g_X).
\end{displaymath}
\end{thm}
The condition (\ref{log-cond})  means that, locally in $Y$
and $X$, sections of $\g'_X$ should be liftable to sections of
$\af_Y$.

\begin{pf} Letting $M^\bullet \in D^b_{s}(\g_Y)$ we have to prove that
  \begin{displaymath}
    R\pi_*(\Dc_{X\leftarrow Y}\otimes^L_{\Dc(\g_Y)}M^\bullet) \in D^b_s(\g_X).
  \end{displaymath}
  a) Let $\Dc^n_{X\leftarrow Y}$ be the right $\Oc_Y$-module that is
  generated by the image of $\pi^{-1}(\Dc^n(\g_X))$ in the
  $(\pi^{-1}(\g_X), \g_Y)$-bimodule $\Dc_{X\leftarrow Y}$. As
  $\Dc^n_{X\leftarrow Y} \g_{Y/X} \subseteq \Dc^n_{X\leftarrow Y}
  \bfr_Y$ and $\bfr_Y\subseteq \g'_Y$, (\ref{log-cond}) implies
  \begin{equation}\label{base-eq}
    \pi^{-1}(\g'_X)\Dc^n_{X\leftarrow Y} \subseteq \Dc^n_{X\leftarrow
      Y} \afr_Y \subseteq \Dc^n_{X\leftarrow
      Y} \g'_Y.
  \end{equation}
  
  b) It suffices to prove the assertion for a generating class for the
  considered category, and since the assertion is local in $X$, it
  suffices to prove that $\pi_+(M)\in D^b_{s}(\g_X)$ when $M\in
  \coh_s(\g_Y)$ is generated by a coherent $\Oc_Y$-submodule $M^0$.
  
  The morphism $\pi$ can be factorized $\pi= p \circ i$ where $p:
  Z=Y\times X \to X$ is the projection on the second factor and $i: Y
  \to Z$ is the graph morphism, which is a closed embedding. We
  contend that it suffices to prove the theorem for $i$ and $p$
  separately. Since $\pi^! = i^! p^!$ \Prop{composed} we then have to
  check that the conditions in the theorem are satisfied for $i$ and
  $p$: Put $\g'_Z=p^+(\g'_X)= p^*(\g'_X)$ and $\g_Z= p^+(\g_X)=
  p^*(\g_X)$. Since $p$ is a projection, by \Lemma{prefib-iso}
  $\pi^+(\g'_X) = i^+(\g'_Z)$, $\g_Y= i^+(\g_Z)$, and $\g_{Y/X}=
  \g_{Y/Z} + i^+(\g_{Z/X})$. Let $\af_Z$ be the Lie algebroid that is
  generated by $\g'_Z + \g_{Z/X}$ and $\af^1_Y$ the Lie algebroid that
  is generated by $\pi^+_1(\g'_Z) + \g_{Y/Z}$; now since
  $\g_{Z/X}\subset \g'_Z$ and $\g_{Y/Z}\subset \pi^+(\g'_X)$ we
  actually have $\af_Z= \g'_Z$ and $\af^1_Y= \pi^+(\g'_X) \subset
  \af_Y$. It should then be clear that (\ref{log-cond}) implies
\begin{displaymath}
  \Imo (i^*(\g'_Z) \to i^*(\g_Z)) = \Imo (\af^1_Y \to
  i^*(\g_Z)) \subseteq \Imo (\af_Y \to i^*(\g_Z))
\end{displaymath}
and
\begin{displaymath}
  \Imo(p^*(\g'_X)\to p^*(\g_X))= \Imo (\af_Z \to p^*(\g_X)).
\end{displaymath}
so (\ref{log-cond}), and hence (\ref{base-eq}) holds for
$i$ and $p$ separately.  Moreover, by \Lemma{pullfree} the
pull-back of a locally free and transitive Lie algebroid again
is locally free and transitive. In the
complex analytic case we also remark that if $M$ is generated by a
coherent $\Oc_Y$-module, then $H^\bullet(i_+(M^\bullet))$ is
generated by a coherent $\Oc_Z$-module, i.e.
$i_+(D^b_{sl}(\g_Y))\subset D^b_{sl}(\g_Z)$. This proves that the
conditions are satisfied for $i$ and $p$ separately.

$i_+(D^b_{s}(\g_Y))\subset D^b_s(\g_Z)$: 
First, since $\g_Z$ is locally free and transitive one can prove as in
\cite[VI, 7.8]{borel:Dmod}, see also c) in the proof of
\Theorem{kashiwara}, that the right $\Dc(\g_Y)$-module
$\Dc_{Z\leftarrow Y}$ is locally free. Therefore, since $i_*$ is exact
one has
\begin{displaymath}
  i_+(M)= i_*(\Dc_{ Z \leftarrow Y} \otimes_{\Dc(\g_Y)}M).
\end{displaymath}
For a coherent $\Oc_Y$-submodule $M^0\subset M$ we let
\begin{displaymath}
  i_*(\Dc^n_{Z\leftarrow Y}\cdot M^0) \subset i_*(\Dc_{Z\leftarrow  Y}\otimes_{\Dc(\g_Y)}M)
\end{displaymath}
be the image of the canonical morphism $i_*(\Dc_{Z\leftarrow Y}
\otimes_{\Oc_Y}M^0) \to i_*(\Dc_{Z\leftarrow Y}\otimes_{\Dc(\g_Y)}M)$;
this is a coherent $\Oc_Z$-module.  If $N \subset i_*(\Dc_{Z
  \leftarrow Y} \otimes_{\Dc(\g_Y)}M)$ is a coherent $\Oc_Z$-submodule
there exists such a coherent submodule $M^0$ and an integer $n$ such
that
\begin{displaymath}
  N \subseteq i_*(\Dc^n_{Z\leftarrow Y}\cdot M^0).
\end{displaymath}
Therefore by (\ref{base-eq})
\begin{equation*}
  \Dc(\g'_Z) N  \subseteq \Dc(\g'_Z) i_*(\Dc^n_{Z\leftarrow  Y}\cdot M^0)  \subseteq i_*(\Dc^n_{Z\leftarrow  Y} \cdot\Dc(\g'_Y)M^0).
\end{equation*}
Since $\Dc(\g'_Y)M^0 \in \coh (\Oc_Y)$ and $\Dc^n_{Z\leftarrow  Y}$ is 
coherent both as right $\Oc_Y$-module and left $\Oc_Z$-module the right 
side is coherent;  the left side is locally  generated by its
sections and $\Oc_Z$ is noetherian, hence $\Dc(\g'_Z) N \in \coh (\Oc_Z)$.

It remains to prove: If $M\in D^b_{sl}(\g_Z)$ and the restriction of
$p: Z \to X$ to $\supp M$ is proper, then $p_+(M)\in D^b_{s}(\g_X)$:
Here $\g_Z = p^*(\g_X)\oplus q^*(T_Z)$, where $q$ is the projection on
$Y$, so $p^*(\g_X)$ is a Lie algebroid on $Z$ in a natural way. The
$\Oc_Z$-module $p^*(\Dc(\g_X))$ is a right $p^{-1}(\g_X)$-module and
left $\g_Y$-module, but since the action of $p^{-1}(\g_X)$ commutes
with that of $q^{-1}(\Oc_Y)$ it follows that $p^*(\Dc(\g_X))$ is a
right module over the Lie algebroid $p^*(\g_X)$.  Changing right to
left we have a $(p^*(\g_X) ,\g_Z)$-bimodule
\begin{displaymath}
  \Dc_{X\leftarrow Z} = \omega_Z \otimes_{\Oc_Z} p^*(\omega_X^{-1})\otimes_{\Oc_Z}p^*(\Dc_X)
\end{displaymath}
which is coherent over $p^*(\g_X)$; hence, since $M\in \coh (\g_Z)$,
one has
\begin{displaymath}
  \Dc_{X\leftarrow Z}\otimes^L_{\Dc(\g_Z)}M\in D^b_{\coh}(p^*(\g_X)).
\end{displaymath}
Now $\Omega^\bullet_{Z/X}\otimes_{\Oc_Z}\Dc(\g_Z)= q^*(\Omega^\bullet
_Y)\otimes_{\Oc_Z}\Dc(\g_Z)$, where $\Omega^\bullet_Y$ is the de~Rham
complex on $Y$, is a left resolution of the bimodule $\Dc_{X\leftarrow
  Z}$ which is free as right $\g_Z$-module (see \cite[VI,
5.3]{borel:Dmod}), so
\begin{displaymath}
 \Dc_{X\leftarrow Z}\otimes^L_{\Dc(\g_Z)}M =
 \Omega^\bullet_{Z/X}\otimes_{\Oc_Z} M.
\end{displaymath}
Then (\ref{base-eq}) implies that $\Dc_{X\leftarrow
  Z}\otimes^L_{\Dc(\g_Z)}M\in D^b_s(p^*(\g_X))$ where the index
$s$ denotes smoothness along $p^*(\g'_X)\subset p^*(\g_X)$. Since
$D^b_s(p^*(\g_X))$ is generated by its homology objects it now
suffices to prove the following: If $N\in \coh_s(p^*(\g_X))$, then
\begin{displaymath}
  Rp_* (N) \in D^b_s(\g_X).
\end{displaymath}
As the restriction of $p$ to $\supp M$ is proper the proof in
\Proposition{direct-coh} implies that each homology module
$R^ip_*(N)\in \coh (\g_X)$. Now the problem is local on $X$ so it
remains to see that if $L_i \subset R^ip_*(N)$ is an $\Oc_X$-coherent
submodule then $\Dc(\g'_X)L_i \in \coh (\Oc_X)$.  To see this first
note that $N$ contains a $p^*(\g'_X)$ module $N_i$ which is coherent
over $\Oc_Z$, and such that
\begin{displaymath}
  L_i \subset  K_i:=\Imo \{R^ip_*(N_i) \to R^ip_*(N)\}.
\end{displaymath}
Then since $K_i\in \coh (\Oc_X)$ is a $\g'_X$-module the assertion follows.
\end{pf}

\begin{remark}
  Is \Theorem{direct-smooth} true for all locally free Lie algebroids?
  The problem occurs when $\pi: Y\to X$ is a closed embedding and
  $\Dc_{X\leftarrow Y}$ is not free as $\g_Y$-module.  We need a
  resolution $F^\bullet \to \Dc_{X\leftarrow Y}$ of the
  $(\pi^{-1}(\g_X),\g_Y)$-bimodule $\Dc_{X\leftarrow Y}$ such that
  each term $F^i$ is a $(\pi^{-1}(\g_X),\g_Y)$-bimodule which is flat
  as $\g_Y$-module and  provided with a filtration by
  $\pi^{-1}(\Oc_X)$-coherent submodules $F^i_n$ such that $\cup_{n\geq
    0}F^i_n= F^i$ and $\pi^{-1}(\g'_X)F^i_n \subset F^i_n \g'_Y$.
\end{remark}

\subsection{Completion}
Let $(A, \mf, k)$ be a noetherian local $k$-algebra and $\g' \subseteq
\g$ Lie algebroids over $A$. Let $\hat M$ be  the $\mf$-adic
completion of an  $A$-module $M$. The Artin-Rees lemma implies that if 
 $M$ is of finite type, then $\hat M = \hat A \otimes_A M$ and $M
 \subseteq \hat M$, in particular $\hat {\g}' = \hat A \otimes_A \g'$, $\hat {\g} = \hat A \otimes_A \g$.  Since $A
\subseteq \hat A$ and derivations of $A$ are continuous in the
$\mf$-adic topology we have that $\Der_k(A) \subseteq \Der_k(\hat A)
$.  More generally, if $M$ is a $\g$-module of finite type, its
defining homomorphism $\rho : \g \to \cf (M)$ lifts to homomorphisms
$\hat {\rho}: \hat {\g} \to \cf (\hat M)$ and $\tilde {\rho} : \hat
{\g} \to \hat A \otimes_A \cf (M)$.  Note that if $M$ is not
of finite type over $A$, then the natural injective mapping $\hat A
\otimes M \to \hat M$ is not surjective, so the $\hat \g$-module of
finite type $\tilde {M}:= \hat A \otimes_A M$ is not $\mf$-adically
complete, and $\hat M$ is not of finite type over $\hat \g$.  For
instance, $\Dc(\hat \g) \cong \hat A \otimes_A \Dc(\g) \ncong \hat
{\Dc}(\g)$.  Still, the following result is very useful for checking
if $M$ is smooth because of the big supply of invariants in $\tilde
M$ (for an application, see \cite[Sec.
4]{kallstrom:2smooth}). \typeout{Referens till II!!}

\begin{prop}\label{complete}
  The following are equivalent:
   \begin{enumerate}
   \item $M$ is smooth along $\g'$.
   \item $\tilde M$ is smooth along $\hat {\g}'$.
   \end{enumerate}
\end{prop}

\begin{pf}
  This is similar to the proof of \Proposition{inv-equiv}.  If $M^0$
  is a $\g'$-module of finite type over $A$, then $\hat M^0 \cong \hat
  A \otimes_A M^0$ is a $\hat {\g}'= \hat A\otimes_A \g'$-module of
  finite type over $\hat A$. It is also clear that $\Dc (\hat \g)
  \cdot \hat M^0 = \hat A \otimes _A M = \tilde M$.  This proves that
  $(1)$ implies $(2)$.
  
  Assume $(2)$. Let $M^0 $ be an $A$-submodule of finite type
  that generates $M$ over $\g$, and set $M_1= \Dc(\g')M^0$.  We
  have $\hat M^0 \cong \hat A\otimes_A M^0 $, and since
  $\Dc(\hat {\g}')= \hat A \otimes_A \Dc(\g')$ it follows that
  $ \tilde M_1 = \hat A \otimes_A M^1 = \hat A \otimes _A
  \Dc(\g')M^0 = \Dc(\hat {\g}') \hat M^0$.  $(2)$ and
  \Lemma{all-coh} implies that $\tilde M_1$ is of finite
  type over $\hat A$. Hence there exists a surjection $\hat f
  : \oplus ^n \hat {A} \to \hat {A} \otimes M_1$ for some
  integer $n$ where $\hat f = \id \otimes f$ for some
  homomorphism $f: \oplus ^n A \to M_1$.  Now since the functor
  $\hat {A}\otimes_{A}\cdot $ is faithfully flat it follows
  that $f$ is surjective. This implies $(1)$.
\end{pf}

\begin{remark} Analogous results hold for holomorphic localisation.
  Let $\Oc_x$ be the local ring at a closed point $x$ of a
  complex variety and $ \Oc^{an}_x$ be the stalk of holomorphic
  functions at $x$. Then $\Oc_x \to \Oc^{an}_x$ is a faithfully
  flat ring extension, so $M_x$ is smooth if and only if
  $M^{an}_x:= \Oc^{an}_x \otimes_{\Oc_x}M_x$ is smooth along
  $\Oc_x^{an}\otimes_{\Oc_x} \g'_x$.
\end{remark}

\subsection{GAGA}
\label{gaga-section}
Let $X$ be a projective complex variety, $X_h$ the associated compact
complex analytic space, and
\begin{displaymath}
  \pi : (\Oc_{X_h},X_h) \to (\Oc_X , X)
\end{displaymath}
the canonical continuous homomorphism of ringed spaces. If
$\g_{X_h}\in \Lie_{X_h}$ we let $\coh^l(\g_{X_h})\subset \coh
(\g_{X_h})$ be  the sub-category of coherent $\g_{X_h}$-modules
$M^h$ which contain a coherent $\Oc_{X_h}$-module $M^h_0$ such that
$\Dc(\g_{X_h})\cdot M^h_0 = M^h$ ($M^h_0$ is a lattice in $M^h$). The
following implication from GAGA \cite{serre:gaga} is well-known in
the case of $T_X$-modules.
\begin{thm}($X$ is a projective complex variety)\label{gaga-eq}
  The functor 
  \begin{displaymath}
    \pi^* : \Mod (\Oc_X) \to \Mod (\Oc_{X_h}), \quad M \mapsto
    \Oc_{X_h}\otimes_{\pi^{-1}(\Oc_X)}\pi^{-1}(M)
  \end{displaymath}
  induces equivalences of categories
  \begin{eqnarray}
    \nonumber \Lie_X & \cong & \Lie_{X_h}, \text{ and }\\
    \nonumber \coh(\g_X) & \cong & \coh^l (\pi^*(\g_X)).
  \end{eqnarray}
\end{thm}
\begin{pf}
  $\pi^*: \Lie_X \to \Lie_{X_h}$ is fully faithful: Let
  $\g^1_X, \g^2_X \in \Lie_X$. Then clearly $ \pi^*( \g^1_X),
  \pi^*( \g^2_X) \in \Lie_{X_h}$ and if $\phi : \g^1_X\to
  \g^2_X $ is a homomorphism the mapping $ \pi^*(\phi)
  :\pi^*(\g^1_X)\to \pi^*( \g^2_X)$ is a homomorphism of Lie
  algebroids on $X_h$.  We have a canonical mapping
  \begin{eqnarray}\nonumber
    && \Hom_{\Lie_X}(\g^1_X, \g^2_X) \xrightarrow{\pi^*\ }
    \Hom_{\Lie_{X_h}}(\pi^*(\g^1_X), \pi^*(\g^2_{X}))\\
    \nonumber & \subseteq & \Hom_{\Oc_{X_h}}( \pi^*(\g^1_X),
    \pi^*(\g^2_X)) \cong \Hom_{\Oc_X}( \g^1_X, \g^2_X).
      \end{eqnarray}
      Here GAGA implies the last isomorphism and that $\pi^*$
      is injective.  Then if $\phi_h= \id \otimes \phi :
      \pi^*(\g^1_X)\to \pi^*(\g^2_X)$ is a Lie homomorphism it
      is obvious that $\phi: \g^1_X \to \g^2_X$ also is a Lie
      homomorphism, so $\pi^*$ is surjective.

$\pi^*: \Lie_X \to \Lie_{X_h}$ is essentially surjective: By
GAGA if $\g_{X_h} \in \Lie_{X_h} \subset \coh(\Oc_{X_h})$ then
there exists a unique $\Oc_X$-module $\g_X$ such that
$\g_{X_h}= \pi^*(\g_X)$ and $\pi^{-1} (\g_X)\subset \g_{X_h}$.
We need to prove that $[\pi^{-1}(\g_X),
\pi^{-1}(\g_X)]\subseteq \pi^{-1}(\g_X)$. The mapping $\ad_h :
\g_{X_h}\to \cf_{X_h}(\g_{X_h})$ is a 1st order differential
operator on the $\Oc_{X_h}$-module $\g_{X_h}$ with values in
$\cf_{X_h}(\g_{X_h})$, hence it can be identified with an
$\Oc_{X_h}$-linear homomorphism $\ad_h: P_{X_h}^1(\g_{X_h})
\to\cf_{X_h}(\g_{X_h})$  \Sec{linear}, hence by GAGA: $P_{X_h}^1(\g_{X_h}) =
P^1_{X_h}(\pi^*(\g_X)) = \pi^*(P^1_X(\g_X))$, and
\begin{displaymath}
  \cf_{X_h}(\g_{X_h}) \subset
  Hom_{\Oc_{X_h}}(P^1_{X_h}(\pi^*(\g_X)), \pi^*(\g_X))=
  \pi^*(Hom_{\Oc_X}(P_X^1(\g_X), \g_X)),
\end{displaymath}
implying that $\cf_{X_h}(\g_{X_h}) = \pi^*(\cf_X(\g_X))$; and
$\ad_h$ is induced by a unique homomorphism $\ad:
P^1_{X}(\g_X) \to \cf_X(\g_X) $; therefore $[\pi^{-1}(\g_X),
\pi^{-1}(\g_X)]= \pi^{-1}(\ad) (\pi^{-1}(\g_X))
(\pi^{-1}(\g_X))\subset \pi^{-1}(\g_X)$.

$\pi^*: \coh(\g_X) \to \coh^{l}(\pi^*(\g_{X}))$ is an
equivalence of categories: Let $M^h$, $ N^h \in
\coh^l(\g_{X_h})$, $\phi^h \in \Hom_{\Lie_{X_h}}(M,N)$, and
$M_0^h \subset M^h$, $N_0^h \subset N^h$ be coherent
$\Oc_{X_h}$-submodules such that 
\begin{displaymath}
  \Dc_{X_h}\cdot M^h_0= M^h, \text{ and } \Dc_{X_h}\cdot N_0^h = N^h.
\end{displaymath}
Since $X_h$ is compact there
exists an integer $l$ so that $\phi^h(M_0^h)\subset
\Dc_{X_h}^lN_0^h$; let $\phi^h_0$ denote the restriction of
$\phi^h$ to $M^h_0$. By GAGA there exist unique $\Oc_X$-modules
$ M_0$, $N^l_0$ and a homomorphism of $\Oc_X$-modules $\phi_0:
M_0 \to N^l_0$ such that $\phi^h= \pi^*(\phi_0): M_0^h=
\pi^*(M_0)\to \Dc_{X_h}\cdot N^h_0 = \pi^* (N^l_0)$. Since
$\Dc_{X_h} = \Dc_{X_h}(\pi^*(\g_X)) = \pi^*(\Dc_X)$ it follows
that there exist coherent $\g_X$-modules $M$, $N$ such that
$M_0 \subset M$, $N_0 \subset N$, $\pi^*(M)= M^h$,
$\pi^*(N)=N^h$, and a homomorphism $\phi: M\to N$, extending
$\phi_0$ and satisfying $\pi^*(\phi)=\phi^h$; since by GAGA
$\pi^*(\phi_0)=0$ implies $\phi_0=0$, this proves that $\pi^*$
is essentially surjective and fully faithful.
\end{pf}

Let us express the parallel for Lie algebroids the
 fact that on a projective complex variety $H^1(X,
\Oc_{X}^*) = H^1(X_h, \Oc^*_{X_h})$ because $\Pic_X=
\Pic_{X_h}$:
\begin{corollary}\label{gaghom}
  Let $X$ be a projective complex variety. Then there is an
  isomorphism of complex vector spaces
  \begin{displaymath}
   \pi^{-1}: H^2(X, \Omega_X^{\geq 1}) \cong H^1(X_h,
   \Omega^{1,cl}_{X_h}).
  \end{displaymath}
\end{corollary}
The proof is immediate from \Theorem{gaga-eq}; see also~\ref{picardliealg}.

Letting $X_h$ be the holomorphic localisation of a projective
algebraic manifold $X$ we now say a few words about non-algebraic
$\g_{X_h}$-modules.  Let $S\subset X_h$ be a hyper-surface, put
$X^0_h= X\setminus S$, and let $M$ be a coherent $\g_{X_h}$-module
such that $M_{X^0_h}\in \coh^l(\g_{X^0_h})$. It is in general a
difficult problem to see when $M\in \coh^l (\g_{X_h})$, but we  note
the following rather trivial fact:
\begin{lemma}\label{discrete}  Let $S\subset X_h$ be  a discrete subset,
  $X^0_h= X_h \setminus S$, $M\in \coh (\g_{X_h})$, and $M_{X^0_h}\in \coh
  (\Oc_{X^0_h})$. Then $M \in \coh^l (\g_{X_h})$.
 \end{lemma}
 \begin{pf}
   Let $x\in S$. Since $M\in \coh (\g_{X_h})$ and $S$ is discrete
   there exists a neighbourhood $\Omega$ of $x$ such that $\Omega \cap
   S= \{x\}$ and an $\Oc_{\Omega}$-coherent submodule
   $M^1_{\Omega}\subset M_{\Omega}$, generating $M_{\Omega}$ over $\Dc
   (\g_{\Omega})$. As $M_{\Omega \cap X^0}\in \coh (\Oc_{X^0\cap
     \Omega})$ and the ring $\Dc(\g_{\Omega})$ is
   noetherian~\ref{noetherianrings}, there exists an integer $k$ such
   that $(\Dc^k (\g_{X_h})_{\Omega} M^1)_{\Omega \cap X_h^0} =
   M_{\Omega \cap X^0_h}$.  Now define an $\Oc_{X_h}$-coherent
   extension of $M_{X^h_0}$ over $\Omega$ by $M^0_{\Omega}= \Dc^k
   M^1_{\Omega}$.  Doing this for each point in $S$ we get a coherent
   $\Oc_{X_h}$-submodule $M^0 \subset M$ such that $\Dc (\g_{X_h}) M^0
   = M$.
 \end{pf}
 Assume that $M = M(*S)$ and that $M$ is coherent over $\Oc_{X_h}(*S)$.  A
 nice example by Deligne in \cite{malgrange:reseaux} shows that $M$
 need not contain a coherent $\Oc_{X_h}$-module $M^0$, generating $M$
 over $\Oc_{X_h}(*S)$; in particular one must expect $\coh
 (\g_{X_h})\nsubseteq \coh^l (\g_{X_h})$.  On the other hand, B.
 Malgrange [loc. \ cit.] proves
 \begin{displaymath}
   \coh (\Oc_{X_h}(*S))\cap \coh (T_{X_h})\subset \coh^l(T_{X_h}),
 \end{displaymath}
 generalizing Deligne's result $\coh (\Oc_{X_h}(*S))\cap
 \coh^{tf}_{rs}(T_{X_h})\subset \coh^l_{rs}(T_{X_h})$, where
 $\coh^{tf}_{rs}(T_{X_h})$ is the category of torsion free coherent
 $T_{X_h}$-modules with (curve) regular singularities
 \Defn{regsing-def}.  Therefore, if $\g_{X_h}$ is {\it transitive}, by
 \Theorem{gaga-eq} we have
 \begin{displaymath}
   \coh(\Oc_{X_h}(*S))\cap \coh (\g_{X_h}) \subset \coh (\g_X)   \quad 
   \quad\text{(algebraic modules)}.
 \end{displaymath}
 The coherent $\Oc_X$-submodule $M^0 \subset M \in \coh
 (\Oc_{X_h}(*S))\cap\coh(T_{X_h})$ is defined by Deligne and Malgrange
 using a certain spectral conditions for the action of $t\partial_t$,
 where $t=0$ is a local equation at a non-singular point of $S$,
 giving `normal Jordan generators' for $M$ as in \cite{levelt}. One
 may speculate whether it is possible apply this idea to modules $M$
 over a non-transitive Lie algebroid $\g_{X_h}$, defining $M^0 \subset
 M \in \coh (\Oc_{X_h}(*S))\cap\coh (\g_X)$ by a spectral condition
 for the action of some distinguished sub-algebroid in $\g_{X_h}$.

\section{Prolonging over $\codim \geq 2$}
\label{prolonging}
\subsection{Coherent extensions}\label{cohprolong}
When $S\subset X$ is a closed subset in a (complex analytic) space we
say that $\codim_X S \geq 2$ if every point in $S$ is contained in a
subspace $S'$, and there exists an open set $\Omega$ such that
$S\cap \Omega \subset S'\cap\Omega$ and $\codim_{\Omega} S' \geq 2$.

Consider the normalization morphism of a reduced space or a general
reduced scheme $X$,
\begin{displaymath}
f:  X' \to X.
\end{displaymath}
So $X'$ is the disjoint union of the normalization of all the
irreducible components of $X$.  This is a finite morphism when $X$ is
a space.  Also, $f$ is almost always finite when $X$ is a locally
noetherian scheme; for instance, it suffices that $\Oc_x$ is a G-ring
at each point $x\in X$.

\begin{lemma}\label{hartog} ($X$ is either a reduced complex analytic space or a
  locally noetherian reduced scheme such that the normalization
  morphism $f$ is finite) Let $S\subset X$ be a closed subset such
  that $\codim_X S \geq 2$. Put $\Omega= X\setminus S$ and let $j:
  \Omega \to X$ be the inclusion.  Then:
  \begin{enumerate}
  \item $j_*(\Oc_\Omega) \in \coh (\Oc_X)$; when $X$ is normal, then
    the canonical mapping $\Oc_X \to j_*(\Oc_\Omega)$ is an isomorphism;
  \item ($X$ is a reduced noetherian scheme) Let $N_\Omega$ be a
    coherent torsion free $\Oc_\Omega$-module (see
    \Section{basic-facts} for the definition of torsion freeness). Then $j_*(N_{\Omega})$
    is coherent;
  \item ($X$ is a reduced complex analytic space) Assume that $N_X$ is
    an $\Oc_X$-module which is locally generated by its sections and
    that $N_\Omega$ is a coherent torsion free $\Oc_\Omega$-module.
    Then $j_*(N_\Omega)$ is coherent.
  \end{enumerate}
\end{lemma}
Here (2) is essentially covered by a result of Grothendieck \cite[Exp.
VIII, p. 6, Prop. 2.3]{grothendieck:cohomologie}, and (3) is a very
special case of a result due to Serre \cite{serre:prolongement}, but
we provide simple direct proofs. When $X$ is non-reduced
$j_*(\Oc_\Omega)$ need not be coherent.
\begin{pf} 
  (1): The second part is in the complex analytic case 
  Riemann's ``Hebbarkeitssatz'', see
  \cite{scheja:fortsetzung}; in the algebraic case we recall
  the proof: One can assume that $X$ is affine; then $\codim_X
  (X\setminus \Omega) \ge 2$ implies that $\Oc_{X}(\Omega)
  \subseteq \Oc_x$ when $\hto (x) =1$; by normality $\Oc_X(X) =
  \cap_{\hto (x)=1}\Oc_x$; hence $\Oc_X(\Omega) \subseteq
  \Oc_X(X)$.
  
  For the first part, consider the Cartesian square
  \begin{equation}\nonumber
\begin{CD}
  \Omega' @>i>> X' \\ @Vg VV @VVf V \\ 
  \Omega @>j>> X.
\end{CD}
\end{equation}
One checks that $\Oc_X \subset f_*(\Oc_X')$ since $X$ is reduced (the
normalization is the union of all irreducible components of $X$, and
these are integral). Then since formation of direct image is local on
the base, so $j^*f_*=g_*i^*$; hence, using $i_*i^*(\Oc_{X'}) =
\Oc_{X'}$,
\begin{displaymath}
  j_*(\Oc_{\Omega}) = j_*j^*(\Oc_X) \hookrightarrow j_*j^*f_*(\Oc_{X'})
  = j_* g_* i^*(\Oc_{X'}) = f_*i_*i^*(\Oc_{X'}) = f_*(\Oc_{X'}).
\end{displaymath}
As $f$ is finite, $f_*(\Oc_{X'})\in \coh (\Oc_X)$; from
$j_*(\Oc_\Omega) = j_*j^*(\Oc_X)$ it follows that $j_*(\Oc_\Omega)$ is
locally generated by its sections; hence since $\Oc_X$ is noetherian
(\ref{noetherianrings}), $j_*(\Oc_\Omega)\in \coh (\Oc_X)$.
  
(2): Put $N^\vee_\Omega= Hom_{\Oc_\Omega}(N, \Oc_\Omega)$. Since $X$
is noetherian, by Gabriel's theorem there exist $N_e, N^\vee_X \in
\coh (\Oc_X)$ such that $j^*(N_e)= N_\Omega, j^*(N^\vee_X) =
N^\vee_\Omega$ \cite[Prop 2.]{borel-serre:riemann-roch}.  As
$N_\Omega$ is coherent and torsion free there exists a non-empty open
subset $\Omega_1 \subseteq \Omega$, with dense intersection with each
irreducible component of $X$, such that $N_{\Omega_1}$ is locally
free; therefore, if $K$ is the kernel of the canonical morphism
$N_\Omega \to N^{\vee \vee}_\Omega$, then $\dim \supp K < \dim
\Omega$, but $N_\Omega$ is torsion free so $K=0$. Hence, $j_*$ being
left exact, the canonical mapping
  \begin{displaymath}
    j_*(N_\Omega) \to j_*(N_\Omega^{\vee \vee})
  \end{displaymath}
  is injective. Since $N_\Omega$ has some coherent extension $N_e$ to
  $X$, so $j_*(N_\Omega) = j_*j^*(N_e)$, implying that $j_*(N_\Omega)$
  is locally generated by its sections. Therefore, $\Oc_X$ being
  noetherian, $j_*(N_\Omega) \in \coh (\Oc_X)$ will follow if
  $j_*(N_\Omega^{\vee \vee}) \in \coh (\Oc_X)$. By adjunction
\begin{displaymath}
  j_*(N^{\vee \vee}) = j_* Hom_{\Oc_\Omega}(j^*(N^\vee_X), \Oc_\Omega) 
  = Hom_{\Oc_X}(N^\vee_X , j_*(\Oc_\Omega));
\end{displaymath}
hence $j_*(N^{\vee \vee}) \in \coh (\Oc_X)$ by (1).

(3): By assumption there exist $N_e =N, N^\vee_X = N^\vee \in \coh (\Oc_X)$ such
that $j^*(N_e)= N_\Omega$ and $j^*(N^\vee_X)= N^\vee_\Omega$. Then
proceed as in $(2)$.
\end{pf}

\begin{remarks} \label{extremark}
  \begin{enumerate}
  \item Serre has proven a stronger result than (3) above; see
    \cite[Th.1]{serre:prolongement}. Namely, the following statements
    are equivalent on a normal complex analytic space $X$ for a
    coherent and torsion free $\Oc_\Omega$-module $M_{\Omega}$: (i)
    for each point $s\in S$ there exists an open set $U\subset X$,
    $x\in U$, such that $M_{\Omega\cap U}$ is generated by its
    sections; (ii) $j_*(M_\Omega)\in \coh (\Oc_X)$.
    
    Some other facts pertaining to the complex analytic case are
    exhibited in {\it loc.\ cit.}\/: (i) Let $X= \Cb^2$ and $S=
    \{0\}$.  There exist invertible sheaves on $X\setminus S$ that are
    not extendible by a coherent sheaf on $X$.  (ii) Let $X = \Cb^n$,
    $n \geq 3$, $S=\{0\}$, and let $i: X\setminus S \to X$ be the
    inclusion. If $M$ is a reflexive $\Oc_{X\setminus S}$-module, then
    if $i_*(M)$ is coherent it is also reflexive, but it need not be
    locally free when $M$ is locally free.
  
  \item Assume that $j:\Omega \to X$ is an open inclusion of complex
    manifolds, and put $S= X\setminus \Omega$, where $\codim_X S \geq
    2$. Deligne proves the following result for the Lie algebroid
    $T_\Omega$ \cite[Cor.  5.8]{deligne:eq}: If a $T_\Omega$-module
    $M_\Omega$ is coherent over $\Oc_\Omega$, then $j_*(M_\Omega)$ is
    coherent over $\Oc_X$.  The proof is based on Hironaka's
    resolution of singularities and Grauert's finiteness theorem.
    Comparing to Serre's results we get: if $U\subset X$ is a sufficiently
    small open neighbourhood of a point in $S$, then $M_{U \setminus
      (U\cap S)}$ is generated by its global sections. Deligne's
    result cannot be generalised to general Picard Lie algebroids,
    which can be seen by taking an invertible sheaf $M_\Omega =
    \lambda_\Omega$ on $\Omega = \Cb^2 \setminus \{0\}$ that is not
    extendible to coherent $\Oc_{\Cb^2}$-module and put $\g_\Omega =
    \cf (\lambda_\Omega)$.
  
  \item (1) in \Lemma{hartog} implies that the canonical mapping
    $T_X(X) \to T_X(\Omega)$ is a bijection when $X$ is normal: if
    $\delta \in T_X(\Omega)$ and $f\in \Oc_X(X)$, then $\delta (f) \in
    \Oc_X(\Omega) = \Oc_X(X)$. In particular, $j_*(T_\Omega)=T_X$.
  \end{enumerate}
\end{remarks}

The generalization to torsion modules is quite hard and we need some
preparation.

 We let $\prof M_x$ be the depth of an
$\Oc_x$-module $M_x$ of finite type.

\begin{thm}\label{groth-coherent}($X$ is a scheme which is locally a 
  closed sub-scheme of a non-singular scheme, or a complex analytic
  space) Let $Y\subset X$ be a closed subset and $j: U \to X$ be the
  open immersion of $U:= X\setminus Y$. Let $M\in \coh (\Oc_X)$. When
  $X$ is a scheme assume:
  \begin{enumerate}
  \item [(1)] $\prof M_x \geq s-1$ for all $x\in U$ such that $\codim
    _Y (\{x\}^-\cap Y) =1$.
  \end{enumerate}
When $X$ is a complex analytic space assume:
\begin{enumerate}
\item[(1')] Each point $y\in Y$ has an open neighbourhood $U_y \subset
  X$ such that $\prof M_x \geq s + \dim_y Y$ when $x\in U_y$.
\end{enumerate}

Then $R^pj_*j^*(M)\in \coh (\Oc_X)$ when $p \leq s-2$.
\end{thm}

\Theorem{groth-coherent} was proven for schemes by Grothendieck
\cite[Exp. VIII, p. 13, Th.3.1]{grothendieck:cohomologie} and when $X$
is a complex analytic space it is due to Siu and Trautmann
\cite{siu:local},\cite{trautmann}.  We include a proof below when $X$
is a scheme which is a little more direct than Grothendieck's proof.

Let $\pdo N$ be the projective dimension of a module $N$ over a
regular ring $R$; let $\dop (R)$ be the Krull dimension of $R$ so
$\prof R = \dop (R)$.  We
recall the Auslander-Buschsbaum relation \cite[Th.  19.1]{matsumura}:
\begin{equation}\label{eq:auslander-buschbaum}
  \pdo N + \prof N = \dop (R).
\end{equation}

That $(1')\Rightarrow (1)$ can be seen using the next lemma, but we
will not need this fact; the opposite implication is false.
\begin{lemma}\label{equivalent-prof} Let $A$ be a local ring which is
  a quotient of a regular ring, let $P \subset A$ be a prime ideal,
  and let $q$ be an integer. Let $N$ be an $A$-module of finite type.
  Then if $\prof N \geq q + \coht P $ it follows that $\prof N_P \geq
  q$.
  \end{lemma}
  Here $\coht P = \dop (A/P)$ is the coheight of a prime ideal. A
  stronger result is true without assuming that $A$ is a quotient of a
  regular ring, namely that $\prof N \leq \prof N_P + \coht P$
  \cite[Lem. 9.3.2]{brodmann-bruns}.
  \begin{pf}   We consider $M$ as a module over a
    regular ring $R$ and $P$ as prime ideal in $R$. Then $\pdo M \leq
    \dop (R) - (q + \coht P) = \dop (R_P) - q$.  Therefore $\pdo M_P =
    \sup \{i: Ext^{i}_{R_P}(M_P, R_P) \neq 0\} \leq \pdo M \leq \dop
    (R_P)-q$. Hence by (\ref{eq:auslander-buschbaum}) $\prof M_x \geq
    q$.
  \end{pf}

  The following fact will be needed:
\begin{lemma}[{\cite[Exp. III, p. 10, Prop. 3.3]{grothendieck:cohomologie}}]\label{eq:depth-local} Let $X$ be a locally noetherian scheme.
  The following are equivalent:
  \begin{enumerate}
  \item $H^i_Z(M)=0$, when $i\leq n-1$;
  \item $\prof M_z \geq n$ for all $z\in Z$.
  \end{enumerate}
\end{lemma}

We begin the proof of \Theorem{groth-coherent}, following
Grothendieck, by showing that one may assume that $X$ is non-singular:
The assertion of the theorem is of local nature, so we can assume that
$X$ can be embedded in a non-singular scheme.  Let $i: X \to X_{ns}$
be a closed immersion where $X_{ns}$ is a regular scheme and put
$V_1=i(V), Y_1=i(Y), U_1 = X_{ns}\setminus Y_1$, so
$\codim_{V_1}Y_1\geq 2$. Let $j_1 : U_1 = X_{ns}\setminus Y_1 \to
X_{ns}$ be the open inclusion and $i_1: U \to U_1=i(j(U))$ be the
restriction of $i$ to $U\subset X$.  Then $j_1\circ i_1 = i\circ j$,
so there exists an isomorphism of functors on $D(\Oc_X)$, $ i_*Rj_*
\cong R{j_1}_* {i_1}_*$, since $i_*$ and ${i_1}_*$ are exact. In
particular, using the fact that the direct image is local on the base,
so ${i_1}_* j^* = j^*_1 i_*$, we get
    \begin{displaymath}
i_*R^sj_*j^*(M) =  R^s{j_1}_*{i_1}_*j^*(M) = R^s{j_1}_* j_1^* i_*(M)
    \end{displaymath}
    for each integer $s$.  Here $i_*(M)$ is a module such that
    $\prof j_1^*(i_*(M))_x \geq s + \dim_y Y_1$. So if one
    proves the theorem     when $X$ is non-singular then $i_*(R^lj_*j^*(M)) \in \coh
    (\Oc_{X_{ns}})$ when $l\leq s-2$; hence, since $i_*$ is fully
    faithful, $R^lj_*j^*(M) \cong i^*i_*(R^lj_*j^*(M)) \in \coh
    (\Oc_X)$.

Until the end of the proof we will now asume that $X$ is a
non-singular scheme.

  \begin{lemma}\label{huvudlemma}
Let $N\in \coh (\Oc_X)$, $S= \supp N$, and $Z= \overline{S\cap
      U} \cap Y$. Let $r\in \Zb_{\geq 0}$ and assume that $\dop
    (\Oc_y) \geq r$ when $y\in Z$. Then 
\begin{displaymath}
  Hom_{\Oc_X}(N, R^{s}j_*j^*(\Oc_X)) \in \coh (\Oc_X)
\end{displaymath}
when $s\leq r-2$.
\end{lemma}
\begin{pf}
  The long exact sequence in homology of the distinguished triangle
\begin{displaymath}
  R\Gamma_Y(\Oc_X)\to \Oc_X \to Rj_*j^*(\Oc_X)\to 
\end{displaymath}
shows that we may instead prove
\begin{displaymath}
  Hom_{\Oc_X}(N, H^{s+1}_Y(\Oc_X))\in \coh (\Oc_X)
\end{displaymath}
when $s\leq r-2$.

Put $N^0 = \Gamma_Y (N)$ and $N^1= N/N^0$. Applying the contravariant
left exact functor
\begin{displaymath}
T:  N \mapsto Hom_{\Oc_X}(N,H^{s+1}_Y(\Oc_X)) 
\end{displaymath}
to the short exact sequence
\begin{displaymath}
  0\to N^0 \to N \to N^1  \to 0
\end{displaymath}
results in a short exact sequence
\begin{displaymath}
0\to T(N^1)\to T(N) \to T(N)'\to 0
\end{displaymath}
where $T(N)':=\Imo (T(N)\to T(N^0)) \subset T(N^0)$ is locally
generated by its sections.  Since $\supp N^0 \subset Y$, clearly
$T(N^0) \in \coh (\Oc_X)$; therefore $T(N)'\in \coh (\Oc_X)$ because
$\Oc_X$ is noetherian. Now it will follow that $T(N)\in \coh(\Oc_X)$
if $T(N^1)=0$. First we have
\begin{displaymath}
  T(N^1)= Hom_{\Oc_X}(N^1, H^{s+1}_Y(\Oc_X))= Hom_{\Oc_X}(N^1,
  \Gamma_Z(H^{s+1}_Y(\Oc_X))).
\end{displaymath}
Second, since $R\Gamma_Z(\Oc_X) = R\Gamma_Z R\Gamma_Y(\Oc_X)$ there
exists a spectral sequence
$H^p_Z(H^q_Y(\Oc_X))\Rightarrow H^{p+q}_Z(\Oc_X)$;
in particular $\Gamma_Z(H^{s+1}_Y(\Oc_X))$ is a subquotient of
$H^{s+1}_Z(\Oc_X)$; since $\prof \Oc_y = \dop (\Oc_y) \geq r$ when
$y\in Z$ it follows that $H^{s+1}_Z(\Oc_X)=0$, and hence
$\Gamma_Z(H^{s+1}_Y(\Oc_X))=0$, when $s+1 \leq r-1$
\Lem{eq:depth-local}.  Therefore $T(N^1)=0$ when $s \leq r-2$.
\end{pf}

\begin{pfof}{\Theorem{groth-coherent} when $X$ is non-singular}    
  Since $M$ is
coherent and $X$ is non-singular there exists a biduality isomorphism
\begin{displaymath}
  M\cong RHom_{\Oc_X}(RHom_{\Oc_X}(M, \Oc_X), \Oc_X).
\end{displaymath}
Then  adjunction gives
\begin{eqnarray*}
  Rj_*j^*(M) &=& Rj_*j^*RHom_{\Oc_X}(RHom_{\Oc_X}(M, \Oc_X), \Oc_X) \\&=&
  Rj_*RHom_{\Oc_U}(j^*RHom_{\Oc_X}(M, \Oc_X), \Oc_U)  
\\ &=&   RHom_{\Oc_X}(RHom_{\Oc_X}(M, \Oc_X), Rj_*(\Oc_U)).
\end{eqnarray*}
giving a spectral sequence
\begin{displaymath}
  Hom_{\Oc_X}(Ext_{\Oc_X}^q(M, \Oc_X), R^{q+p}j_*j^*(\Oc_X)) \Rightarrow R^pj_*j^*(M).
\end{displaymath}
It now suffices to prove that the left side is coherent when $p\leq
s-2$.

Put $S^q= \supp Ext_{\Oc_X}^q(M, \Oc_X)$ and $Z^q = \overline
{(S^q\cap U)}\cap Y$.

b) $\dop (\Oc_y) \geq q+s$ when $y\in Z^q$: For any $y\in Y$ there
exists $x\in S^q \setminus Y$ such that $y\in \{x\}^- \cap Y$ and
$\dop (\Oc_y)= \dop (\Oc_x)+1$. Then $\prof M_x \geq s-1$ by $(1)$,
so by (\ref{eq:auslander-buschbaum})
\begin{displaymath}
 \dop (\Oc_y)= \dop (\Oc_x) +1 = \pdo M_x + \prof M_x +1 \geq q + s-1 + 1 = s+q,
\end{displaymath}
where the inequality is because $x\in S^q$.

c) By b), $\prof \Oc_y = \dop (\Oc_y) \geq q+s$ when $y\in Z^q$.
Hence by \Lemma{huvudlemma}, with $N:=Ext^q_{\Oc_X}(M,\Oc_X)$,
\begin{displaymath}
  Hom_{\Oc_X}(Ext_{\Oc_X}^q(M,\Oc_X), R^{q+p}j_*j^*(\Oc_X)) \in \coh(\Oc_X)
\end{displaymath}
when $q+p \leq q+s -2$, i.e. \ when $p \leq s-2$.
\end{pfof}

Let $M\in \coh (\Oc_X)$ and put $V= \supp M$, $k= \dim V$. We say that
$M$ is (k-) {\it pure}\/ if every coherent $\Oc_X$-submodule $N\subset
M$ satisfies $\dim \supp N = k$.  When $M$ is not pure we will
consider its filtration by supports (= the Krull filtration), which we
recall.  Let $S(M)$ be the maximal coherent submodule $N\subset M$
such that $\dim \supp N < \dim \supp M$. Then setting $M_0=M$ and
inductively $M_{i+1}=S(M_i)$ one has a decreasing filtration of length
$\leq \dim X$
\begin{displaymath}
  0 = M_{n}  \subset M_{n-1} \subset \cdots   \subset M_0 =M
\end{displaymath}
where successive quotients $L_i= M_{i}/M_{i+1}$ are pure.

\begin{corollary}\label{coherent} ($X$ is of the same type as in \Theorem{groth-coherent})
  Let $M$ be a pure coherent $\Oc_X$-module with support $V\subset X$.
  Let $Y\subset V$ be a closed subset such that $\codim_V Y \geq 2$
  and let $j: U \to X$ be the associated open immersion where $U=
  X\setminus Y$. Then $j_*j^*(M)\in \coh (\Oc_X)$.
\end{corollary}
\begin{pf} When $X$ is a scheme one may consider $M$ as a module on
  $V$ over the ring $\Oc_V:=\Oc_X/ \Ann (M)$, so one may assume that
  $V=X$. Then since $M$ is pure with support $X$ we have  $\prof M_x
  \geq 1$ when $\codim_Y (\{x\}^- \cap Y)$, satisfying $(1)$ in
  \Theorem{groth-coherent}, so we may take $s=2$.
  
  To prove the assertion in the complex analytic case we employ
  $(1')$, which turns out a bit more complicated in spite of
  $(1')\Rightarrow (1)$!  Put $d= \dim V$ and let
  \begin{displaymath}
    S=  \{x\in \prof M_x < d(M_x) \} \subset V 
  \end{displaymath}
  be the locus of points where $M$ is not Cohen-Macaulay. We assert
  that this is a closed (analytic) subset such that $\codim_V S \geq
  2$. Assuming this for a while we put $Y_1 = S \cup Y$, so $\codim_V
  Y_1\geq 2$ and $M_{U_1}$ is Cohen-Macaulay, where $U_1:= X \setminus
  Y_1$.  Letting $j_1: U_1 \to X$ be the open immersion
  \Theorem{groth-coherent} implies that $(j_1)_*(j_1)^*(M)\in \coh
  (\Oc_X)$; since $j_*j^*(M)\subset (j_1)_*(j_1)^*(M) $ is a submodule
  which is locally generated by its sections and $\Oc_X$ is
  noetherian, it follows that $j_*j^*(M)\in \coh(\Oc_X)$.
  
  It remains to prove $\codim_V S \geq 2$. Replacing the structure
  sheaf on $V= \supp M$ we may assume that $X= V$ and $n=d$, and since
  the assertion is local in $X$ we may assume that $X\subset X_{ns}$ where
  $X_{ns}$ is either a complex manifold or a non-singular scheme; we thus
  consider $M$ as an $\Oc_{X_{ns}}$-module with support on $X$. Let $m\geq 1$
  be an integer and put
  \begin{eqnarray*}
    D_m&:=& \{x\in X : \prof M_x \leq \dop(M_x)-  m\}\\ 
&=& \{x\in X : \pdo    M_x \geq \dop (\Oc_x) - \dop (M_x) + m\} \\
&=& \{x\in X : \pdo M_x \geq n-d +m\}
  \end{eqnarray*}
  where the second line follows from (\ref{eq:auslander-buschbaum}).
  Hence
  \begin{displaymath}
    D_m = \bigcup_{r\geq n-d +m} \supp Ext^r_{\Oc_{X_{ns}}}(M, \Oc_{X_{ns}}).
  \end{displaymath}
  Since $M\in \coh (\Oc_{X_ns})$ this shows that $S=D_1$ is closed in
  $X$ (see also \cite[exercise 24.2]{matsumura}), and a complex
  analytic subspace when $X$ is a complex analytic space. Now
  $\codim_{X_{ns}} Ext^r_{\Oc_{X_{ns}}}(M,\Oc_{X_{ns}})\geq r \geq n-d
  +1$, so $\codim_X S \geq 2$ will follow if the $\Oc_x$-module
  $Ext^{n-d+1}_{\Oc_{X_{ns}}}(M,\Oc_{X_{ns}})_x \cong
  Ext^{n-d+1}_{\Oc_x}(M_x, \Oc_x)$ has dimension $\leq d-2$ ($M$ is
  coherent). Put $A= \Oc_x$, $N=M_x$, and let $p\subset A$ be a prime
  ideal of height $n-d+1$. Then since $N$ is a pure $A$-module of
  finite type, $\prof N \geq 1 $ and $Ext^{n-d+1}_A(N, A)_p =
  Ext^{n-d+1}_{A_p}(N_p,A_p)$; hence by
  (\ref{eq:auslander-buschbaum}), $\pdo N_p \leq n-d+1 -1 = n-d$, and
  therefore $Ext^{n-d+1}_A(N, A)_p=0$; whence
  \begin{displaymath}
    \dim  Ext^{n-d+1}_{\Oc_{X_{ns}}}(M, \Oc_{X_{ns}})_y \leq d-2.
  \end{displaymath}
\end{pf}
\begin{remark}
  Note that in our definition a variety is always locally a closed
  subscheme of a non-singular scheme, so \Theorem{groth-coherent} and
  \Corollary{coherent} are applicable to spaces.
\end{remark}
\subsection{Prolongation of smoothness}
Let $\pi :Y\to X$ be a proper homomorphism of spaces, $S$ a closed
subset of $X$, and put $X_1= X \setminus S$.  The mapping $j: X_1 \to
X$ is the canonical inclusion and $\pi_1 : Y_1\to X_1$ the
corresponding base change of $\pi$. So the projection on the other
coordinate $i:Y_1 \to Y$ is an open inclusion and we  have a
Cartesian diagram:
\begin{equation}\label{diagram:prolong}
  \begin{CD}
    Y_1 @>i>>Y \\ 
@V{\pi_1}VV @VV{\pi}V
    \\ X_1 @>j>> X.
  \end{CD}
\end{equation}
If $\g_Y \in \Lie_Y$ we put $\g_{Y_1}= i^+(\g_Y)$
(see~\ref{pullbacklie}). For a a sub-space $V\subset X$ the $V$-{\it
  component } \/of a coherent $\g_X$-module $M$ is the maximal
coherent sub-module of $M$ whose support is contained in $V$.
\begin{thm}\label{relative}  
  Let $\g'_Y\subset \g_Y$ be a sub-algebroid and $M$ a coherent
  $\g_Y$-module containing, locally in $X$, a coherent
  $\Oc_Y$-submodule $M^0 \subset M$ such that $M= \Dc(\g_Y)M^0$; this
  always holds in the algebraic case.  Put $V=\supp M$ and let $V=
  \cup^k_i V_i$ be an irreducible decomposition of $V$. Assume that
  $V_i\cap V_j \cap \pi^{-1}(S)=\emptyset$ when $i\neq j$, and if
  $V_i\cap \pi^{-1}(S)\neq \emptyset$, then the $V_i$-component of $M$
  is pure and $\dim \pi (V_i)\geq \dim S + 2$. Assume also:
\begin{enumerate}
\item The canonical homomorphism $\pi^{*}\pi_*(M)_y \to
  M_y$ is surjective when $y \in \pi^{-1}(S)$;
  
\item The canonical homomorphism $\pi^*\pi_*(\g'_Y)_y
\to \g'_y$ is surjective when $y\in \pi^{-1}(S)$.
\end{enumerate} 

Then if $i^!(M)$ is point-wise smooth along $\g'_{Y_1}$ it follows
that $M$ is point-wise smooth along $\g'_Y$.
\end{thm}

We shall consider non-relative versions of \Theorem{relative}.  First
the case when $M$ is torsion free is singled out, which may be of
special interest, and we then need only the more elementary
\Lemma{hartog}.
\begin{corollary}\label{pure}
  Let $\g'_X$ be a sub-algebroid of a Lie algebroid $(X, \Oc_X ,
  \g_X)$ and $M$ a torsion free $\g_X$-module.
  \begin{enumerate}
  \item  ($X$ is a variety) The following are equivalent:
  \begin{enumerate}
  \item $M$ is smooth at all points.
  \item $M$ is smooth at all points of height $1$.
  \end{enumerate}
\item ($X$ is a  complex analytic space) The following
are equivalent:
\begin{enumerate}
\item $M$ is smooth at all points.
\item There exists an open set $\Omega \subset X$ with
$\codim_X (X\setminus \Omega) \geq 2$ such that $M$ is smooth
at all points $p \in \Omega$.
  \end{enumerate}
  \end{enumerate}
\end{corollary}
\begin{pf} We need to prove  $(b)\Rightarrow (a)$, and this
  follows from \Theorem{relative} by letting $\pi : X \to X$ be the
  identity mapping. Here follows a direct proof: Let $p$ be a point in
  $X$. Then there exists an affine (Stein) neighbourhood $\Omega$ of
  $p$ and a coherent $\Oc_\Omega$-submodule $M_\Omega^0 \subset
  M_\Omega$ such that $\Dc (\g_\Omega)M^0_\Omega = M_\Omega$. Put
  $M'_\Omega := \Dc (\g'_\Omega)M^0_\Omega$. As $\Omega$ is affine
  (Stein) and $\Dc(\g'_\Omega)$ is a union of coherent
  $\Oc_\Omega$-modules, it follows that $\Dc(\g'_\Omega)$ generated by
  its sections (Cartan A in the analytic case) and since $M'_\Omega$
  is generated by its sections over $\Dc (\g'_\Omega)$ it is also
  generated by its sections over $\Oc_\Omega$. Now, since the strong
  support of the $\g'_{\Omega}$-module $M'_\Omega$ is closed
  \Prop{point}, $(b)$ implies both in the algebraic and the analytic
  case that there exists an algebraic (analytic) subset $Z\subset
  \Omega$, $\codim_\Omega Z\geq 2$, such that $M'_{\Omega\setminus
    Z}\in \coh (\Oc_{\Omega\setminus Z})$.  Letting $j:
  \Omega\setminus Z \to \Omega$ be the inclusion mapping we have a
  canonical mapping
  \begin{displaymath} M'_p \to j_*(M'_{\Omega\setminus Z})_p
  \end{displaymath}
  which is injective because $M'_\Omega$ is pure and $\supp M'_\Omega
  = \Omega$. Then  \Lemma{hartog} implies that $j_*(M'_{\Omega\setminus
    Z})_p$ is the stalk of a coherent module, thus it is of finite
  type over $\Oc_p$; since $\Oc_p$ is noetherian it follows that
  $M'_p$ is of finite type.
\end{pf}

\begin{cor}\label{purestcor}Let $\g'_X \subset \g_X$ be a sub-algebroid,
  $L\in \coh (\g_X)$ be pure, $V= \supp L$, and $Y \subset V$ a
  subspace such that $\codim_V Y \geq 2$; put $U= X\setminus Y$. Then
  if $L_U$ is point-wise smooth along $\g'_U$ it follows that $L$ is
  point-wise smooth along $\g'_X$.
\end{cor}
\begin{pf}
  The proof is similar to that of \Corollary{pure}. Namely, the
  problem is local at a point $p\in Y$ so one may assume that $L$
  contains a coherent $\Oc_X$-submodule $L^0 \subset L$ such that $L=
  \Dc(\g_X)L^0$. Put $L'= \Dc(\g'_X)L^0$. By assumption $L'_U\in
  \coh(\Oc_U)$ and since $L'$ is pure
  \begin{displaymath}
    L' \subset j_*(L_U')
  \end{displaymath}
  The right side is coherent by \Corollary{coherent}; $\Oc_X$ is
  noetherian (\ref{noetherianrings}) and $L'$ is locally generated by
  its sections; hence $L'\in \coh (\Oc_X)$.
\end{pf}

\begin{corollary}\label{point-cor}($X$ is a space) Let $M\in \coh
  (\g_X)$ be torsion free, and $S\subset X$ a divisor such that $M$ is
  smooth along $\g'_X \subseteq \g_X$ at all points in $X\setminus S$.
  Then if $M_x$ is smooth at one point $x$ in each connected component
  of the non-singular locus of $S$, it follows that $M$ is point-wise
  smooth in $X$.
\end{corollary}
Using \Corollary{purestcor} it is straightforward to modify the proof
below to find a corresponding result for any pure module; we omit this
generalisation.
\begin{remark} Let $X$ be a complex analytic manifold, $S\subset X$ a
  divisor, $\g_X = T_X$ and $\g'_X= T_X(I_S) $ its sub-sheaf of
  derivations that preserve the ideal $I_S$ of $S$.  Let $M$ be a
  $T_X$-module, coherent and torsion free over $\Oc_X(*S)$.  In
  \cite[Th\'{e}or\`{e}me 4.1(i)$\Rightarrow$(ii)]{deligne:eq} it is
  proved that $\phi^!(M)$ is smooth along $\phi^+(T_X(I_S)) =
  T_C(I_{\phi^{-1}(S)})$ for each curve $\phi : C\to X$ if this holds
  for each curve $\phi : C\to \Omega$ where $\Omega\subset X$ is a
  subset with $\codim_X (X\setminus \Omega) \geq 2$.  By
  \Corollary{curve-testc} this curve test is equivalent to $M$ being
  smooth along $T_X(I_S)$ (this also follows from [loc. cit]); hence
  Deligne's result follows from \Corollary{point-cor}.
\end{remark}
\begin{pf} The assertion being local in $X$ we can clearly assume that $M$ contains a coherent
  $\Oc_X$-submodule $M^0$ such that $M= \Dc (\g_X)M^0$. Put
  $M^1= \Dc (\g'_X)M^0 \in \coh (\g'_X)$ and set
  \begin{eqnarray}\nonumber
    S_s & = & \{x \in S : M_x \text{ is smooth along
    $\g'_x$}\}\\ \nonumber & = & \{x\in S : M^1_x\text { is of
    finite type over $\Oc_x$}\}.
  \end{eqnarray}
  Now $\ssupp M^1$ is a union of irreducible spaces \Lem{strong-var}, $\ssupp M^1 \subset S$, and $S_s = S
  \setminus \ssupp M^1 \subset S$ is an open subset containing a point
  from each component of $S$, implying that
  \begin{displaymath}
    \codim_X (\ssupp M^1) \geq 2.
  \end{displaymath}
  Thus if $\Omega= X\setminus \ssupp M^1$, then $M_\Omega$ is
  point-wise smooth, so the assertion follows from
  \Corollary{pure}.
\end{pf}

The following is a sharpening of \Theorem{direct-smooth} for pure modules. 
\begin{corollary}\label{direct-smooth-pure} Make the same
  assumptions and use the same notation as in
  \Theorem{direct-smooth}. Let $M$ be a pure $\g_Y$-module which
  is smooth along $\g'_Y$ and put $V= \supp M$. Let $V_1$ be a
  subset of $V$ so that $\codim_V V_1 \geq 2$. Then if
  \begin{equation}
    \label{eq:log-cond-pure}
     \Imo (\pi^*( \g'_X)_y \to \pi^*(\g_X)_y) \subseteq  \Imo (\afr_y \to \pi^*(\g_X)_y)
  \end{equation}
  when $y\in V\setminus V_1$, it follows that $\pi_+(M)\in D^b_{s}(\g_X)$.
\end{corollary}
\begin{pf}
  The proof is close to that of \Theorem{direct-smooth} so we
  will avoid too much repetition. This time
  (\ref{eq:log-cond-pure}) implies
  \begin{equation}\label{base-eq-pure}
    (\pi^{-1}(\g'_X)\Dc^n_{X\leftarrow Y})_y  \subseteq (\Dc^n_{X\leftarrow
      Y} \afr_Y)_y \subseteq (\Dc^n_{X\leftarrow
      Y} \g'_Y)_y,
  \end{equation}
  when $y \in Y- V_1$.  We factorize the proper morphism $\pi: Y \to
  X$ into a closed embedding $i :Y \to Z= X\times Y$ and a projection
  $p :Z\to X$ where the restriction of $p$ to the support of $i_+(M)$
  is proper.  The right $\Dc(\g_Y)$-module $\Dc_{Z\leftarrow Y}$ is
  locally free, hence $i_+(M) = i_*(\Dc_{Z\leftarrow
    Y}\otimes_{\Dc(\g_Y)} M)$ and this is a pure $\Oc_Z$-module. Also,
  (\ref{eq:log-cond-pure}) holds for $i$ and $p$ separately where
  $\g'_Z = p^+(\g'_X)$. We may therefore prove the assertion when
  either $\pi$ is a closed immersion or $\pi$ is a projection such
  that the restriction of $\pi$ to $\supp M$ is proper
  \Prop{composed}.
  
  $\pi$ is a closed embedding: Put $W= \pi(V)$ and $W_1 = \pi (V_1)$,
  so $\pi_+(M)$ is pure with support $W$.
  Using the same argument as in the proof of \Theorem{direct-smooth}
  by (\ref{base-eq-pure}) $\pi_+(M)= \pi_*(\Dc_{X\leftarrow
    Y}\otimes_{\Dc(\g_Y)} M)$ is smooth along $\g'_X$ in $X\setminus
  W_1$; hence it is smooth everywhere \Cor{purestcor}.
  
  $\pi: Y= Z\times X \to X$ is the projection on the second coordinate
  and $M$ is a pure $\g_Y$-module such that the restriction $\pi' : V
  \to X$ of $\pi$ to the support $V=\supp M$ is proper: Again
  the proof is similar to that of \Theorem{direct-smooth}. We have
  \begin{displaymath}
    \Dc_{X \leftarrow Y}\otimes^L_{\Dc(\g_Y)} M =    \Omega^\bullet_{Y/X}\otimes_{\Oc_Y} M
  \end{displaymath}
  so this is a complex of pure coherent $\pi^*(\g_X)$-modules. Now
  (\ref{base-eq-pure}) implies that this complex  is smooth along $\pi^*(\g'_X)$
   at all points in  $Y\setminus V_1$, hence it is smooth everywhere
  \Cor{purestcor}. Now the proof can be finished as in
  \Theorem{direct-smooth}.
\end{pf}

\subsection{Proof of \Theorem{relative}}
In the proof, which  is divided into lemmas, we will use the notation in the
Cartesian diagram (\ref{diagram:prolong}).

\begin{lemma}\label{limsection}
  Let $M^i\in \coh (\Oc_Y)$, $i\in I$, be an inductive system
  over a directed partial ordered set $I$. Then
  \begin{enumerate}
  \item $\dirlim \ M^i$ is locally generated by its sections.
  \item $\pi_*(\dirlim \ M^i) = \dirlim \ \pi_*(M^i)$.
  \end{enumerate}
\end{lemma}
\begin{pf} $(1)$: Let $y \in Y$. In the algebraic case let $\Omega$ be
  an affine neighbourhood of $y$, in the complex analytic case
  $y$ is contained in a Stein neighbourhood $\Omega\subset Y$
  and `Theorem A' holds for Stein spaces (see
  \cite{grauert-remmert:stein}). Hence in either case $M^i_y
  \cong \Oc_y\otimes_{\Oc_Y(\Omega)} M^i(\Omega)$.  Since
  inductive limits commute with $\Oc_y\otimes \cdot$ (see
  \cite[II, 1.11]{godement:faisceaux}), we get
  \begin{displaymath}
    \Oc_y\otimes_{\Oc_Y(\Omega)} (\dirlim \ M)(\Omega) \cong
    \dirlim \ (\Oc_y \otimes_{\Oc_Y(\Omega)}M^i(\Omega)) \cong
    \dirlim \ M^i_y.
  \end{displaymath}
  
  $(2)$: This follows by formal reasons if one admits that by
  properness the  derived functor $R\pi_*$ has a right adjoint
  $\pi^!(\cdot)$; see \cite{hartshorne-res}, \cite{banica}.
  
  Here follows a more self-contained proof.  We have a canonical
  homomorphism $\phi : \dirlim \ \pi_*(M^i) \to \pi_*(\dirlim \ M^i)$.
  This is an isomorphism if it induces an isomorphism between stalks
  over any point $x\in X$. Therefore, by the last line in the
  proof of (1), it suffices to prove that the canonical mapping
  \begin{displaymath}
    \phi_x :\dirlim \ \pi_*(M^i)_x \to \pi_*(\dirlim \ M^i)_x
  \end{displaymath}
  is an isomorphism. Let $m\in \pi_*(\dirlim \ M^i)_x$ be
  represented by the section
  \begin{displaymath}
     m_\Omega \in \pi_*(\dirlim \ M^i)(\Omega) = (\dirlim \
     M^i) (\pi^{-1}(\Omega))
  \end{displaymath}
  where $\Omega\subset X$ is an open neighbourhood of $x$. If $U$ is
  an open subset in $\pi^{-1}(\Omega)$ we denote by $m_U$ the
  restriction of $m_{\Omega}$ to $U$. Every point $y\in
  \pi^{-1}(\Omega)$ has an open neighbourhood $U_y \subset Y$ and an
  integer $n_y$ such that $m_{U_y}$ is represented by a section
  $m^{n_y}_{U_y}\in M^{n_y}(U_y)$.  Since $\pi$ is proper it suffices
  with a finite number of the $U_y$ to cover $\pi^{-1}(x)$; let $U$ be
  their union and $n$ be the maximum of the corresponding integers
  $n_y$. Since $I$ is directed it follows that $m_U$ is represented by
  a section $m^n_U \in M^n(U)$.  Now $\pi : Y \to X$ is a continuous
  proper mapping between locally compact topological spaces, hence it
  is closed; hence the open sets $\pi^{-1}(\Omega')$, where $\Omega'$
  is an open subset of $X$, form a basis for the open neighbourhoods
  of $\pi^{-1}(x)$.  So there exists an open neighbourhood
  $\Omega'$ of $x$ such that $\Omega' \subset \pi (U) \subset \Omega$,
  and therefore $m^n_{\pi^{-1}(\Omega')}\in \pi_*(M^n)(\Omega')$; let
  $m^n_x \in \pi_*(M^n)_x$ be its germ at $x$.  Letting $i_n :
  \pi_*(M^n)_x \to \dirlim \ \pi_*(M^i)_x $ be the canonical mapping
  we have a mapping $ \psi_x : \pi_*(\dirlim \ M^i)_x \to \dirlim \ 
  \pi_*(M^i)_x $, well-defined by $m_x \mapsto i_n(m^n_x)$.  One
  checks that $\phi_x \circ \psi_x=\id$ and $\psi_x \circ \phi_x =
  \id$.
\end{pf}

\begin{lemma}\label{lattice} Let $M^i \in \coh (\Oc_Y)$, $i\in I$, be
  as above, $\phi_i : M^i \to M:= \dirlim \ M^i$ the canonical
  homomorphism, and $\psi_i : \pi^* \pi_* (M^i)\to M$ be the
  composition of the canonical homomorphisms $\pi^*\pi_* (M^i) \to M^i
  \to M$.  Assume that the canonical homomorphism $\pi^*\pi_*(M) \to
  M$ is surjective. Then for any $i\in I$ there exists, locally in $X$,
  an index $n = n(i)\in I$ such that
  \begin{displaymath}
    \Imo (\phi_i) \subset \Imo (\psi_{n(i)}).
  \end{displaymath}
\end{lemma}
\begin{pf}  \Lemma{limsection} and the fact that $\pi^*$ commutes with
  $\dirlim \ $ \cite[II, 1.11]{godement:faisceaux} implies that
  we have a surjection
\begin{displaymath}
  \dirlim \ \pi^* \pi_*(M^i) \cong \pi^*\pi_* \dirlim \ M^i \to
  \dirlim \ M^i.
\end{displaymath}
It is obvious that locally near a point $y\in Y$ one can find such an
index $n_y=n(i,y)$. Then, since $\pi$ is proper, for any $x\in X$ one
can reason as in the proof of \Lemma{limsection} to see that there
exists an index $n(i,x)\in I$ and a neighbourhood $\Omega$ of $x$ such
that
\begin{displaymath}
  \Imo (\phi_i)_{\pi^{-1}(\Omega)} \subset \Imo
  (\psi_{n(i,x)})_{\pi^{-1}(\Omega)}.
\end{displaymath}
\end{pf}

Of course, a proper and surjective morphism may blow up subsets of
codimension 2 to subsets of codimension 1. On the other hand, in the
complex analytic case Hartogs' theorem states (see e.g.  \cite[p.
42]{banica}): If $Y$ is a reduced complex Stein analytic space, $K$ a
compact subset such that $Y\setminus K$ has no relatively compact
irreducible components in $Y$, and the depth of $\Oc_y$ is $\geq 2$
for all $y\in K$, then $\Oc_Y(Y) = \Oc_Y(Y\setminus K)$. The main
ingredient in the proof of \Theorem{relative} is the following
coherence result which maybe can be thought of as Hartogs' type along
the fibre and Riemann's Hebbarkeitssatz on the base.

\begin{lemma}\label{rel-hartog}    Let $M^i \in \coh (\Oc_Y)$, $i \in
  I$, be an inductive system of modules such that, with $M_1= \dirlim
  \ M^i$, one has $i^*(M_1)\in \coh (\Oc_{Y_1})$.  Assume that $V=
  \supp M_1\subset X$ is a subspace and let $V= \cup^k_i V_i$ be an
  irreducible decomposition. Assume that $V_i\cap V_j \cap
  \pi^{-1}(S)=\emptyset$ when $i\neq j$, and if $V_i\cap
  \pi^{-1}(S)\neq \emptyset$ then the $V_i$-component $M_1^{V_i}
  \subseteq M_1$ is pure and $\dim \pi (V_i)\geq \dim S + 2$.  Put
\begin{displaymath} 
   N = \Imo (\pi^{*}\pi_*(M_1) \to M_1).
  \end{displaymath}

  Then if $M'$ is an $\Oc_{Y}$-submodule of $N$ which is locally
  generated by its sections, it follows that $M' \in \coh (\Oc_{Y})$.
\end{lemma}

Note that $\pi (V_i)$ will be a subspace since $\pi$ is proper, by
Remmert's theorem in the complex analytic case.

\begin{pf} Let $y \in \pi^{-1}(S)$. Since the problem is local near
  $\pi (y)$  we can and will shrink $X$ to any neighbourhood of
  $\pi(y)$ when necessary.\  
  
  a) It suffices to prove the assertion when $M_1$ is pure and $V$ is
  irreducible: As the $\Oc_X$-submodules $M^{V_i}_1 \subset M$ are
  pure we have $\oplus^k_{i=1} M_1^{V_i} \subset M_1$, and since
  $V_i\cap V_i \cap \pi^{-1}(S)=\emptyset$ when $i\neq j$, one may
  shrink $X$ so that $V_i\cap V_j = \emptyset$ when $i\neq j$; hence
  $\oplus^k_i M^{V_i}_1 = M_1$. Therefore, if $M'\cap M_1^{V_i} \in
  \coh (\Oc_Y)$, $i=1,\ldots k$, it will follow that $M' \in \coh
  (\Oc_Y)$.

  b) If $M_1$ is pure, then $\pi_*(M_1)$ is pure: Factorize $\pi =
  p\circ i$ where $i: Y \to Y \times X$ is the graph morphism and $p:
  Z= Y\times X\to X$ the projection on the second factor. Clearly
  $N=i_*(M)$ is pure with irreducible support $V_1=i(V)$. Since
  $\pi_*(M_1)= p_*i_*(M)$ it therefore suffices to prove: If $N$ is a
  coherent $\Oc_Z$-module such that the restriction of $p$ to a
  morphism $\supp N \to X$ is proper, then $p_*(N)$ is pure.
  
  By properness $W := \supp p_*(N)\subset X$ is a subspace. Suppose
  the contrary, that there exists a non-zero pure submodule $L \subset
  p_*(N)$ such that $\dim \supp L < \dim W$. Let $\phi:p^*(L)\to N$ be
  the composed morphism $p^*(L)\to p^*p_*(N)\to N$. Then $\dim \supp
  \Imo \phi < \dim \supp N$; $N$ is pure; hence $\phi = 0$; therefore,
  projections being flat, the canonical morphism $p^*p_*(N)\to N$ is
  $0$; hence $p_*(N)=0$, contradicting the assumption that $L\neq 0$.

  c) If $M_1$ is pure and $V$ is irreducible, then $H^0_{[T]}(M_1)=0$
  when $T= \pi^{-1}(S)$: Suppose the contrary, that
  $H^0_{[T]}(M_1)\neq 0$. Since $M_1$ is pure with irreducible support
  $V$ the submodule $H^0_{[T]}(M_1)\subset M_1$ also has support $V$;
  hence $V\subset T$, contradicting the relation $\dim \pi (V) \geq
  \dim \pi (T) +2$.

  d) $\pi_*(M_1)\in \coh(\Oc_X)$: By c)  and since $\pi_*$ is left exact we
  have an inclusion
\begin{displaymath}
  \pi_*(M_1) \subseteq \pi_*i_*i^* (M_1) = (\pi\circ  i)_*i^*(M_1) = (j\circ \pi^1)_*i^*(M_1) =
  j_*(\pi^1_*(i^*(M_1))).
\end{displaymath}
By \Lemma{limsection} $\pi_*(M_1)$ is locally generated by its
sections; then since $\Oc_X$ is noetherian (\ref{noetherianrings}) it
will follow that $\pi_*(M_1)\in \coh
(\Oc_X)$ if  we prove 
\begin{displaymath}
  j_*(\pi^1_*(i^*(M_1)))\in \coh (\Oc_X).  
\end{displaymath}

Since $i^*(M_1) \in \coh (\Oc_{Y_1})$ and properness is preserved
under base change, by Grauert's theorem \cite{grauert} or
Grothendieck's theorem \cite[3.2.1]{EGA3} 
\begin{displaymath}
  \pi^1_*(i^*(M_1)) \in \coh (\Oc_{X_1}).
\end{displaymath}
Also by properness and since $i^*(M_1)\in \coh (\Oc_{Y_1})$, there
exists a coherent $\Oc_Y$-submodule $M^0_1\subset M_1$ in a
neighbourhood of $\pi^{-1}(S)$ such that $i^*(M^0_1)= i^*(M_1)$. As
the formation of direct image is local on the base, so $j^*\pi_* =
\pi^1_* i^*$, it follows that $M:=\pi_*(M^0_1)$ is a coherent
extension of $\pi_1^*(i^*(M_1))$.  We thus have a coherent
$\Oc_X$-module $M$, such that $j^*(M)= \pi^1_*(i^*(M_1))$ is pure by
b), whose support $\pi (V)$ contains $S$ and $\codim_{\pi(V)} S \geq
2$; hence by \Corollary{coherent} $j_*(\pi^1_*(i^*(M_1)))\in
\coh(\Oc_X)$.

e) $M'\in \coh(\Oc_Y)$: Let $L^i$ and $N^i$ be the kernel and
image of the canonical homomorphism $\pi^*\pi_* (M^i)\to M^i$,
so we have a short exact sequence
 \begin{displaymath}
 0 \to L^i \to \pi^* \pi_*(M^i) \to N^i \to 0.
 \end{displaymath}
 Now $M^i \in \coh (\Oc_Y)$ and $\pi$ is proper so
 $\pi_*(M^i)\in \coh(\Oc_X)$ (\cite{grauert}, resp.
 \cite[3.2.1]{EGA3}), $\pi^*$ is right exact and $\Oc_Y$ is
 coherent, hence $\pi^*\pi_*(M^i)\in \coh (\Oc_{Y})$; therefore
 $L^i, N^i \in \coh(\Oc_Y)$ (\cite[ 0. Cor.  5.3.4]{EGA1}).
 Clearly we have canonical homomorphisms $N^i\to N^j$, and $L^i
 \to L^j$, $i \leq j$, so put $L = \dirlim \ L^i$ and $N =
 \dirlim \ N^i$ (the sheaves associated to the presheaves). Now
\begin{displaymath}
   \pi^* \pi_* (M_1) = \pi^*\dirlim \ \pi_*(M^i) = \dirlim \
   \pi^*\pi_*(M^i),
\end{displaymath}
where the first step follows from \Lemma{limsection} and the
second from \cite[II, 1.11]{godement:faisceaux}, and since the
index set $\Nb$ is directed the functor $\dirlim \ $ is exact,
 we get a short exact sequence
\begin{displaymath}
  0 \to L \to \pi^* \pi_* (M_1) \to N \to 0.
\end{displaymath}
By d) $ \pi^* \pi_* (M_1)\in \coh (\Oc_Y)$ (for the same reason
above that the $\pi^*\pi_*(M^i)\in \coh (\Oc_Y)$) and $L^i \in
\coh(\Oc_Y)$, so $L$ is locally generated by its sections
\Lem{limsection}; as $\Oc_Y$ is noetherian this implies that
$L\in \coh (\Oc_Y)$; hence $N\in \coh (\Oc_Y)$ \cite[ 0. Prop
5.4.2]{EGA1}, and again since $\Oc_Y$ is noetherian it follows
that $M' \in \coh (\Oc_Y)$.
\end{pf}

\begin{pfof}{\Theorem{relative}}
  Let $y\in \pi^{-1}(S)$. Since the problem is local near $x= \pi (y)$
  we can and will shrink $X$ to any small neighbourhood of $x$ when
  necessary and we can assume that $M= \Dc(\g_Y)M^0$.  Arguing as in
  a) in the proof of \Lemma{rel-hartog} it suffices to prove that
  $M':=\Dc(\g'_Y)M^0 \in \coh(\Oc_Y)$ when $M$ is pure and $V$ is
  irreducible, and we then know as in b) that $\pi_*(M)$ is pure.

  The modules $M^{i}= \Dc^{i}(\g_Y)M^0 \in \coh (\Oc_Y)$, $i\in \Nb$,
  form an inductive system by inclusion $M^{i}\subset M^{j}$ when
  $i\leq j$, and it follows from \Lemma{lattice}, possibly after
  shrinking $X$ around $x$, that there exists an integer $n_0$ such
  that $M^0$ is contained in the image of the canonical homomorphism
  $\pi^*\pi_*(M^{n_0}) \to M^{n_0}$; set $M_1= \Dc(\g'_Y)M^{n_0}$. 
  
  Now note that $\pi^* \pi_* (\g'_Y)$ is a Lie algebroid on $Y$ in a
  natural way and that $(2)$ is a homomorphism of Lie algebroids;
  since $N := \Imo (\pi^*\pi_*(M_1)\to M_1)$ clearly is a
  $\pi^*\pi_*(\g_Y')$-module $(2)$ implies that
  \begin{displaymath}
    \g'_Y\cdot N \subset N,
  \end{displaymath}
 and since $M^0 \subset N$ it follows that
  \begin{displaymath}
    M' \subset N.
  \end{displaymath}
  The $\Oc_Y$-coherent submodules $M_1^i= \Dc^i(\g'_Y) M^{n_0} \subset
  M_1$ form an inductive system by inclusion $M_1^i \subset M_1^{j}$,
  $j\geq i$, its limit $M_1 = \dirlim \ M_1^i$ is pure with
  irreducible support $V$, and $\dim \pi (V)\geq \dim S +2$.  By
  assumption $i^*(M_1)\in \coh (\Oc_{Y_1})$; therefore since $M' =
  \dirlim \ \Dc^i(\g'_Y)M^0$ is locally generated by its sections
  \Lem{limsection} it follows by \Lemma{rel-hartog} that $M' \in \coh
  (\Oc_Y)$.
\end{pfof}

\section{Regular singularities}
\label{regularity} 
\subsection{Torsion free modules}
Say that a set of vectors $\delta_i$ in a Lie algebroid $(A,\g
, \alpha)$ generates $\g$ if the subspace $\sum A\delta_i
\subseteq \g$ is not contained in any proper sub-algebroid of
$\g$. 

The following theorem follows from a straightforward application of
the integrability theorem \Th{gabber}.
\begin{thm} \label{fibre:1}($X$ is a  space) 
  Let $\g'_X $ be a sub-algebroid in a Lie algebroid $\g_X$ and $M$ a
  $\g_X$-module. Let $x$ be a  point in $X$, $M^0_x$ an
  $\Oc_x$-submodule of finite type that generates $M_x$ over $\g_x$,
  and $M_x' \subseteq M_x$ the $\g'_x$-module that is generated by
  $M^0_x$.  Let $\delta_i$ be generators of the Lie algebroid $\g'_x$.
  Then the following are equivalent:
  \begin{enumerate}
  \item $M$ is smooth at $x$;
  \item Each $\delta_i$ generates an $\Oc_x$-module of finite
type
    \begin{displaymath}
      \sum_k \delta^k_i \cdot M_x^0;
    \end{displaymath}
  \item $M'_x$ is smooth along all the vectors $\delta_i$ in
$\g'_x$ \Defn{along};
  \item Put inductively ${M'}_x^n = \g'_x {M'}_x^{n-1} +
{M'}_x^{n-1}$ for $n\geq 1$. The increasing sequence of mappings
        \begin{displaymath}
      k_x \otimes_{\Oc_x}{M'}^n_x \to k_x
      \otimes_{\Oc_x}{M'}^{n+1}_x\to \cdots
    \end{displaymath}
    is eventually surjective (here $k_x= \Oc_y/\mf_y$, while  $k_x=\Cb$ when
    $X$ is a complex analytic space).
  \end{enumerate}
\end{thm}

\begin{pf}  Clearly $M_x' = \cup_{n\geq 0} {M'}_x^n$, and  $M$ is smooth at
  $x$ if and only if $M_x'$ is of finite type over $\Oc_x$.
  
  The implications $(1)\Rightarrow (2)$,
  $(2)\Longleftrightarrow (3)$, $(1) \Rightarrow (4)$ are
  evident.  $(4)\Rightarrow (1)$: It suffices to note that by
  Nakayama's lemma, if  the mapping $
  k\otimes_{\Oc_x}{M'}_x^n \to k\otimes_{\Oc_x}{M'}_x^{n+1}$ is
  surjective  then ${M'}_x^n =  {M'}_x^{n+1}$.
  
  $(2)\Rightarrow (1)$: See \Section{diffoperators} for the notation.
  Let $\Dc_x = \Dc_x(\g'_X)$ be the enveloping ring of differential
  operators of $\g'_x$, $(\So_x , \cdot , \{\cdot , \cdot \})$ the
  Poisson algebra that is defined on $\So_x:= \gr \Dc_x$, and let
  $\sigma : \g'_x \to \So_x$ be the canonical morphism. Let $G(M_x')$
  be the graded module over $\So_x$ that is associated to the
  canonical filtration ${M'}^n_x = {\Dc}^n_x {M'}^0_x \subseteq
  {M'}^{n+1}_x$, $\cup {M'}^{n}_x = {M'}_x$.  Then $G(M'_x)$ is a
  module of finite type over the noetherian algebra $\So_x$, and the
  characteristic ideal $\Jo({M'})_x \subseteq \So_x$ is involutive
  \Th{gabber}; by assumption $\sigma (\delta_i) \in \Jo (M') _x$,
  hence $\sigma ([\delta_i , \delta_j]) = \{ \sigma (\delta_i), \sigma
  (\delta_j) \} \in \Jo ({M'})_x$; the $\delta_i$ generate $\g'_x$;
  hence $\sigma (\g'_x) \subseteq \Jo ({M'})_x$.  Therefore there
  exists an integer $n$ such that the $\So_x$-module of finite type
  $G(M'_x)$ is a module of finite type over $R_x:=\So_x/(\sigma
  (\g'_x))^n$; since clearly $R_x$ is of finite type over $\Oc_x$ it
  follows that $G(M'_x)$ is of finite type over $\Oc_x$. This implies,
  by an easy argument using the short exact sequences $0 \to
  {M'}^{n-1}_x \to {M'}^n_x \to G^n(M'_x) \to 0$, that $M'_x$ is of
  finite type over $\Oc_x$.
  \end{pf}
  
  We shall now work out curve tests for smoothness.  All curves $C$
  are assumed to be irreducible and non-singular, and morphisms $\phi
  :C\to X$ are assumed to be locally closed embeddings of spaces; when
  we say that $C$ is a curve in $X$ we mean such a morphisms $\phi$.
  An analytic curve in a variety $X$ is a morphism of schemes $\phi :
  C\to X$ where $C= \Spec k[[x]]$. If $\phi: C \to X$ is a curve we
  set \Sec{pullbacklie}
\begin{displaymath}
\g_C  :=   \phi^+(\g_X), \quad \text{and }\quad \g'_C  :=  \phi^+(\g'_X).
\end{displaymath}
We thus have homomorphisms of Lie algebroids $(\phi , \phi'):
\g_C \to \g_X$, and $\g'_C \to \g'_X$.

Let $\delta \in \g'_X(X_{reg})$ be a section defined over the
non-singular locus of $X$. Locally one can assume that $\delta = h
\delta^1$, where $\alpha (\delta^1) \in T_X(\Omega)$ is non-zero
outside an algebraic (analytic) subset $X_{crit}\subset X$ such that
$\codim_X X_{crit} \geq 2$.  Now if $X$ is a complex analytic space
and $x$ is a point in the non-singular locus of $\Omega$ there exists
locally integral curves for $\delta$. More precisely, by Frobenius'
theorem there exists a locally closed embedding of a curve $\phi : C
\to \Omega \subset X$ such that $\alpha(\delta)(I_{C_x})\subseteq
I_{C_x}$, where $I_{C_x}$ is the ideal of a germ $\phi (C)_x$ of the
curve at a point $x\in X$.  When $X$ is algebraic there exist analytic
(formal) solutions $\phi$, i.e.\ there exist ideals $\hat I (x)
\subset \hat \Oc_x$ so that $\alpha (\delta)(\hat I (x)) \subseteq
\hat I(x)$ and $\hat \Oc_x/\hat I(x)= k[[x]]$, but algebraic integral
curves in general do not exist.  If one insists that everything be
algebraic we say that $\g'_X\in \Lie_X$ has {\it good generators} \/
at a point $x\in X$ if:
\begin{description}
\item[(G)] There exist an open neighbourhood $U$ of $x$ and generators
  $\delta_1, \dots \delta_r\in \g'_X(U)$ of $\g'_U$ (as Lie algebroid)
  such that for every point $y\in U$ there exists a locally closed
  embedding of an algebraic curve $\phi^i : C^i \to U$, $y\in
  \phi^i(C^i)$, satisfying $\alpha(\delta_i(x))(I^i(x))\subset I^i(x)$
  ($I^i(x)$ is the ideal of $\phi^i(C^i)_x$).
\end{description}

For example, when $X$ is an algebraic manifold and $S$ a non-singular
divisor whose ideal is $I_S$, then $T_X(I_S)$ has good generators at
each point $x\in X$.

We say that a vector $\delta \in \g_x$  is {\it transversal} \/ to a divisor
$S$ at $x\in S$, if there exists an ideal $I_x\subset \Oc_x$ so that
$I_S + I_x=\mf_x$, $I_S \nsubseteq I_x$, and $\alpha (\delta)_x \in T_x(I_x)$ ($I_S$ is the
ideal of $S$); it is transversal to $S$ if it is transversal at each
point $y\in S$.

Say that a space $X$ is non-singular in codimension $1$ (or
non-singular outside $\codim \geq 2$) if there exists a closed subset
$C\subset X$ such that $\codim_X S\geq 2$ such that $X_0=X\setminus S$
is a non-singular space. When $X$ is a scheme this means that $\Oc_x$
is a regular ring when $x$ is a point of height $1$. This holds for
example when $X$ is normal.

Below we let $\phi^!(M)$ denote the non-derived inverse image
(\ref{pullbackmod}). In the notation above we have:
\begin{thm}\label{curve-test}($X$ is a space which is non-singular in
  codimension $1$; when $X$ is a variety we assume that the field $k$
  is uncountable) Let $M$ be a coherent torsion free $\g_X$-module
  which is point-wise smooth along $\g'_X$ outside a divisor $S$ in
  $X$. Let $S_0$ be an open subset of the non-singular locus $S_{reg}$
  of $S$, intersecting each component of $S_{reg}$.  Assume that
  $\g'_X$ locally has generators which are transversal to $S_0$. Then
  the following are equivalent:
\begin{enumerate} 
\item $M$ is point-wise smooth in $X$;
  \item $\phi^!(M)$ is smooth along $\g'_C$ for any analytic
curve $\phi : C \to X$ such that $\phi (C)\cap S \subseteq
S_0$;
  \item $\phi^!(M)$ is smooth along $\g'_C$ for any analytic
curve $\phi : C \to X$.
  \end{enumerate}
 Assume that $X$ is a
  variety and that $\g'_X$ has good generators at each point $x\in
  S_0$. Then $(1-3)$ are equivalent to:
  \begin{enumerate}\addtocounter{enumi}{3}
  \item $\phi^!(M)$ is smooth along $\g'_C$ for any curve $\phi : C \to X$
such that $\phi (C)\cap S \subseteq S_0$.
  \end{enumerate}

\end{thm}
The annoying condition that $\g'_X$ be generated by vectors that
preserve curves that are transversal to $S_0$ can be removed in the
following case:
  \begin{cor}\label{curve-testc}
    Make the same assumptions as in \Theorem{curve-test} except that
    $\g_X'$ need not have transversal generators. Assume instead that
    \begin{displaymath}
      T_X(I)\subseteq \alpha(\g_X')
    \end{displaymath}
($I$ is the ideal of $S$). Then the conclusion
    of  \Theorem{curve-test} holds.
  \end{cor}
\typeout{Check reference to Smooth II}
  \begin{pf}
    By \Theorem{curve-test} we first conclude that $M$ is smooth along
    the sub-algebroid that is generated by all transversal vectors.
    Then the assertion follows from \cite[Th.4.1]{kallstrom:2smooth}
    (the proof is by eliminating the transitive directions in
    $\g_X$).\typeout{Referens till II!}
  \end{pf}

\begin{remark} Let  $S\subset X$ be a divisor in a complex manifold
  and $M$ a $T_X$-module which is coherent as $\Oc_X(*S)$-module; let
  $I$ be the ideal of $S$.  By \cite[Th. 4.1]{deligne:eq} $M$ is
  smooth along $T_X(I)$ if for any curve $\phi :C\to X$ which is
  transversal to $S$ the inverse image $\phi^!(M)$ is smooth along
  $T_C(I_{\phi^{-1}(S)})$. The proof in {\it loc. \ cit.} \/ is based
  on complex analysis in two ways: (i) The connection in $X\setminus
  S$ is determined by its local system of horizontal sections; (ii)
  One proofs that $S$ is a countable union of certain {\it closed} \/
  subsets, so one of these will have an interior point by Baire's
  theorem. Our use of category argument below is more complicated.
\end{remark}

\subsection{Proof of \Theorem{curve-test}}

As in \Section{preserve} we put $\tilde {M} = \hat A \otimes_A M$ for
an $A$-module $M$, where $\hat A$ is the $\mf$-adic completion of a
local ring $(A,\mf)$. The proof of \Theorem{curve-test} is based on
the following lemma. We use the notation in \Theorem{curve-test}.
\begin{lemma}
\label{fibre:0}
Let $\Omega\subseteq X $ be an open subset such that $S_0=S\cap
\Omega$ is a non-empty open subset of $S$ defined by the zero locus of
a function $h\in \Oc_X(\Omega)$, and let $\delta \in \g'_X(\Omega)$.
Assume that for every point $x\in S_0$ that there exists an ideal
$\hat I(x)\subset \hat {\Oc}_x$ such that:
  \begin{enumerate}
\item $h\notin \hat {I}(x)$;
  \item $\delta_x \in \g'_x(\hat I(x))$;
  \item $\tilde M_x/\hat I(x)\tilde M_x$ is smooth along
$\delta_x$ in $\g_x'(\hat I(x))$ \quad (see~\ref{along}).
  \end{enumerate}
  
  Then there exists a (closed) point $y\in S$ such that $M_y$ is
  smooth along $\delta_y$ in $\g'_y$.
\end{lemma}

As a first preparation for the proof of \Lemma{fibre:0} we recall how
one may reduce certain local questions about analytic sheaves to
questions about modules over noetherian rings.  Let $X$ be a complex
analytic space, put $n=\dim X $, and let $x$ be a point in $X$. Then
there exists a neighbourhood $X'$ of $x$ and a finite morphism $f:X'
\to \Cb^n$.  Letting $K_0 \subset \Cb^n$ be a closed polydisc such
that $x\in K:=f^{-1}(K_0)$, then the ring $\Oc_{\Cb^n}(K_0)$ of germs
of holomorphic functions in neighbourhoods of $K_0$ is noetherian
\cite{frisch}; hence $A:=\Oc_X(K)$ is a noetherian ring since it is of
finite type over $\Oc_{\Cb^n}(K_0)$ \cite{grauert}. Moreover, $K$ is a
compact Stein set \cite[Ch IV, \S 1, Th.  3]{grauert-remmert:stein},
hence by Cartan's theorems A and B the functor  \begin{displaymath}
  \Gamma : \coh (\Oc_K) \to  A\text{-mod}_f,\quad  M \mapsto M(K)
\end{displaymath}
is an equivalence between the category of coherent $\Oc_X$-modules on
$K$ and the category of $A$-modules of finite type, with quasi-inverse
$M(K)\mapsto \Oc_K \otimes_{A}M(K)$; in particular, the stalk $M_x\cong
\Oc_x \otimes_A M(K)$.

Set $A= \Oc_X(K)$, where $K$ is either an affine variety or a compact
Stein set $K\subset X$ such that $A$ is noetherian (as above), $X_0
:=\mSpec A$, and $X:=\Spec A$. In the complex analytic case $X_0=K$.
The Nullstellensatz for either category implies that $X$ is an
 Jacobson scheme (for this notion, see \cite[Ch. 0, \S
2]{EGA1}), i.e.\ the canonical inclusion $i:X_0 \to X$ is a
quasi-homeomorphism, meaning that the mapping $i^{-1}$ defines a
bijection between the set of open subsets of $X$ and the set of open
subsets of $X_0$, and letting $U$ be an open subset of $X$, there is
an isomorphism $i^{\#}:\Oc_X(U) \to \Oc_{X_0}(i^{-1}(U))$ with the
$k$-algebra of regular algebraic functions on $i^{-1}(U)\subset X_0$,
defined by $\phi \mapsto i^{\#}(\phi): x \mapsto \phi \mod \mf_x =k$
($k$ is algebraically closed).  In particular, if $S_0 \subset X_0$ is
dense then $\cap_{x\in S_0}\mf_x =0$.

\begin{lemma}\label{nullstellen} 
  Let $\pfr =(h) \subset \Spec A$ be a principal prime ideal, and $S$
  be a dense subset of $\mSpec A/\pfr$. For each $\mf_x \in S$ choose
  an ideal $I(x) \in \Spec A$ such that $I(x) + \pfr \subset \mf_x$
  and $\pfr \nsubseteq I(x)$. Then if  $G$ is an $A$-module of finite
  type, we have
  \begin{displaymath}
    \pfr \notin \supp (\cap_{x\in S} I(x)G).
  \end{displaymath}
\end{lemma}
\begin{pf}   
  The $A$-module $G$ is noetherian so there exists an element $f\in A$
  such that $G_f\cong  A^{r}_f$ for some integer $r$, or $G_f=0$.
  Let $\bar G$ and $G^{(f)}$ be the image and kernel of the canonical
  mapping $G \to G_f$. Assume that $\bar G \neq 0$. Since
  $h\notin I(x)$ the set $S_1=\{x\in \mSpec A_f : \cap I(x)_f \subset
  \mf_x \}$ is dense in the Jacobson scheme $\Spec A_f$, so
\begin{displaymath}
  \bigcap_{x\in S} (I(x) \bar G) \subset \bigcap_{x\in S}(I(x)
  G_f) \cong \oplus (\bigcap _{x\in S} I(x))_f =0;
  \end{displaymath} 
  hence 
  \begin{displaymath}
    \cap_{x\in S} (I(x)G) \subset G^{(f)}.
  \end{displaymath}
Since any non-zero $\phi \in
  I(x)$ acts by an injective mapping on $\bar G \subset A_f^r$ it
  follows that 
  \begin{displaymath}
    \cap_{x\in S} (I(x)G) = \cap_{x\in S} (I(x) G^{(f)}).
  \end{displaymath}
  One may therefore assume that $G=G^{(f)}$. If $f\notin \pfr$, then
  \begin{displaymath}
      (\cap_{x\in S} (I(x) G^{(f)}))_{\pfr}\subset G^{(f)}_{\pfr}=0,
  \end{displaymath}
so $\pfr  \notin \supp \cap_{x\in S} (I(x)G)$. It remains to consider the case $f\in
  (h)$. As $G= G^{(f)}$ there exists an element $g\notin \pfr = (h)$
  such that $G_g= G_g^{(h)}$, and since $G_g$ is a noetherian
  $A_g$-module it has a finite filtration with successive quotients of
  the form $A_g/(h)$.  Now $\pfr+ I(x) \subset \mf_x \in \mSpec
  A/(h)$, and since $X$ is an  Jacobson scheme, implying
  \begin{displaymath} \bigcap_{x\in S} (I + \pfr) \subseteq
    \bigcap_{\pfr \subset\mf }\mf = \pfr,
  \end{displaymath} one gets $\cap_{x\in S} (I(x) \frac {A_g}{(h)}) =
  0$. Let now
\begin{displaymath}
  0 \to \frac{A_g}{(h)} \to N \to\frac {A_g}{(h)} \to 0
\end{displaymath}
be a short exact sequence of $A_g$-modules. As $h\notin I(x)$ any
non-zero element $\phi\in I(x)$ acts as an injective mapping on
$A_g/(h)$ (the right copy in the above exact sequence),  implying
\begin{displaymath}
  \bigcap_{x\in S} (I(x) N) = \frac {A_g}{(h)} \cap  \bigcap_{x\in S} (I(x) N)  =
  \bigcap_{x\in S}  (I(x) \frac {A_g}{(h)}) =0.
\end{displaymath}
An induction over the length of the filtration of $G_g$ shows that
\begin{displaymath}
  \cap_{x\in S} (I(x)G_g)=0,
\end{displaymath}
and since $G^{(g)}_{\pfr} =0$ when $g\notin \pfr$ it follows that
$(\cap_{x\in S} (I(x)G))_{\pfr}=0$.
\end{pf}

\begin{lemma}\label{inclusion}
  Let $S\subset V(\pfr)$ be a dense subset. For each $x\in S$ choose
  an ideal $\hat {I}(x)\subset \hat A_x$ such that $\hat {I}(x) + \hat
  {\pfr} \subset \hat {\mf}_x$ and $\pfr \nsubseteq \hat I (x)$. Let
  $K$ and $N$ be submodules of finite type in an $A$-module of finite
  type $F$. Then if $K\subset \hat {N}+ \hat {I}(x)\hat F$ when $x\in
  S$, we have
  \begin{displaymath}
    K_{\pfr}\subseteq N_{\pfr}.
  \end{displaymath}
\end{lemma}
\begin{pf}
Letting $\bar K$ be the image of $K$ in $G:=F/N$, we have to prove that
$\bar K_\pfr =0$. 

a) There exists an element $g\in A$ with the following properties: 
\begin{enumerate}\roman{enumi}
\item $g \notin \pfr = (h)$;
\item For any point $x\in S\setminus g^{-1}(0)$ there exists an ideal
  $I(x)\subset A$ such that $h\notin I(x)$, $I(x)+ (h)=\mf_x$, and
  $G_g\cap (\hat I(x)G_g) \subseteq I(x)G_g$.
\end{enumerate}
The proof is similar to the proof of the previous lemma and we use a
similar notation. There exists $f\in A$ such that $G_f$ is free over
$A_f$. Setting $I(x)_f= A_f \cap \hat I(x)_f \subset A_f$, then $G_f
\cap \hat I(x) G_f = I(x)_fG_f$.  If $f\notin (h)$ we are done letting
$g=f$, so now suppose that $f\in (h)$. Then there exists $g\in
A\setminus (h)$ such that $G^{(f)}_{g}= \Ker (G_{g}\to G_{fg})$ is a
successive extension of modules of the type $A/h$. Now choose any
ideal $I(x)\subset A$ such that $\hat I(x)\cap A\subset I(x)$, $I(x)+
(h)=\mf_x$, and $h\notin I(x)$. Then
\begin{eqnarray*}
G^{(f)}_g\cap \hat I (x) G^{(f)}_g  \subset I(x) G^{(f)}_g,\\  
G_{fg} \cap \hat I(x) G_{fg}  \subset I(x)G_{fg}
\end{eqnarray*}
implying $(\hat I(x)G_g) \cap G_{g} \subset I(x)G_g$.

b) We have  by  Artin-Rees' lemma $G=F/N \subset \hat F / \hat N$. Since
$\bar K \subseteq \hat {I}(x)\hat F/\hat N$,  by a) there exists $g\in 
A$ satisfying $(1-2)$ above, hence
\begin{displaymath}
  \bar K_g \subset  G_g \bigcap \hat I (x)G \subset I(x) G , \quad  x\in S,
\end{displaymath}
hence $\bar K_g \subset \cap_{x\in S}I(x)G_g$.  Since $\bar K_{\pfr} =
 (\bar K_g)_{\pfr}$, the proof can be finished with
\Lemma{nullstellen}.
\end{pf}

  \begin{pfof}{\Lemma{fibre:0}}
    Clearly one can assume that $\Omega = X$ and that there exists an
    $\Oc_X$-coherent submodule $M_X^0\subset M_X$ such that
    $M_X=\Dc(\g_X)M_X^0$. We need to prove that there exists a
    (closed) point $y\in S$ and an integer $n$ such that
    $\delta^nM^0_y \subseteq \Dc^{n-1}(\g'_y)M^0_y$.
    
    a) $S$ is not a countable union of nowhere dense subsets: When $X$
    is a variety the normalization theorem shows that we can reduce to
    either   $S$ being discrete or  affine $k$-space and then the
    assertion is obvious as the field $k$ is uncountable.  When $S$
    is a complex analytic space it is a locally compact Hausdorff
    space, so the statement follows from Baire's
    theorem.\ 
    
    b) Let $x$ be a point in $S$. By Artin-Rees lemma $\hat
    {\Oc}_x$ is faithfully flat over $\Oc_x$, so 
    \begin{eqnarray*}
      \Dc^{n_x-1}(\hat
    \g'_x ) \hat M^0_x = \hat {\Oc}_x\otimes_{\Oc_x}\Dc^{n_x-1}(\g'_x)
    M^0_x, \\
\hat I(x) \hat M^{m_x}_x = \hat {I}(x)\otimes_{\Oc_x}M^{m_x}_x,
    \end{eqnarray*}
    and by $(3)$ there exist positive integers $n_x, m_x$ such that
\begin{displaymath}
   \delta^{n_x} M_x \subset \Dc^{n_x-1}(\hat \g'_x)\hat M^0_x + \hat I(x) \hat M^{m_x}_x.
\end{displaymath}
Setting
    \begin{displaymath}
      S_{n,m}= \{x\in S : \delta^n M_x \subset \Dc^{n-1}(\hat {\g}'_x)
      M^0_x + \hat I(x) \hat {M}^{m}_x\}
    \end{displaymath}
    we thus have $S= \cup_{n,m}S_{n,m}$. By a) there exist integers
    $n_0,m_0$ such that the closure of $S_0:=S_{n_0,m_0}$ contains an
    interior point. Shrinking $X$ if necessary we can therefore assume
    that $X$ is either a compact Stein set or affine, so that $A=
    \Oc_X(X)$ is noetherian, and that $S_0$
    is dense in $S$.
    
    Put $M^0= M^0(X)$, $K= \delta^{n_0}M^0$, $N= \Dc^{n_0-1}(\g')M^0$
    and $F = \Dc^{m_0}(\g) M^0(X)$; these are $A$-modules of finite
    type where $N$ and $F$ generate the coherent $\Oc_X$-modules
    $\Dc^{n_0-1}(\g') M^0_X$ and $ M^{m_0}_X$.
    
    Since $\hat A_x$ is faithfully flat over $A$ (Artin-Rees) we have
\begin{displaymath}
  K\subset \hat A_x \otimes_A N + \hat I(x)\otimes_A F, \quad x\in S_0,
\end{displaymath}
so by \Lemma{inclusion} $K_\pfr \subset N_\pfr$. By coherence this
inclusion holds in a neighbourhood of $\pfr \in \Spec A$, and since
the closed points are dense in $S= V(\pfr)$ there exists a closed point $y\in S$
such that $K_y \subset N_y$.
\end{pfof}

\begin{pfof} {\Theorem{curve-test}}  $(1)\Rightarrow (3)$,
  $(1)\Rightarrow (4)$: This follows from the non-derived version of
  \Proposition{preservation}; it is not necessary that $\g_X$ be
  locally free.  $(3)\Rightarrow (2)$: o.k. It remains to prove $(2)
  \Rightarrow (1)$ and $(4)\Rightarrow (1)$: When $\dim X =1$ it
  follows that $X=X_0$, hence it is a non-singular curve, so we can
  clearly assume that $\dim X \geq 2$.
  
  a) Since the assertion is local in $X$ it will be clear that all
  locally defined sections can be assumed to be globally defined on
  $X$.  One can clearly assume that $S$ and hence $S_0$ are connected.
  As $M$ is torsion free it suffices to prove that $M_x$ is smooth
  along $\g'_x$ at one point $x$ in $S_0$ \Cor{point-cor}, and by
  \Theorem{fibre:1} it suffices to prove that $M_x$ is smooth along
  each element in a set of generators of $\g'_x$. One can thus, by
  shrinking $X$ if necessary, assume that $S=S_0$ and that the ideal
  of $S$ is a principal prime ideal $\pfr \subset \Oc_X$ generated by
  a function $h\in \Oc_X(X)$. Let $(\hat h)_y \subset \hat {\Oc}_y$
  denote the principal ideal of $h\in \Oc_y$ in the completion $\hat
  {\Oc}_y$ at a point $y\in X$.
  
  b) The assumption is that $\g'_X$ has globally defined
  generators 
  \begin{displaymath}
    \delta(1), \delta(2),\dots , \delta(k),
  \end{displaymath}
  which in the case $(4)$ satisfy the condition (G). By removing
  subsets of $\codim \geq 2$ from $X$ \Cor{pure} we on the one hand
  ensure that $X=X_0$, so $X$ is non-singular, and on the other ensure
  that $\alpha (\delta(i))\in T_X$ is of the form $g_i \partial_i$
  where ${\partial_i}_y)\notin \mf_y T_y$ for every point $y\in X$ (it
  is non-critical).  Moreover we have assumed that we can choose the
  $\delta (i)$ such that $\alpha (\delta (i)_y) (h)\notin (\hat h)_y$.
  Therefore, by a formal version of Frobenius' theorem, each generator
  $\delta (i)$ has for each point $y\in S$ an analytic integral curve
  passing through $y$, defined by an ideal $\hat I^{(i)}(y)\subset
  \hat {\Oc}_y$, satisfying:
\begin{eqnarray*}
\alpha(\delta(i)_y)(\hat I^{(i)}(y)) &\subseteq& \hat I^{(i)}(y);\\
h &\notin& \hat {I}^{(i)}(y);\\
\hat {I}^{(i)}(y)+ (\hat h)_y  &\subset& \hat {\mf}_y.  
\end{eqnarray*}
Therefore $(2)$ or $(4)$ implies, by \Lemma{fibre:0}, that there
exists a point $x_i\in S$ such that $M_{x_i}$ is smooth along
$\delta(i)_{x_i}$ in $\g'_{x_i}$.

c) There exists an open neighbourhood $U_{1}$ of $x_1$ such that for
every point $y\in U_{1}$ the stalk $M_y$ is smooth along $\delta(1)_y$
in $\g'_y$ \Prop{point}.  By b), with $X= U_1$, there exist a point
$x_2\in U_1 \cap S_0$ such that $M_{x_2}$ is smooth along
$\delta(2)_{x_j}$ in $\g'_{x_2}$.  Again there exists open
neighbourhood $U_{2}\subset U_{1}$ of $x_2$ such that $y\in U_2$ then
$M_y$ is smooth along $\delta(2)_y$ in $\g'_{y}$. Iterating we get an
open subset $U = U_k \subseteq U_{k-1}\subseteq \cdots \subseteq U_1$
such that $M_x$ is smooth along each generator $\delta(i)_x \in \g'_x$
when $x\in U\cap S_0$.
\end{pfof}

\subsection{Regular singularities} 
\label{reg-sing}
{\it Unless otherwise stated we assume that $X$ is a complex analytic
  or algebraic manifold and $\g_X$ is a locally free transitive Lie
  algebroid. When $\pi : Y \to X$ is a morphism of spaces we tacitly
  assume that $\g_Y$  denotes $\pi^+(\g_X)$.}  
\medskip

First let $X$ be any space, $\g_X\in \Lie_X$, $M\in \coh (\g_X)$ be
pure with support $V$, and let $J$ be the ideal of $V$. Let $S\subset
X$ be a subspace whose ideal is $I$. We say that a pure $\g_X$-module
$M$ has {\it regular singularities along $S$} \/ if $M$ is point-wise
smooth along the sub-algebroid
\begin{displaymath}
   \g_X(J) \cap \g_X(I) \subset \g_X.
\end{displaymath}
Consider a decomposition
\begin{displaymath}
  S= D \cup D'
\end{displaymath}
where $D$ is the union of the irreducible components $S_i\subseteq S$
such that $S_i \cap V$  has pure codimension $1$ in $V$, and $D'$ is the
union of the irreducible components $S_j$ such that $S_j \cap V$ has
codimension $\geq 2$ in $V$; let $I_1$ be the ideal of $D$ and $I_2$
the ideal of $D'$.  Now since $\g_X(I)_{|X\setminus D'} =
\g_X(I_1)_{|X\setminus D'}$ it follows by \Corollary{purestcor} that
$M$ is smooth along $\g_X(I_1)\cap \g_X(J)$ if and only if $M$ is
smooth along $\g_X(I)\cap \g_X(J)$.  We may therefore in practice
assume that $S=D$ when checking if $M$ has regular singularities along 
$S$.

If $S$ is a subspace of a space $S'$ it may not be obvious that $M$
has regular singularities along $S'$ when it has regular singularities
along $S$.  That this is so follows from the following lemma, which
also implies that $M$ has regular singularities along $S$ if it has
regular singularities along each component $S_i $ in an irreducible
decomposition $S= \cup S_i$.

\begin{lemma}\label{no-inclu}($X$ is a space)
  Let $I\subset \Oc_X$ be a proper ideal and $\sqrt{I}$ its
  radical.  Then:
  \begin{enumerate}
  \item If $I=I_1 \cap I_2  \cdots \cap I_r$ is a minimal primary decomposition we have 
  \begin{displaymath}
    T_X(I)= \T_X(I_1)\cap T_X(I_2) \cap \cdots \cap T_X(I_r);
  \end{displaymath}
\item   $T_X(I)\subset T_X(\sqrt{I})$.
  \end{enumerate}
\end{lemma}

\begin{remark}   One also has $T_X(\sqrt {I})=
  T_X((\sqrt{I})^n)$ for any positive integer $n$, but in general
  $T_X(\sqrt{I})\nsubseteq T_X(I)$. For example, let
  $I=(x+y,y^2)\subset k[x,y]$, then $T(\sqrt{I})= T((x,y)) =
  k[x,y]y\partial_y + k[x,y]y\partial_y$, but $x\partial_x (x+y)= x  \notin I $.
\end{remark}
\begin{pf}
  (1): It is evident that $T_X(I_1)\cap \cdots \cap T_X(I_r)\subset
  T_X(I)$. Now let $\partial \in T_X(I)$ and $x_1\in I_1$. Since we
  are considering a minimal primary decomposition there exists an
  element $y\in I_2 \cap \cdots I_r \setminus \sqrt{I_1}$. Then
  \begin{displaymath}
    y \partial (x)= \partial(yx) - x\partial(y)\in I_1
  \end{displaymath}
and since $I_1$ is $\sqrt{I_1}$-primary, $\partial(x)\in I_1$. In the
same way  $\partial (I_i)\subset I_i$, $i=1, \ldots r$, implying the assertion.

(2): First assume that $P= \sqrt{I}$ is prime and let $\partial \in
T(I)$, $x\in P$, and $x^n\in I$. If $n=1$ we have $\partial(x)\in
I\subset P$ so assume that $x\notin I$. Then $nx^{n-1} \partial (x)\in
I$. Supposing that $\partial(x)\notin P$, since $I$ is $P$-primary we
get $nx^{n-1}\in I$, and iterating gives $n! x \in I$. Since the
characteristic is $0$ we reach a contradiction. Thus $\partial (x)\in
P$. This shows that $T(I)\subseteq T(P)$. Hence $T_X(I_i)\subset
T_X(\sqrt{I_i})$.  By (1)
\begin{eqnarray*}
  T_X(I)= T_X(I_1) \cap \cdots \cap T_X(I_r)\\
T_X(\sqrt{I}) = T_X(\sqrt{I_1})\cap \cdots \cap T_X(\sqrt{I_r});
\end{eqnarray*}
hence $T_X(I)\subset T_X(\sqrt I)$ for any ideal $I$.
\end{pf}

Now let $M$ be any coherent $\g_X$-module and
\begin{equation}\label{filt-supp}
  0 = M_{n}  \subset M_{n-1} \subset \cdots   \subset M_0 =M
\end{equation}
be its Krull filtration (see \Section{cohprolong}), where successive quotients
$L_i= M_{i}/M_{i+1}$ are pure.

\begin{definition}\label{regsing-def}
  A coherent $\g_X$-module $M$ has regular singularities if each  pure
  component $L_i$ has regular singularities as above.   Let $\coh_{rs}(\g_X)$
  be the category of such modules.
\end{definition}

\begin{remarks}
  \begin{enumerate}
\item \label{discrete-case} If $M\in \coh (\g_X)$ and $V= \supp M$ is
  discrete, then $M\in \coh_{rs}(\g_X)$. To see this, let $x\in \supp
  M$ and $M^0_x \subset M_x$ be an $\Oc_x$-submodule of finite type.
  Then $\mf_x^n M^0_x =0$ for some big integer $n$. If $\phi \in
  \mf_x$, $\mu \in M^0_x$, and $\delta \in \g(\mf_x)$, then since
  $\alpha (\delta)(\phi)\in \mf_x$, so $\phi^n \delta \mu =0$, it
  follows that $\Dc(\g(\mf_x))M^0_x$ is of finite type over $\Oc_x$.
  Thus $M_x$ is smooth along $\g_x(\mf_x)$.
  
\item By \Proposition{reg-D} below a torsion free coherent
  $T_X$-module $M$ on a space $X$ has regular singularities if and
  only if $M$ has regular singularities in the sense of \cite[Th
  4.1]{deligne:eq}.
  \end{enumerate}
\end{remarks}

\begin{prop}\label{reg-inv}
 Consider a short  exact sequence of coherent $\g_X$-modules
  \begin{displaymath}
    0\to K \to M \to N \to 0.
  \end{displaymath}
Then $M$ has regular singularities  if and only if $K$
and $N$ have regular singularities. 
\end{prop}
\begin{pf}  Consider the Krull filtration in (\ref{filt-supp}).

  a) If $M\in \coh_{rs}(\g_X)$ then $K,N\in \coh_{rs}(\g_X)$: We have
  a short exact sequence
\begin{equation}\label{eq:pureseq}
  0\to \frac {K\cap M_i}{K\cap M_{i+1}} \to \frac {M_i}{M_{i+1}}\to
  \frac {\bar M_i}{\bar M_{i+1}} \to 0
\end{equation}
where $\bar M_i$ is the image of $M_i$ in $N$. Here $K\cap M_i/K\cap
M_{i+1}$ is $0$ or a maximal pure sub-quotient of $K$, and all maximal
pure sub-quotients of $K$ occur in this way, but $\bar M_i/\bar
M_{i+1}$ need not be pure. Now if $J_i \subset I_i$ are radical
ideals such that $M_i/M_{i+1}$ is smooth along $\g_X^i=\g_X(J_i)\cap
\g_X(I_i)$, $i=0,\ldots ,m$, then $K\cap M_i/K\cap M_{i+1}$ and $\bar
M_i/\bar M_{i+1}$ also are smooth along $\g^i_X$ \Prop{extension}.
Hence $K\in \coh_{rs}(\g_X)$, and every pure subquotient of $\bar M_i
/\bar M_{i+1}$ is smooth along $\g^i_X$ \Prop{extension}, so $N \in
\coh_{rs}(\g_X)$.
  
b) If $K,N \in \coh_{rs}(\g_X)$ then $M\in \coh_{rs}(\g_X)$: Consider
the short exact sequence (\ref{eq:pureseq}). As $\bar M_i \subset N$
and $M_i\cap K \subset K$, by a) $\bar M_i, M_{i}\cap K\in
\coh_{rs}(\g_X)$, and again by a) $(K\cap M_i)/(K\cap M_{i+1}), \bar
M_i/ \bar M_{i+1}\in \coh_{rs}(\g_X)$.  It therefore suffices to prove
that $M\in \coh_{rs}(\g_X)$ when $M$ is pure. Then $K$ is also pure,
and by \Corollary{purestcor} we can assume also that $N$ is pure. Put
$V_K=\supp K$, $V_M=\supp M$ and $V_N=\supp N$, and let $J_K, J_M, J_N
\subset \Oc_X$ be their ideals.  Let $I_K, I_N\subset \Oc_X$ be ideals
such that $K$ is smooth along $\g_X(J_K)\cap \g_X(I_K)$ and $N$ is
smooth along $\g_X(J_N)\cap \g_X(I_N)$; let $S$ and $T$ be the spaces
of $I_K$ and $I_N$. As $M$ is pure we have $\dim V_K = \dim V_M$, and
we can assume that $\dim V_M - 1 \leq \dim V_N \leq \dim V_M$
\Cor{purestcor}. So we have two cases.
  
Suppose first that $\dim V_N = \dim V_M-1$. Then $\codim_{V_M} T \geq
2$, hence we may then assume that $N$ is smooth along $\g_X(J_N)$
\Cor{purestcor}. Moreover, \Corollary{purestcor} also implies that one
may assume that all irreducible components of $S$ and $V_N$ have the
same dimension. Letting $I_M$ be the ideal of $S\cup V_N$ we claim
that $M$ is smooth along
  \begin{displaymath}
    \g'_X:=\g_X(J_M)\cap\g_X(I_M),
  \end{displaymath}
  so $M\in \coh_{rs}(\g_X)$. By \Proposition{extension} this follows
  if $N$ and $K$ are smooth along $\g'_X$.  We have assumed that all
  irreducible components of $S$ and $V_N$ have the same dimension;
  hence the prime components of $I_K\cap J_N$ is the union of the
  prime components of $I_K$ and $J_N$; hence \Lemma{no-inclu} implies
  \begin{displaymath}
  \g'_X \subset \g_X(I_M)= \g_X(I_K \cap J_N)  \subset \g_X(J_N).
  \end{displaymath}
  Thus $N$ is smooth along $\g'_X$.  Since $V_K \subset V_M$ and $\dim
  V_K = \dim V_M$ we may assume that the prime components of $J_K$ is
  a part of the prime components of $J_M$ \Cor{purestcor} so
  $\g_X(J_M)\subset \g_X(J_K)$ \Lem{no-inclu}; in the same way we get
  $\g_X(I_M)\subset \g_X(I_K)$. Hence
   \begin{displaymath}
     \g'_X \subseteq \g_X(J_K)\cap \g_X(I_K \cap J_N)  \subseteq \g_X(J_K)\cap \g_X(I_K);
   \end{displaymath}
    and therefore $K$ is    smooth along $\g'_X$.   
   
    Now suppose that $\dim V_N = \dim V_M$. Then we may assume that $
    S$ and $T$ have pure codimension $1$ in $V_M$ \Cor{purestcor}; let
    $I_M$ be the ideal of $S\cup T$.  Using a very similar argument as
    above, based on \Lemma{no-inclu}, one sees that $K$ and $N$ are
    smooth along $\g_X(J_M)\cap \g_X(I_M)$; hence $M$ is smooth along
    $\g_X(J_M)\cap \g_X(I_M)$ \Prop{extension}.  
\end{pf}

The subcategories of $D^b_{\coh}(\g_X)$ that are of concern in this
section are denoted $D^b_{rs}(\g_X)$, $D^b_{crs}(\g_X)$, and
$D^b_{\hol}(\g_X)$.  They are defined as follows: $D^b_{rs}(\g_X)$ is
the derived category of bounded complexes $M^\bullet$ of
$\g_X$-modules whose homology $H^\bullet(M^\bullet)\in
\coh_{rs}(\g_X)$.  It follows from \Proposition{reg-inv} that
$D_{rs}^b(\g_X)$ is a triangulated sub-category of $D_{\coh}^b(\g_X)$;
i.e.\ if two vertices of a distinguished triangle in
$D_{\coh}^b(\g_X)$ belong to $D_{rs}^b(\g_X)$, then all three vertices
belong to $D_{rs}^b(\g_X)$.  $D^b_{crs}(\g_X)$ is the sub-category of
complexes $M^\bullet$ such that for any curve $\pi : C \to X$ one has
$\pi^!(M^\bullet)\in D^b_{rs}(\g_C)$; we then say that $M^\bullet$ is
{\it curve regular}.  Since $D_{rs}^b(\g_C)$ is a triangulated
category and $\pi^!$ preserves triangles, it follows that
$D^b_{crs}(\g_X)$ is a triangulated category.  A coherent
$\g_X$-module $M$ is holonomic when $\dim \SSo(M) = \dim X$
(\ref{diffoperators}).  $D^b_{\hol}(\g_X)$ is the derived category of
bounded complexes with holonomic homology modules; it is easy to see
that sub-quotients of holonomic modules are holonomic, implying that
$D^b_{\hol}(\g_X)$ is a triangulated category. We refer to
\cite{borel:Dmod} and \cite{bjork:analD} for the theory of holonomic
$T_X$-modules on a non-singular space $X$.

\begin{remark}\label{complication}
  Let $\phi : (Y,\g_Y) \to (X, \g_X)$ be a morphism of Lie algebroids
  on non-singular spaces, where  $\g_Y=
  \phi^+(\g_X)$. Letting $I\subset \Oc_X$ be the  ideal of a divisor $S$ on
  $X$ and $J\subset \Oc_Y$ the ideal of $\pi^{-1}(S)$    we have a canonical morphism 
  \begin{displaymath}
   \phi^+(\g_X(I))  \to  \g_Y(J),
  \end{displaymath}
  but this need not be surjective when $S$ is singular and $\dop \phi
  : T_Y \to \phi^*(T_X)$ is not surjective.  For instance, let $I=
  (x_1^3+x_2^3 + x_3^3)\subset A:=\Cb\{x_1, x_2 , x_3\}$ be the ideal
  of the germ of the cubic cone in $(\Cb^3,0)$. Then the $A$-module of
  derivations that preserve $I$ is $T_{(\Cb^3,0)}(I) = AE$ where $E=
  x_1\partial_{x_1}+ x_2\partial_{x_2}+ x_3 \partial_{x_3}$
  \cite{bernstein-gege:cubic}.  Therefore, if $Y \subset (\Cb^3 ,0)$
  is a germ of an irreducible curve in $(\Cb^3,0)$ with ideal $L
  \subset A$ such that $E\notin T_{\Cb^3,0}(L)$, then one easily sees
  that $T_{\Cb^3,0}(I)\cap T_{\Cb^3,0}(L) = LE$. Hence, if $\phi : Y
  \to (\Cb^3,0)$ is the inclusion morphism then the induced
  homomorphism $\phi^+(T_{\Cb^3,0}(I)) =
  \phi^*(T_{\Cb^3,0}(I))\times_{\phi^*(T_{\Cb^3,0})}T_Y\to T_Y$ is
  $0$.
\end{remark}

When $X$ is non-singular we will generalize the proposition below
considerably in \Theorem{equivalent}.

\begin{proposition}\label{reg-D}  ($X$ is a space of the same type as in 
  \Theorem{curve-test}) Let $M$ be a holonomic $\Dc_X$-module without
  $\Oc_X$-torsion and $I$ be the ideal of its strong support (see
  \Lemma{strong-var}).  The following are equivalent:
  \begin{enumerate}
  \item $\phi^!(M)$ is  regular along $\phi^+(T_X(I))$ when $\phi :
    C\to X$ is a curve (when $X$ is singular we only consider
    the non-derived inverse image);
\item $M$ is smooth  along $T_X(I)$.
  \end{enumerate}
\end{proposition}
\begin{remark}
  \cite[VII, Props. 11.6-7 and Cor. 11.8]{borel:Dmod} are attained by
  employing GAGA\@, but they also follow from \Proposition{reg-D} and
  \Theorem{curve-test} without GAGA.
\end{remark}

\begin{pf} 
  $(2)\Rightarrow (1)$: Let $\phi : C \to X$ be a curve. By
  \Theorem{curve-test} we know that the non-derived inverse image
  $\phi^!(M)\in \coh_s(T_C)$.  In the non-singular case the
  assertion follows from \Proposition{preservation}, but we will see
  below  that we actually have $\phi^!(M)\in
  D^b_{rs}(T_C)$ \Th{inv-rs}.
  
  $(1) \Rightarrow (2)$: Since $M$ is holonomic, $\SSo M \subset T^*X$
  is a conic algebraic (analytic) Lagrangian subset; hence, $M$ being
  torsion free, $S:=\ssupp M \subset X$ is a divisor
  \Defn{strong-supp}. As the space $X$ is non-singular in codimension
  $1$ we may assume that $S$ and $X$ are non-singular spaces
  \Cor{pure}.  Let $x$ be a point in $S$ and $(t=x_1, x_2, \dots ,
  x_n)$ be a regular system of parameters of the local ring $\Oc_x$
  such that $t=0$ is a local equation for $S$ at $x$. Since $\Oc_x$ is
  either isomorphic to the ring of convergent power series in the
  given parameters, or a localisation of a ring of finite type, there
  exist derivations $\partial_t = \partial_{x_1}, \partial_{x_2},
  \dots , \partial_{x_n}\in \Der_k(\Oc_x)$ such that
  $\partial_{x_i}(x_j)= \delta_{ij}$ (the Kronecker delta). As $\SSo
  M$ is an algebraic (analytic) Lagrangian set, each irreducible
  component of the nonsingular locus of $\dSSo M$
  (see~\ref{diffoperators}) belongs to the conormal set of $S$;
  therefore $M_x$ is smooth along $\partial_{x_2}, \dots
  \partial_{x_n}$. Now $(1)$ and \Theorem{curve-test} implies that
  $M_x$ is smooth along $t\partial_t$; hence $M_{x}$ is point-wise
  smooth along $\Oc_x t \partial_t + \Oc_x\partial_{x_2}+ \cdots +
  \Oc_x \partial_{x_n}= T_x(I_S)$, proving $(2)$.
\end{pf}

The following result is due to Bernstein.
\begin{prop}\label{bernstein1}($X$ is a non-singular space)
  Let $j : Y \to X$ be an open embedding of spaces such that $S:=
  X\setminus j(Y)$ is a  subspace.
  Let $M^\bullet \in D^b_{\coh}(\g_X)$ and consider the distinguished
  triangle
  \begin{equation}\label{loc-seq}
    R\Gamma_{[S]} (M^\bullet) \to M^\bullet \to M^\bullet (*S)\to .
  \end{equation}
Then the following are equivalent:
\begin{enumerate}
\item $M\in D^b_{\hol}(\g_X)$;
\item $R\Gamma_{[S]}(M^\bullet)$ and $M^\bullet(*S)$ belong to
  $D^b_{\hol}(\g_X)$.
\end{enumerate}
\end{prop}
Using a Mayer-Vietoris argument, the main case in the proof is
$(1)\Rightarrow (2)$ when $S$ is a divisor and this follows from the
existence of  functional equations  containing a Bernstein-Sato
polynomial; see e.g.  \cite[Th.  3.2.13]{bjork:analD}.

\begin{prop}\label{holonomic}
  $D^b_{rs}(\g_X)\subset D^b_{hol}(\g_X)$.
\end{prop}
\begin{pf}
  $D^b_{rs}(\g_X)$ and $D^b_{hol}(\g_X)$ are generated by its homology
  objects, so it suffices to see that a pure $\g_X$-module $M$ with
  regular singularities is holonomic. Let $J$ be the ideal of the
  support $V$ of $M$ and $I \subset \Oc_X$ be an ideal such that $M$
  is smooth along $\g_X(J)\cap \g_X(I)$. Let $S$ be the space of $I$.
  Since $\Dc(\g_X)$ is noetherian (\ref{noetherianrings}) the union
  $N$ of all holonomic sub-modules of $M$ is holonomic. Then clearly
  $N|_{X\setminus S} = M|_{X\setminus S}$, so the quotient $K=M/N$ is
  a coherent $\g_X$-module whose support belongs to $S\subset V$;
  hence $H^0(K(*S)^\bullet)=0$, so by the first part of the long exact
  sequence induced by the distinguished triangle $N(*S)^\bullet \to
  M(*S)^\bullet \to K(*S)^\bullet \to $ one gets $H^0(N(*S))=
  H^0(M(*S))$; by \Proposition{bernstein1} $H^0(N(*S)^\bullet)$ is
  holonomic; hence $M\subset H^0(M(*S)^\bullet)$ is holonomic.
\end{pf}

We next have a version of Bernstein's theorem for regular holonomic
$\Dc$-modules; see \cite[Th. 5.4.1]{bjork:analD}.
\begin{thm}\label{bernstein}
  Let $j : Y \to X$ be an open embedding such that $S:= X\setminus
  j(Y)$ is a a subspace. Let
  $M^\bullet \in D^b_{\coh}(\g_X)$ and consider the distinguished
  triangle (\ref{loc-seq}).  
  Then the following are equivalent:
\begin{enumerate}
\item $M^\bullet \in D^b_{rs}(\g_X)$.
\item $R\Gamma_{[S]}(M^\bullet)$ and $M(*S)^\bullet$ belong to
  $D^b_{rs}(\g_X)$.
\end{enumerate}
\end{thm}

\begin{pf} $(2) \Rightarrow (1)$ follows since $D^b_{rs}(\g_X)$ is a
  triangulated category. 
  
  To prove the converse one may assume that $M^\bullet$ is a single
  degree complex, consisting of a pure $\g_X$-module $M$ with regular
  singularities; let $K$ be the ideal of $V=\supp M$ and $I$ an ideal
  such that $M$ is smooth along $\g_X(K)\cap \g_X(I)$.  Letting $T$ be 
  the space of $I$, by the  discussion before \Lemma{no-inclu} we can assume that every
  irreducible component   of $ V\cap  T  $ has pure codimension $1$ in $V$.
  
  By \Propositions{bernstein1}{holonomic} 
  \begin{displaymath}
    R\Gamma_{[S]}(M)^\bullet, M(*S)^\bullet \in D^b_{\coh}(\g_Y)
  \end{displaymath}
  so it remains to prove that $R\Gamma_{[S]}(M^\bullet)^\bullet$ and
  $M(*S)^\bullet$ have regular singularities.
  
  a) We first prove that $M(*S)\in \coh_{rs}(\g_X)$ when $S$ is a
  divisor, locally defined by a function $f\in \Oc_X$. So in this case
  $M(*S)^\bullet$ is a single degree complex given by a coherent
  $\g_X$-module $M(*S)$.  If $S\cap V = \emptyset$ then $M(*S)=0$, so
  we assume $S\cap V \neq \emptyset$, and then the discussion before
  \Lemma{no-inclu} shows that one may assume that $S\cap V$ is of pure
  codimension $1$ in $V$. Then $M(*S)$ is again a pure $\g_X$-module,
  and it suffices to see that it is smooth along 
  \begin{displaymath}
    \g'_X= \g_X(K)\cap \g_X(I\cap J),
  \end{displaymath}
  where $J= (f)$ is the ideal of $S$.  Now there exists (locally) a
  coherent $\Oc_X$-submodule $M_0\subset M$ such that $\Dc(\g_X)M_0
  =M$ and $\g_X(K)\cap \g_X(I)\cdot M_0 \subset M_0$ and an integer
  $k$ such that $\Dc(\g_X)M_0f^{-k}= M(*S)$. It suffices now to see
  that $\g'_X \cdot M_0f^{-k}\subseteq M_0f^{-k}$. Since the prime
  components of $I\cap J$ is the union of the prime components of $I$
  and the prime components of $J$, \Lemma{no-inclu} implies
  \begin{displaymath}
    \g'_X = \g_X(K)\cap \g_X(I)\cap \g_X(J) \subset \g_X(J).
  \end{displaymath}
Therefore, if $\delta \in \g'_X$, we have  
\begin{eqnarray*}
  \delta \cdot M_0 f^{-k}&\subseteq& ( \delta \cdot M_0)f^{-k} -k \alpha(f)M_0
  f^{-k-1}\\
&\subseteq& (\delta\cdot M_0)f^{-k} -kJM_0 f^{-k-1}\subseteq M_0f^{-k},
\end{eqnarray*}
  implying the assertion.\ 

  b) The general case. Let $(f_1, \dots , f_l)$ be local generators of
  the ideal $I$ of $S$. Letting $S_i$ be the locally defined divisor
  that is defined by $f_i$, by a) $\oplus_i M(*S_i)\in
  D^b_{rs}(\g_X)$, and inductively $M(*\cup_i S_i)\in D^b_{rs}(\g_X)$.
  Therefore $\oplus_i R\Gamma_{[S_i]}(M)\in D^b_{rs}(\g_X)$ and
  $R\Gamma_{[\cup_i S_i]}(M)\in D^b_{rs}(\g_X)$ (\ref{loc-seq}), so using the
  distinguished triangle of   Mayer-Vietoris type
  \begin{displaymath}
    R\Gamma_{[S]}(M) \to \oplus R\Gamma_{[S_i]}(M) \to R\Gamma_{[\cup_i  S_i]}(M) \to 
  \end{displaymath}
  this implies that the third vertex $R\Gamma_{[S]}(M)\in
  D^b_{rs}(\g_X)$; hence the third vertex $M(*S)^\bullet$ in the
  distinguished triangle (\ref{loc-seq}) also belongs to
  $D^b_{rs}(M)$.
\end{pf}

\begin{thm}\label{inv-rs}
  Let $\pi : Y\to X$ be a morphism of non-singular spaces. Then
  $\pi^!$ is a functor from $D^b_{rs}(\g_X)$  to $D^b_{rs}(\g_Y)$.
\end{thm}

Note that \Proposition{preservation} does not directly imply
\Theorem{inv-rs}, in the light of \Remark{complication}.
  
\begin{pf}    The triangulated category $D^b_{rs}(\g_X)$ is generated by
  $\coh_{rs}(\g_X)$ and $D^b_{rs}(\g_Y)$ is stable under extensions,
  so it suffices to prove that $\pi^!(M)\in D^b_{rs}(\g_X)$ when $M\in
  \coh_{rs}(\g_X)$.  Moreover, by induction over the length of the
  filtration by supports of $M$ one can assume that $M$ is pure. Let
  $I$ be the ideal of $\supp M$ and $J$ an ideal such that $M$ is
  smooth along $\g'_X:=\g_X(J)\cap \g_X(I)$. Now $\pi$ can be
  factorized into a closed embedding and a projection, and by
  \Proposition{composed} it suffices to treat these cases separately.
  
  A projection $\pi$ is flat, hence
  \begin{eqnarray*}
    \pi^!(M) &=& \pi^*(\Dc(\g_X))\otimes^L_{\pi^{-1}(\Dc(\g_X))}\pi^{-1}(M)\\ &=&
    \Oc_Y \otimes^L_{\pi^{-1}(\Oc_X)}\pi^{-1}(M) \\&=&
  \Oc_Y\otimes_{\pi^{-1}(\Oc_X)}\pi^{-1}(M)
  \end{eqnarray*}
  and this is a pure coherent $\g_Y$-module. Letting $I_1, J_1 \subset
  \Oc_Y$ be the ideals that are generated by $I$ and $J$ one sees that
  $\pi^!(M)$ is smooth along $\g'_Y = \pi^+(\g'_X)= \g_Y(\pi^*(I))\cap
  \g_Y(\pi^*(J))$.
  
  Now assume that $\pi$ is a closed embedding.  By
  \Theorems{kashiwara}{bernstein} $\pi_+\pi^!(M)\in D^b_{rs}(\g_X)$.
  So the proof will be complete if one proves: if $N\in \coh (\g_Y)$
  and $\pi_+(N)\in \coh_{rs}(\g_X)$, then $N\in \coh_{rs}(\g_Y)$. To
  prove this one may again assume that $N$ is pure; then $\pi_+(N)$ is
  also pure with support in $Y\subset X$. Suppose that $\pi_+(N)$ is
  smooth along $\g_X(I)\cap \g_X(J)$ where $J$ is the ideal of
  $V:=\supp \pi_+(M)$ and $I$ is the a ideal of a subspace $S \subset
  V\subset X$. Let $L\subset \Oc_X$ be the ideal of $\pi(Y)\subset X$.
  By \Theorem{kashiwara}
\begin{displaymath}
  N= \pi^!\pi_+(N)[\dop_{Y,X}] \cong \pi_+(N)^L \subset \pi_+(N),
\end{displaymath}
where $\pi_+(L)^L$ is the $\Oc_X$-submodule of $L$-invariants of
$\pi_+(N)$. Letting $I_1$ and $J_1$ be the ideals of $S\subset V
\subset Y$ in $\Oc_Y= \Oc_X/L$ we need to prove that $N$ is smooth
along $\g_Y(I_1)\cap \g_Y(J_1)$, for then $N\in \coh_{rs}(\g_Y)$.
Since $L \subset J \subset I$ and, $Y$ being non-singular, $T_Y=
T_X(L)/LT_X$, it follows that the canonical mapping $T_X(I)\cap T_X(J)
\cap T_X(L) \to   T_Y(I_1)\cap T_Y(J_1) $ is surjective. Hence  the canonical mapping
\begin{displaymath}
\g_X(I)\cap \g_X(J)\cap \g_X(L) \to \g_Y(I_1)\cap \g_Y(J_1) 
\end{displaymath}
is surjective. Now $ N = \pi_+(N)^L $ is a $\g_X(L)$-module which is
smooth along the sub-algebroid $\g_X(I)\cap\g_X(J) \cap \g_X(L)$. This
readily implies that $N$ is smooth along $\g_Y(I_1)\cap \g_Y(J_1)$.
\end{pf}

Let $(\g_X, \alpha)$ be a transitive locally free Lie algebroid on
$X$; then $\g_Y:= \pi^+(\g_X)$ is also locally free and transitive on
$Y$ \Lem{pullfree}.
\begin{thm}\label{directreg} Let $\pi :Y \to X$ be a  proper morphism  of
  spaces.  Then
  \begin{displaymath}
    \pi_+(D^b_{rs}(\g_Y))\subset D^b_{rs}(\g_X).
\end{displaymath}
  \end{thm}
\begin{pf}  
  We will use the following criterion for pure modules.  
  \begin{description}
  \item [(C)]  Assume that $M$ is a pure $\g_Y$-module.  Let $J$ be the ideal
  of $V$ and $I$ the ideal of a subspace $S\subset V \subset Y$ such
  that $M$ is smooth along $\g'_Y = \g_Y(J)\cap \g_Y(I)$; then clearly
  $\Ker (\g_Y \to T_Y)\subset \g'_Y$.  Putting $\g_{Y/X}= \alpha^{-1}
  (T_Y/X)$ we let $\af_Y$ be the Lie algebroid that is generated by
  $\g'_Y + \g_{Y/X}$. If there exist a subspace $V_1 \subset V$, so
  that $\codim _V V_1 \geq 2$, and radical ideals $J_1, I_1 \subset
  \Oc_X$ such that
  \begin{equation}\label{main:eq}
    \Imo (\pi^*(\g_X(J_1)\cap \g_X(I_1))_y \to \pi^*(\g_X)_y) \subseteq
    \Imo (\af_y \to \pi^*(\g_X)_y)
  \end{equation}
  when $y\in V\setminus V_1$, then $\pi_+(M)\in D^b_{rs}(\g_X)$
  \Cor{direct-smooth-pure}.
  \end{description}

  We now start the proof. Let $M^\bullet \in D^b_{rs}(\g_Y)$ and set $V=
  \supp M^\bullet$. To prove that $\pi_+(M^\bullet)\in D^b_{rs}(\g_X)$
  we may assume that $M^\bullet = M$ is a single degree complex
  consisting of a coherent $\g_Y$-module with regular singularities.
  We factorize $\pi$  as  $\pi = p \circ i$, where $i : Y \to Y\times
  X$ is the graph embedding and $p: Y\times X \to X$ is the projection
  on the second factor.  \

   a) $i_+(M)\in D^b_{rs}(\g_{Y\times X})$: To see this one may assume
  that $M$ is pure and apply (C), letting $J_1$ and $I_1$ be the
  ideals of $i(V)$ and $i(S)$.

  b) Assume that $M$ is torsion free. Let $S\subset Y$ be a divisor
  such that $M$ is smooth along $\g_Y(I')$, where $I'$ is the ideal of
  $S$.  Let $T$ be the critical locus of $\pi$, i.e.\  the locus of
  points $y\in Y$ such that $\pi$ is not submersive; this also is a
  divisor on $Y$. Let $I$ be the ideal of $S\cup T$. Then $M$ is
  smooth along $\g_Y(I)$. Letting $V_1$ be the singular locus of
  $S\cup T$ and $V= Y$ we have $\codim_V V_1 \geq 2$, and it is
  straightforward to check that (\ref{main:eq}) is satisfied; hence
  $\pi_+(M)\in D^b_{rs}(\g_X)$.
  
  c) The proof now is by induction over $n = \dim Y$. When $n =0$ the
  assertion follows from \Remark{discrete-case}, so assume that the
  theorem holds for all proper morphisms $\pi_1 :Y_1 \to X_1$ when
  $\dim Y_1 <n$.
  
  d) First assuming that $\dim V \le n-1$ we prove that $\pi_+(M)\in
  D^b_{rs}(\g_X)$ by induction over $\dim V$. If $\dim V =0$ the
  assertion follows from \Remark{discrete-case}, so assume that $\dim
  V \geq 1$.  Let $r: V' \to Y\times X $ be a resolution of $\supp
  i_+(M) = i(V)$. We then have a distinguished triangle \Th{adj-prop}
\begin{equation}
  r_+r^!(i_+(M)) \to i_+(M) \to C^\bullet \to 
\end{equation}
and applying $p_+$, using $p_+r_+ = (p\circ r)_+$ and $\pi_+= p_+i_+$
\Prop{direct-compose}, we also have the distinguished triangle
\begin{equation}
  (p\circ r)_+r^!(i_+(M)) \to \pi_+(M) \to p_+(C^\bullet) \to.
\end{equation}
By a) $i_+(M)\in D^b_{rs}(\g_{Y\times X})$, hence $r^!(i_+(M))\in
D^b_{rs}(\g_{V'})$ \Th{inv-rs}. Since $\dim V' \leq n-1$, by induction
$(p\circ r)_+r^!(i_+(M))\in D^b_{rs}(\g_{X})$ and $r_+r^!(i_+(M)) \in
D^b_{rs}(\g_{Y\times X})$. It will now follow that $\pi_+(M)\in
D^b_{rs}(\g_X)$ if we prove that $p_+(C^\bullet)\in D^b_{rs}(\g_X)$.

First, $C^\bullet \in D^b_{rs}(\g_{Y\times X})$ since $r_+r^!(i_+(M)),
i_+(M)\in D^b_{rs}(\g_{Y\times X})$. Second, $i$ is a closed embedding
and $\supp C^\bullet \subset i(Y)$ so $C^\bullet =
i_+i^!(C^\bullet)[-\dim X]$ \Th{kashiwara}. Hence
\begin{equation}
  \label{eq:indstep}
  p_+(C^\bullet) = p_+i_+i^!(C^\bullet)[-\dim X]=  \pi_+(i^!(C^\bullet)[-\dim X]) \quad \Prop{direct-compose}.
\end{equation}
Since $r$ is birational to the support of $i_+(M)$ 
\begin{displaymath}
  \dim \supp i^!(C^\bullet)= \dim \supp C^\bullet \leq \dim V -1;
\end{displaymath}
hence by induction and (\ref{eq:indstep}) it follows
that $p_+(C^\bullet)\in D^b_{rs}(\g_X)$.

e) Assume that $\dim V= n$. There exists a divisor $S\subset Y$ so
that in the distinguished triangle
\begin{displaymath}
  R\Gamma_{[S]}(M) \to M \to M(*S) \to
\end{displaymath}
the $\g_Y$-module $M(*S)$ is torsion free. Applying $\pi_+$ we get the
distinguished triangle
\begin{displaymath}
  \pi_+R\Gamma_{[S]}(M)\to \pi_+(M) \to \pi_+(M(*S)) \to.
\end{displaymath}
Now   $R\Gamma_{[S]}(M), M(*S)\in D^b_{rs}(\g_Y)$ \Th{bernstein},
hence by b) $\pi_+(M(*S))\in D^b_{rs}(\g_X)$ and by d) $\pi_+(R\Gamma_{[S]}(M))\in D^b_{rs}(\g_X)$. Therefore the third vertex
$\pi_+(M)\in D^b_{rs}(\g_X)$.
\end{pf}

\begin{remark} Using an  argument similar to b) in the proof  one can 
  prove the theorem for pure modules whose support $V$ is regular in
  codimension $1$ (for example when $V$ is normal) without using
  Hironaka's resolution of singularities. Still, I was unable to avoid
  the general resolution of singularities altogether by instead in d)
  let $r: V' \to Y\times X$ be the normalization of $V$.
\end{remark}
  
\begin{thm}\label{equivalent}
\begin{displaymath}
  D^b_{crs}(\g_X) = D^b_{rs}(\g_X).
\end{displaymath}  
\end{thm}
The proof of \Theorem{equivalent} is based on two lemmas, where the
first is the counterpart for $D^b_{crs}$ of
\Theorem{bernstein}. 
\begin{lemma}\label{div-curve}
  Let $S\subset X$ be a non-singular divisor, where $i: S\to X$ is the
  inclusion.  Consider the distinguished triangle \Th{kashiwara}
  \begin{displaymath}
    i_*i^!(M^\bullet)[\dop_{S,X}] \to M^\bullet \to M^\bullet (*S) \to.
  \end{displaymath}
  Then the following are equivalent:
\begin{enumerate}
\item $M\in D^b_{crs}(\g_X)$;
\item $R\Gamma_{[S]}(M^\bullet)$ and $M^\bullet(*S)$ belong to
  $D^b_{crs}(\g_X)$.
\end{enumerate}
\end{lemma}
\begin{pf}   $(2)\Rightarrow (1)$: Apply  \Theorem{inv-rs} and \Proposition{reg-inv}.
  
  $(1)\Rightarrow (2)$:  Pull back the distinguished triangle
  (\ref{loc-seq}) to a curve $\phi : C \to X$,
  \begin{displaymath}
    \phi^!(R\Gamma_{[S]} (M^\bullet)) \to \phi^!(M^\bullet) \to \phi^!(M^\bullet (*S))\to 
  \end{displaymath}
  Then either $\phi (C)\subset S$, or the support of the homology of $
  \phi^!(R\Gamma_{[S]} (M^\bullet))$ is discrete. In the latter case,
  $\phi^!(R\Gamma_{[S]} (M^\bullet)) \in D^b_{rs}(\g_C)$ by
  \Remark{discrete-case}, and since by assumption
  $\phi^!(M^\bullet)\in D^b_{rs}(\g_C)$, all vertices in the above
  distinguished triangle belong to $D^b_{rs}(\g_C)$. In the first, the
  assertion follows if $\phi^!(M^\bullet (*S))\cong
  0$, and this follows if $R\Gamma_{[C]}(M^\bullet (*S))\cong 0$
  \Th{kashiwara}.  The last assertion is well-known, but not
  finding an accurate reference we include the argument.
  
  First note that if $C\subset S$ is an inclusion of closed subsets of
  $X$, then the functor $R\Gamma_{[C]}(\cdot)$ is equivalent to
  $R\Gamma_{[C]}\circ R\Gamma_{[S]}(\cdot )$. To see this, the
  non-derived functors satisfy for any $\Oc_X$-module $M$ the
  composition property $H^0_{[C]}(H^0_{[S]}(M)) = H^0_{[C]}(M)$, where
  $H^0_{[S]}(M)= \{m\in M: I^n_S \cdot m=0, n\gg 1 \}$ ($I_S$ is the ideal of
  $S$); $H^0_{[C]}(\cdot)$ is defined similarly. Now the assertion
  follows since $H^0_{[S]}(\cdot)$ takes injective $\Oc_X$-modules to
  injective modules.
  
  Second, applying $R\Gamma_{[C]}$ to the distinguished triangle
\begin{displaymath}
  R\Gamma_{[S]}(M^\bullet) \to M^\bullet \to M^\bullet(*S)\to 
\end{displaymath}
one  gets $R\Gamma_{[C]}(M^\bullet(*S))\cong 0$.
\end{pf}

\begin{lemma}\label{split-div}
  Let $M^\bullet \in D^b(\g_X)$ and $x\in X$. There exists an open
  neighbourhood $X_0$ of $x$ and a divisor $S\subset X_0$
  such that  $M_{X_0}^\bullet(*S)\cong N^\bullet$ (in $D^b(\g_{X_0})$) where
  $N^\bullet$ is a complex of  free  $\Oc_{X_0}(*S)$-modules
  and $H^\bullet(M^\bullet  (*S))$ is a   free $\Oc_{X_0}(*S)$-module.
\end{lemma}
\begin{pf}
  This follows from \cite[VII, Lem. 9.3]{borel:Dmod}.
\end{pf}

\begin{pfof}{\Theorem{equivalent}}
  The inclusion $D^b_{rs}(\g_X) \subseteq D^b_{crs}(\g_X)$ follows
  from \Theorem{inv-rs}. Now letting $M^\bullet \in
  D^b_{crs}(\g_X)$ we have to prove that $M^\bullet \in
  D^b_{rs}(\g_X)$.
 
  The proof uses an induction over $d =\dim X$. The assertion is
  obvious when  $d=1$, so assume that $d\geq 2$ and that the theorem
  is true for all non-singular spaces $Y$ with $\dim Y < d$.
  
  a) The assertion is true when $\dim \supp M^\bullet < d$: The proof
  is by induction over the dimension $k <d$ of $Y_1 = \supp
  M^\bullet$.
  
  If $\dim Y_1 = 0$ the assertion follows from \Remark{discrete-case}.
 Assume that the assertion is true for
  all $M^\bullet$ such that $\dim \supp M^\bullet < k$.  Let
  $\pi: Y \to X$ be a desingularisation of $Y_1$ and consider the
  distinguished triangle in $D^b(\g_X)$ \Th{adj-prop}
\begin{displaymath}
\pi_+\pi^!(M^\bullet)[\dop_{Y,X}] \to M^{\bullet} \to    N^\bullet \to 
  \end{displaymath}
  \Propositions{composed}{preservation} imply 
  \begin{displaymath}
    \pi^!(M^\bullet)\in D^b_{crs}(\g_{Y}).
  \end{displaymath}
  As $\dim Y < d$, by induction $\pi^!(M^\bullet)\in
  D^b_{rs}(\g_{Y})$, so \Theorem{directreg} implies
  \begin{displaymath}
    \pi_+\pi^!(M^\bullet)\in D^b_{rs}(\g_X)\subseteq D^b_{crs}(\g_X).
  \end{displaymath}
  By assumption $M^\bullet \in D^b_{crs}(\g_X)$, so the third vertex
  $N^\bullet \in D^b_{crs}(\g_X)$; since $\dim \supp N^\bullet < k$
  \Th{kashiwara}, by induction we get $N^\bullet \in D^b_{rs}(\g_X)$. Then since
  $\pi_+\pi^!(M^\bullet)\in D^b_{rs}(\g_X)$ the third vertex  $M^\bullet\in
  D^b_{rs}(\g_X)$.
  
  b) The assertion holds when $\dim \supp M^\bullet =d$: The statement
  being local in $X$ one can, by \Lemma{split-div},  assume that
  $S\subset X$ is a divisor such that $M^\bullet(*S)$ and $H^\bullet
  (M^\bullet(*S))= H^\bullet (M^\bullet)(*S)$ are locally free over
  $\Oc_X(*S)$. 
  
  $M^\bullet(*S) \in D^b_{rs}(\g_X)$: We have to prove that $H^\bullet
  (M^\bullet(*S)) \in \coh_{rs}(\g_X)$, and
  for this, since $H^\bullet(M^\bullet(*S))$ is torsion free, one can
  assume that $S$ is non-singular \Cor{pure}.

As $\Oc_X(*S)$ is locally flat over $\Oc_X$ it follows
  that $M^\bullet(*S)$ is a complex of flat
  $\Oc_X$-modules. Therefore, if $\phi : C \to X$ is a curve,
  \begin{displaymath}
    H^\bullet(\phi^!(M^\bullet(*S))) \cong \phi^!(H^\bullet(M^\bullet(*S))).
  \end{displaymath}
  Since $M^\bullet(*S)\in D^b_{crs}(\g_X)$ it follows that
  $H^\bullet(M^\bullet(*S))\in D^b_{crs}$; $H^\bullet(M^\bullet(*S))$
  being torsion free \Corollary{curve-testc} implies that
  $H^\bullet(M^\bullet (*S))\in \coh_{rs}(\g_X)$.
  
  $M^\bullet \in D^b_{rs}(\g_X)$: Since $R\Gamma_{[S]}(M^\bullet)\in
  D^b_{crs}(\g_X)$ \Lem{div-curve}, and 
  \begin{displaymath}
    \dim \supp  R\Gamma_{[S]}(M^\bullet) \leq \dim S < d,
  \end{displaymath}
  by induction $R\Gamma_{[S]}(M^\bullet) \in D^b_{rs}(\g_X)$.
  Therefore the third vertex $M^\bullet \in D^b_{rs}(\g_X)$.
\end{pfof}

When $M^\bullet \in D^b_{rs}(\g_X)$ \Proposition{reg-inv} implies that
any sub-quotient of $H^\bullet(M^\bullet) $ belongs to
$\coh_{rs}(\g_X)$ (full regularity). If we had known that the category
$D^b_{crs}(\g_X)$ has the same property then \Lemma{split-div} would
be unnecessary in the above proof.
\begin{corollary}\label{full-reg}
The following are equivalent:
\begin{enumerate}
\item $M^\bullet \in D^b_{crs}(\g_X)$;
\item $H^\bullet (M^\bullet)\in D^b_{crs}(\g_X)$.
\end{enumerate}
If these conditions are satisfied, then any sub-quotient of
$H^\bullet(M^\bullet)$ belongs to $D^b_{crs}(\g_X)$.
\end{corollary}

\begin{remark}  
 When $\g_X=T_X$ we get \cite[VII, Cor.
  12.8]{borel:Dmod} and \cite[Th. 5.3.4]{bjork:analD}.
\end{remark}

\subsection{Completely regular complexes}
{\it In this section all spaces $X$ are quasi-projective complex algebraic
manifold, and $X_h$ is the associated complex analytic manifold.  All
Lie algebroids are locally free and transitive.}

\medskip Let $(j',j): (X,\g_X, \alpha) \to (\bar X, \g_{\bar X}, \bar
{\alpha})$ be a morphism of locally free transitive Lie algebroids
\Defn{lie-morphism}, where the underlying morphism $j: X\to \bar X$ is
an open embedding of algebraic manifolds. We say
that $(j',j)$  is a {\it completion}\/ of $(X, \g_X , \alpha)$ if  $\bar X$ is projective,
$\bar \alpha$ is surjective, and $j^+(\g_{\bar X}) = \g_X$.
Completions do not always exist.
\begin{example}
  Let $(\Cb^*,\g_{\Cb^*}, \alpha)$ be a transitive Lie algebroid on
  $\Cb^* = \Cb^1 \setminus \{0\}$ of the form $\g_{\Cb^*}=
  \bfr_{\Cb^*}\oplus T_{\Cb^*}$ where $\bfr_{\Cb^*}$ is a commutative
  $\Oc_{\Cb^*}$-Lie algebra, invertible as $\Oc_{\Cb^*}$-module, and
  moreover a $T_{\Cb^*}$-module with the global generator $\mu$
  satisfying the relation $(t^l\partial_t-\gamma)\mu =0$, where $l$ is a
  non-negative integer.  The structure of Lie algebroid on
  $\bfr_{\Cb^*}\oplus T_{\Cb^*}$ is $[(b_1,\partial_1),
  (b_2,\partial_2)] = (\partial_1(b_2)-\partial_2(b_1),[\partial_1,
  \partial_2])$, $(b_i,\partial_i)\in \bfr_{\Cb^*}\oplus T_{\Cb^*}$.
  The unique non-singular completion of $\Cb^*$ is the open embedding
  $j: \Cb^* \to \Pb^1$. Now  if $l> 1$, $\gamma \neq 0$, or $l=1$,
  $\gamma \notin \Zb$, then $\T_{\Cb^*}$-module
  $\Oc_{\Cb^*}\mu$ is not the restriction of  a $T_{\Pb^1}$-module
  which is coherent over $\Oc_{\Pb^1}$; hence  $(\Cb^*,\g_{\Cb^*}, \alpha)$ does not
  have a (transitive) completion.
\end{example}

\begin{remark} There always exist non-transitive completions of a Lie
  algebroid $(X,\g_X, \alpha)$.  Let $I$ be the ideal of $\bar X
  \setminus X$.  $\bar X$ is a noetherian scheme so there exists a coherent
  $\Oc_{\bar X}$-submodule $\g^1_{\bar X}\subseteq j_*(\g_X)$ such that
  $j^*(\g^1_{\bar X})= \g_X$.  For integers $n\gg 1$ consider the
  submodule $ I^n \g^1_{\bar X}\subset \g^1_{\bar X}$ and let $\bar
  \alpha$ denote the restriction of $j_*(\g_X)\to j_*(T_X)$ to
  $\g^1_{\bar X}$. When $n\gg 1$ one checks that 
  $[ I^n \g^1_{\bar X}, I^n \g^1_{\bar X}]\subset I^n \g^1_{\bar X}$
  and $\alpha(I^n \g^{1}_{\bar X}) \subset T_{\bar X}$. Thus for
  sufficiently big $n$ the $\Oc_{\bar X}$ -module $\g_{\bar X} =
  I^n\g^1_{\bar X}$ is a completion of $\g_X$ except that the morphism
  $\bar \alpha$ is not surjective.
\end{remark}

We will {\it assume} \/ that the Lie algebroids under consideration do
have a completion. For example, for the Lie algebroid $T_X$ any open
embedding $j: X\to \bar X$ to a proper algebraic manifold $\bar X$
will do, and more generally if $M$ is a locally free $\Oc_X$-module which has a
locally free extension to  $\bar X$, then the linear Lie
algebroid $\cf_X(M)$ has a completion.

Let $(j',j)$ be a completion.
\begin{definition}  
  An object $M^\bullet \in D^b_{rs}(\g_X)$ is {\it completely regular}
  \/ if $M^\bullet= j^!(N^\bullet)$ for some $N^\bullet \in
  D^b_{rs}(\g_{\bar X})$, where $(\bar X, \g_{\bar X}, \alpha)$ is a
  completion of $(X, \g_X, \alpha)$.
  \end{definition}
  Let $\bar D^b_{rs}(\g_X) \subset D^b_{rs}(\g_{X})$ be the
  sub-category of completely regular complexes. 
  
  We also say that $N^\bullet$ is a (regular) completion of
  $M^\bullet$.  By \Theorem{bernstein} there exist such $N^\bullet $
  if and only if $j_+(M^\bullet) = N^\bullet(*S) \in D^b_{rs}(\g_{\bar
    X})$, where $S= \bar X \setminus j(X)$, so it is natural to check
  if $j_+(M^\bullet)\in D^b_{rs}(\g_{\bar X})$.

Next  lemma shows that the definition of $\bar D^b_{rs}(\g_X)$
is intrinsic. 
\begin{lemma}\label{intrinsic}
  Let 
  \begin{displaymath}
    (j_1,j'_1):(X, \g_{X}, \alpha) \to (X_1, \g_{X_1}, \alpha_1)
  \end{displaymath}
 and
 \begin{displaymath}
(j_2, j'_2) : (X, \g_{X}, \alpha) \to (X_2, \g_{X_2},  \alpha_2)
\end{displaymath}
be completions. Put $S_2= X_2 \setminus j_2(X)$ and let $M^\bullet \in
D^b_{rs}(\g_X)$. Suppose that there exist $N^\bullet_i \in
D^b_{rs}(\g_{X_i})$, $i=1,2$, satisfying $j^!(N^\bullet_i)\cong
M^\bullet$. Then if $N^\bullet_1 \in D^b_{rs}(\g_{X_1})$ and
$R\Gamma_{[S_2]}(N_2^\bullet)\in D^b_{rs}(\g_{X_2})$ it follows that
$N^\bullet_2 \in D^b_{rs}(\g_{X_2})$.
\end{lemma}

\begin{pf} By \Theorem{bernstein} it
  suffices to prove that $N_2^\bullet(*S)\in D^b_{rs}(\g_{X_2})$.
  
  Let $X_3\to X_1\times X_2$ be a desingularisation of the closure of the image of the
  morphism $(j_1,j_2): X \to X_1\times X_2$, so there exist  canonical
  proper birational morphisms $p : X_3 \to X_1$, $q: X_3 \to X_2$, and
  letting $j: X \to X_3$ be the inclusion morphism of $X$ in $X_3$ we
  have  $q\circ j = p\circ j= \id_X$. Let $S_3, S_2$ be the algebraic  sets $X_3
  \setminus j(X)$ and $X_2\setminus j_2(X)$.

We have $N^\bullet_3 := p^!(N_1^\bullet)\in  D^b_{rs}(\g_{X_3})$
\Th{inv-rs} and $j^!(N^\bullet_3)=  j_2^!(N^\bullet_2)=M^\bullet$, so 
  \begin{eqnarray*}
    N^\bullet_2(*S_2) &=& (j_2)_+j_2^!(N^\bullet_2)\\
&=& (q\circ    j)_+j^!(N^\bullet_3)  \\
&=& q_+(j_+j^!(N^\bullet_3)) \quad \text{\Prop{direct-compose}} \\
&=& q_+(N^\bullet_3(*S_3)).
  \end{eqnarray*}
  By \Theorem{bernstein} $N_3(*S_3)\in D^b_{rs}(\g_{X_3})$, noting
  that the pull-back of a transitive Lie algebroid is transitive,
  hence by \Theorem{directreg} $N^\bullet_2(*S_2)\in
  D^b_{rs}(\g_{X_2})$.
\end{pf}

In the following two corollaries to \Section{reg-sing} $\pi : Y \to X$
is a morphism of quasi-projective algebraic manifolds. 

We first note  that $\g_Y$ has a completion if $\g_X$ has one.  By
assumption there exists a completion $(i',i):(X, \g_X, \alpha )\to
(\bar X, \g_{\bar X}, \bar \alpha)$. Let $Y\subseteq Y_1$ where $Y_1$
is projective, define $g:  Y \to X\times Y_1$, $g(y)= (\pi(y), y)$ and
let $\phi : \bar Y \to X\times Y_1$ be a desingularisation of the closure
of $g(Y)$. Now  there exists an open embedding $j: Y \to \bar Y$ such
that $\phi \circ j = g$ and a Cartesian diagram:
\begin{equation}\label{cart-complete}
  \begin{CD}
     Y @>j>> \bar Y \\ 
      @V{\pi}VV @VV{\bar \pi}V
    \\ X @>i>> \bar X.
  \end{CD}
\end{equation}
Since $\g_{\bar X}$ is locally free and transitive it follows that
${\bar \pi}^+(\g_{\bar X})$ is locally free and transitive
\Lem{pullfree}; hence
\begin{displaymath}
  (\pi^+(i'), j):(Y, \g_Y= \pi^+(\g_X), \pi^+(\alpha)) \to
(\bar Y , \pi^+(\g_{\bar X}), \pi^+(\bar \alpha))
\end{displaymath}
is a completion.

We have for  $M^\bullet \in D^b(\g_Y)$
 \begin{equation}\label{openemb}
    j_+\pi^! (M^\bullet) = \bar \pi^! i_+ (M^\bullet)
 \end{equation}
since $i$ and $j$ are open embeddings.
\begin{cor}\label{inv-complete}
 $\pi^!(\bar D^b_{rs}(\g_X))\subset \bar D^b_{rs}(\g_Y)$.
\end{cor}
\begin{pf}
As $i_+(M^\bullet)\in D^b_{rs}(\g_{\bar X})$ the result follows
from (\ref{openemb})  and \Theorem{inv-rs}.
\end{pf}

\begin{cor}\label{completedirect} Assume that $\pi$ can be factorized
  into  an open embedding which is a completion of $(Y, \g_Y,\alpha)$ 
  and a proper morphism. Then
  \begin{displaymath}
    \pi_+(\bar D^b_{rs}(\g_Y)) \subset \bar D^b_{rs}(\g_X).
  \end{displaymath}
\end{cor}

\begin{pf} Let $M^\bullet  \in \bar D^b_{rs}(\g_Y)$.
  If $\pi = p\circ j$ where $j : Y \to \bar Y$ is a completion we know
  that $j_+(M^\bullet)\in D^b_{rs}(Y)$, and since $p$ is proper
  $\pi_+(M^\bullet)= p_+ (j_+(M^\bullet))\in D^b_{rs}(\g_X)$
  \Th{directreg}. Since $p(\bar Y) $ is a closed  subspace of
  $\bar X$ the restriction of $i$ to $p(\bar Y)$ is proper. The proof
  of \Theorem{directreg} is now applicable to prove that
  $i_+\pi_+(M^\bullet)\in D^b_{rs}(\g_{\bar X})$ only using that the
  restriction of $i$ to the support of $p_+(j_+(M^\bullet))$ is
  proper.
\end{pf}
When $\g_X=T_X$ all open embeddings are completions so
\Corollary{completedirect} holds for any morphism $\pi$; in this case
it is proven for algebraic $\Dc_X$-modules by Borel \cite[VII, Th.
12.2]{borel:Dmod} using curve regularity as the definition of complete
regularity, cf. \Corollary{complete-equiv} below. The proof in [loc.\ 
cit.] is, as is ours, inspired by Deligne's proof that the Gauss-Manin
connection has regular singularities \cite[II. \S 7]{deligne:eq}. In
\cite[Ch. V, 5.5.28]{bjork:analD} the theorem is proven for analytic
$\Dc_X$-modules, by proving the comparison condition
    \begin{displaymath}
      R\Gamma_{[x]}(\Omega (T_X,\pi_+(M^\bullet))) \cong R\Gamma_{x}(\Omega (T_X ,\pi_+(M^\bullet))),
    \end{displaymath}
    where the left side is the tempered local cohomology of the de
    Rham complex of $\pi_+(M^\bullet)$, and the right side is the full
    local cohomology at any point $x\in X$, when the analogous
    comparison condition holds for $M$ at all points of $Y$.

Let $\bar D^b_{crs}(\g_X)\subset D^b_{crs}(\g_X)$ be the sub-category
of complexes $M^\bullet$ such that $\phi^!(M^\bullet)\in \bar
D^b_{rs}(\g_C)$ whenever $\phi : C \to X$ is a curve in $X$.

\begin{cor}\label{complete-equiv}
  \begin{displaymath}
\bar D^b_{crs} (\g_X) = \bar D^b_{rs}(\g_X).
\end{displaymath}
\end{cor}

\begin{pf}
  Let $(j,j'): (X, \g_X, \alpha ) \to (\bar X, \g_{\bar X}, \bar
  \alpha )$ be a completion, and $S$ be the algebraic set $\bar X
  \setminus j(X)$.
  
  $\bar D^b_{rs}(\g_X)\subseteq \bar D^b_{crs}(\g_X)$: A morphism $\phi
  : C\to X$, where $C$ is curve, can be completed to a morphism $\bar
  \phi :\bar C \to \bar X$, where $\bar C$ is a proper curve, by the
  valuative criterion of properness.  Let $i: C\to \bar C$ be the
  corresponding open inclusion of (non-singular) curves, so $\bar \phi \circ i =
  j\circ \phi$.  Then if $N^\bullet \in D^b_{rs}(\g_{\bar X})$, we
  have $\bar \phi^! (N^\bullet )\in D^b_{rs}(\g_{\bar C})$
  \Cor{inv-complete} and $i^!\bar \phi^!(N^\bullet) = \phi^!j^!(N^\bullet)$
  \Prop{composed}; hence $\phi^!\circ j^!(N^\bullet) \in \bar
  D^b_{cr}(\g_C)$; hence $j^!(N^\bullet)\in \bar D^b_{crs}(\g_X)$.

  $\bar D^b_{crs}(\g_X)\subseteq \bar D^b_{rs}(\g_X)$: Let $M^\bullet
  \in \bar D^b_{crs}(\g_X)$ and $N^\bullet \in D^b_{\coh}(\g_{\bar
    X})$ be such that $j^!(N^\bullet)= M^\bullet$ and $N^\bullet =
  N^\bullet(*S)$.  By \Theorem{equivalent} we have to prove that
  $N^\bullet \in D^b_{crs}(\g_{\bar X})$. Let $\bar \phi : \bar C \to
  \bar X$ be a curve.  Then either $\bar \phi (\bar C)$ is contained
  in $S$, or $C:=\bar \phi^{-1}(j(X)) \subset \bar C$ is Zariski
  dense. In the first case, $\bar \phi^!(N^\bullet(*S))=0$, see the
  proof of \Lemma{div-curve} or \cite[VI, Cor.  8.5]{borel:Dmod};
  hence we can assume that $i: C \to \bar C$ is an inclusion of a
  Zariski open dense subset of $C$. Now, by assumption, $i^!\bar
  \phi^!(N^\bullet) = \phi^!  j^!(N^\bullet) \in \bar D^b_{rs}(\g_C)$;
  hence $\bar \phi^!(N^\bullet) \in D^b_{rs}(\g_{\bar C})$.
\end{pf}
One can improve  \Proposition{inv-equiv} when considering  $\bar
D^b_{rs}(\g_X)$:
\begin{prop}\label{inv-equiv-complete}  Let $M^\bullet \in D^b_{\coh}(\g_X)$.  Assume that
  $\pi$ is surjective on $\supp M^\bullet \subset X$. Then the
  following are equivalent:
\begin{enumerate} 
\item $M^\bullet \in  \bar D^b_{rs}(\g_X)$.
\item $\pi^!(M^\bullet )\in \bar D^b_{rs}(\g_Y)$.
  \end{enumerate}
\end{prop}
When $\g_X=T_X$ we recover \cite[Th. 5.6.5]{bjork:analD} and
\cite[VII, Prop. 12.9]{borel:Dmod}. \Proposition{inv-equiv-complete}
can be used to prove the complete regularity of certain equivariant
$\g_X$-modules occurring as localizations of Harish-Chandra
representations of Lie algebras.
  \begin{pf} $(1)\Rightarrow (2)$ follows from \Corollary{inv-complete},
    so we need to prove $(2)\Rightarrow (1)$.
    
    Letting $\alpha : C\to X$ be a curve one has to prove that
    $\pi^!(M^\bullet)\in \bar D^b_{rs}(\g_Y)$ implies
    $\alpha^!(M^\bullet)\in \bar D^b_{rs}(\g_C)$ \Cor{complete-equiv}.
    Let $p:Y_1 \to C$ be a desingularisation of the base change
    $\alpha$ of $\pi$, and $\phi : Y_1 \to Y$ be the second
    projection.  Since $p^! \circ \alpha^! \cong \phi^!  \circ \pi^!$
    \Prop{composed} we have that $p^!  \alpha^!  (M^\bullet)\in \bar
    D^b_{rs}$.  It therefore suffices to prove the proposition when
    $\pi$ is a morphism to a curve $X$.  Moreover, there exists a
    curve $\alpha_1 :C_1 \to Y$ such that the composed morphism
    $\pi\circ \alpha_1 : C_1 \to X$ is dominant, which moreover can be
    assumed to be etale after restricting to an open subset of $C_1$.
    Since $\alpha_1^!$ preserves $\bar D^b_{rs}(\g_Y)$
    \Cor{inv-complete}, by \Lemma{intrinsic} we can then assume that
    $\pi : Y\to X$ is an etale morphism of non-singular curves.  Then
    as complex of $\Oc_Y$-modules the inverse image $\pi^!(M^\bullet)$
    coincides with $\pi^*(M^\bullet)$, the inverse image in the
    category of $\Oc_Y$-modules, and if $N\in \bar D^b_{rs}(\g_Y)$,
    then $\pi_+(N^\bullet)$ coincides with $\pi_*(N^\bullet)$, the
    direct image in the category of $\Oc_X$-modules. Moreover, that
    $\pi$ is etale also implies that the canonical morphism $M^\bullet
    \to \pi_+\pi^!(M^\bullet) = \pi_*\pi^*(M^\bullet)$ is a
    monomorphism in a triangulated category, hence it is split.
    Consider a completions of $X$ and $Y$ forming a Cartesian diagram
    of the type (\ref{cart-complete}).  Then $\bar \pi_+ j_+
    \pi^!(M^\bullet) \cong i_+\pi_+\pi^!(M^\bullet) \cong
    i_*\pi_*\pi^*(M^\bullet)$. Hence we have a split distinguished
    triangle \begin{displaymath} i_+(M^\bullet) \to \bar \pi_+j_+
      \pi^!(M^\bullet) \to C^\bullet \xrightarrow{\ 0\ },
    \end{displaymath}
so $H^\bullet (i_+(M^\bullet)) \subset H^\bullet(\bar \pi_+j_+
\pi^!(M^\bullet))$; by assumption   and \Theorem{directreg}  $H^\bullet(\bar
    \pi_+j_+\pi^!(M^\bullet))\in \coh_{rs}(\g_{\bar X})$; hence
    $H^\bullet(i_+(M^\bullet))\in \coh_{rs}(\g_{\bar X})$ \Prop{reg-inv}.
  \end{pf}
  
  Let us make a comment about the role of $\bar D^b_{rs}(\g_X)$ in the
  Riemann-Hilbert correspondence; compare also to the rather different
  argument in \cite[VIII]{borel:Dmod}.  Fix a completion $(X, \g_X,
  \alpha)\to (\bar X, \g_{\bar X},\bar \alpha)$. Put $S= \bar X
  \setminus j(X)$ and let $D^b_{rs}(S, \g_{\bar X})\subset
  D^b_{rs}(\g_{\bar X})$ be the full sub-category of objects
  $N^\bullet\in D^b_{rs}(\g_{\bar X})$ such that $N^\bullet =
  N^\bullet(*S)= Rj_*j^*(N^\bullet)$.  Then the functor
  \begin{displaymath}
    j_+ : \bar D^b_{rs}(\g_X) \to D^b_{rs}(S,\g_{\bar X}), \quad
    M^\bullet \mapsto j_+(M^\bullet)
  \end{displaymath}
is an equivalence of categories, and using GAGA (\ref{gaga-section}) 
\begin{displaymath}
g: D^b_{rs}(S, \g_{\bar X})\cong D^b_{rs}(S, \g_{\bar X_h}).
\end{displaymath}
A main step is to prove that the restriction functor defines an
equivalence
\begin{displaymath}
  r:  D^b_{rs}(S, \g_{\bar X_h}) \cong  D^b_{rs}(\g_{X_h}) .
\end{displaymath}
Since this is a local problem in $X_h$ the well-known proof when
$\g_{X_h}= T_{X_h}$ works, based on Hironaka's resolution of
singularities; see Malgrange's account in \cite[IV]{borel:Dmod} of
Deligne's proof for the case of connections.  Assume now that $\g_X\in
\LPic_{X}$ and recall the discussion in (\ref{twisted}).  Put $ \phi =
t([\g_{X_h}]) \in H^2(X_h, \Cb^*)$, so one has the invertible
$\phi$-twisted sheaf $\lambda(\phi)\in \Mod(\g_{X_h}, \phi)$, and the
solution functor
\begin{displaymath}
  Sol_\phi: D^b(\g_{X_h})\to D^b(\Cb_{X_h}, \phi), \quad M \to
  RHom_{\g_{X_h}}(M, \lambda(\phi)).
\end{displaymath}
The usual Riemann-Hilbert correspondence states that $Sol_0$ is
an equivalence of categories between $D^b_{rs}(T_{X_h})$ and the
derived category $D^b({\Cb_{X_h}})$ of complexes of sheaves whose
homology sheaves are constructible, while the abelian sub-category
$\coh_{rs}(T_{X_h})$ is equivalent to the category of perverse sheaves
on $X_h$ (\cite{kashiwara:constr}, \cite{kashiwara:riemann},
\cite{mebkhout:riemann}, \cite{mebkhout:bidual}, \cite{mebkhout:six}).
Now  $Sol_\phi$ defines  an equivalence  between
$D^b_{rs}(\g_{X_h})$ and the derived category of complexes of
$\phi$-twisted sheaves whose homology $\phi$-twisted sheaves are
constructible; see also \cite{kashiwara:rep}. Hence by the previous
equivalences there exists an equivalence between an algebraic and
a topological category:
\begin{displaymath}
 Sol_\phi \circ r \circ g\circ j_+ : \bar D^b_{rs}(\g_X) \cong  D^b(\Cb_{X_h}, \phi).
\end{displaymath}

\bibliographystyle{alpha} 
\bibliography{mindatabas}

\begin{thebibliography}{BGK{\etalchar{+}}87}

\bibitem[Ati57]{atiyah:connections}
Michael Atiyah.
\newblock Complex analytic connections in fibre bundles.
\newblock {\em Trans. Amer. Math. Soc.}, 85:181--207, 1957.

\bibitem[BB81]{beilinson-be}
Alexander Beilinson and Joseph Bernstein.
\newblock Localisation de {$\g$}-modules.
\newblock {\em C.R. Acad. Sci. (Paris)}, 292(Ser. I.), 1981.

\bibitem[BB93]{bei-ber:jantzen}
Alexander Beilinson and Joseph Bernstein.
\newblock A proof of {J}antzen conjectures.
\newblock {\em Adv. Sov. Math.}, 16(1):1--50, 1993.

\bibitem[BBD82]{be-be-de}
Alexander Beilinson, Joseph Bernstein, and Pierre Deligne.
\newblock Fasciuex pervers.
\newblock {\em Asterisque)}, 100:1--100, 1982.

\bibitem[BGG72]{bernstein-gege:cubic}
I.~N. Bern{\v{s}}te{\u\i}n, I.~M. Gelfand, and S.~I. Gelfand.
\newblock Differential operators on a cubic cone.
\newblock {\em Uspehi Mat. Nauk}, 27(1(163)):185--190, 1972.

\bibitem[BGK{\etalchar{+}}87]{borel:Dmod}
A.~Borel, P.-P. Grivel, B.~Kaup, A.~Haefliger, B.~Malgrange, and F.~Ehlers.
\newblock {\em Algebraic ${D}$-modules}, volume~2 of {\em Perspectives in
  Mathematics}.
\newblock Academic Press Inc., Boston, MA, 1987.

\bibitem[Bj{\"o}93]{bjork:analD}
Jan-Erik Bj{\"o}rk.
\newblock {\em Analytic $D$-modules and applications.}
\newblock Kluwer Academic Publishers, 1993.

\bibitem[BK81]{brylinski-ka}
Jean-Luc Brylinski and Masaki Kashiwara.
\newblock Kazhdan-{L}usztig conjecture and holonomic systems.
\newblock {\em Invent. Math.}, 64(387), 1981.

\bibitem[BS58]{borel-serre:riemann-roch}
Armand Borel and Jean-Pierre Serre.
\newblock Le th\'eor\`eme de {R}iemann-{R}och (d'apr\`es des r\'esultats
  in\'edits de {A}. {G}rothendieck).
\newblock {\em Bull. Soc. Math France}, 86:97--136, 1958.

\bibitem[BS76]{banica}
Constantin B{\u{a}}nic{\u{a}} and Octavian St{\u{a}}n{\u{a}}{\c{s}}il{\u{a}}.
\newblock {\em Algebraic methods in the global theory of complex spaces}.
\newblock Editura Academiei, Bucharest, 1976.
\newblock Translated from the Romanian.

\bibitem[BS88]{beilinson-schechtman}
Alexander Beilinson and Vadim Schechtman.
\newblock Determinant bundles and virasoro algebras.
\newblock {\em Commun. Math. Phys.}, 118:651 -- 701, 1988.

\bibitem[BS98]{brodmann-bruns}
M.~P. Brodmann and R.~Y. Sharp.
\newblock {\em Local cohomology: an algebraic introduction with geometric
  applications}, volume~60 of {\em Cambridge Studies in Advanced Mathematics}.
\newblock Cambridge University Press, Cambridge, 1998.

\bibitem[Del70]{deligne:eq}
Pierre Deligne.
\newblock {\em \'Equations diff\' erentielles \`a points singuliers r\'
  eguliers}, volume 163.
\newblock Springer Math. L.N, 1970.

\bibitem[Fri67]{frisch}
Jacques Frisch.
\newblock Points de platitude d'un morphisme d'espaces analytiques complexes.
\newblock {\em Invent. Math.}, 4:118--138, 1967.

\bibitem[Gab81]{gabber:integrability}
Ofer Gabber.
\newblock The integrability of the characteristic variety.
\newblock {\em Am Jour Math}, 103(3):445--468, 1981.

\bibitem[GD61]{EGA3}
Alexander Grothendieck and Jean Dieudonne.
\newblock \'etude cohomolique des faisceaux coh\'erent.
\newblock {\em Publ. IHES}, (4), 1961.

\bibitem[GD71]{EGA1}
Alexander Grothendieck and Jean Dieudonne.
\newblock {\em El\'ements de g\'eom\'etrie algebrique. I}.
\newblock Berlin-Heidelberg-New York: Springer-Verlag, 1971.

\bibitem[Gir71]{giraud}
Jean Giraud.
\newblock {\em Cohomologie non ab\'elienne}.
\newblock Springer-Verlag, Berlin, 1971.
\newblock Die Grundlehren der mathematischen Wissenschaften, Band 179.

\bibitem[God73]{godement:faisceaux}
Roger Godement.
\newblock {\em Topologie alg\'ebrique et th\'eorie des faisceaux}.
\newblock Hermann, Paris, 1973.
\newblock Troisi\`eme \'edition revue et corrig\'ee, Publications de l'Institut
  de Math\'ematique de l'Universit\'e de Strasbourg, XIII, Actualit\'es
  Scientifiques et Industrielles, No. 1252.

\bibitem[GR79]{grauert-remmert:stein}
Hans Grauert and Reinhold Remmert.
\newblock {\em Theory of Stein spaces}, volume 236 of {\em Grundlehren der
  Mathematischen Wissenschaften}.
\newblock Springer-Verlag, 1979.

\bibitem[Gra60]{grauert}
Hans Grauert.
\newblock Ein theorem der analytischen garbentheorie und die modulr{\"a}ume
  komplexer strukturen.
\newblock {\em Publ. IHES}, (5), 1960.

\bibitem[Gro62]{grothendieck:cohomologie}
Alexander Grothendieck.
\newblock {\em Cohomologie locale de faisceaux coh\'erents et Th\'eor\`emes de
  Lefshetz locaux et globaux (SGA 2)}, volume~2.
\newblock North-Holland, advanced studies in pure mathematics edition, 1962.

\bibitem[Gro67]{EGA4}
Alexander Grothendieck.
\newblock \'{E}l\'ements de g\'eom\'etrie alg\'ebrique. {I}{V}. \'{E}tude
  locale des sch\'emas et des morphismes de sch\'emas {I}{V}.
\newblock {\em Inst. Hautes \'Etudes Sci. Publ. Math. No.}, 32:361, 1967.

\bibitem[Gro71]{SGA1}
Alexander Grothendieck{, et al.}
\newblock {\em Rev\^etements \'etales et groupe fondamental}.
\newblock Springer-Verlag, Berlin, 1971.
\newblock S\'eminaire de G\'eom\'etrie Alg\'ebrique du Bois Marie 1960--1961
  (SGA 1), Dirig\'e par Alexandre Grothendieck. Augment\'e de deux expos\'es de
  M. Raynaud, Lecture Notes in Mathematics, Vol. 224.

\bibitem[Gun90]{gunning:2}
Robert~C. Gunning.
\newblock {\em Introduction to holomorphic functions of several variables.
  {V}ol. {I}{I}}.
\newblock The Wadsworth \& Brooks/Cole Mathematics Series. Wadsworth \&
  Brooks/Cole Advanced Books \& Software, Monterey, CA, 1990.

\bibitem[Har61]{harper}
R.~Jr. Harper, Laurence.
\newblock On differentiably simple algebras.
\newblock {\em Trans. Amer. Math. Soc.}, 100:63--72, 1961.

\bibitem[Har66]{hartshorne-res}
Robin Hartshorne.
\newblock {\em Residues and duality}.
\newblock Springer-Verlag, Berlin, 1966.
\newblock Lecture notes of a seminar on the work of A. Grothendieck, given at
  Harvard 1963/64. With an appendix by P. Deligne. Lecture Notes in
  Mathematics, No. 20.

\bibitem[Jan79]{jantzen:mod}
Jens~Carsten Jantzen.
\newblock {\em Moduln mit einem h{\"o}chsten {G}ewicht}, volume 750 of {\em
  Lect. Notes Math.}
\newblock Springer-Verlag, 1979.

\bibitem[K{\"a}l96]{kallstrom:gldim}
Rolf K{\"a}llstr{\"o}m.
\newblock Homological dimensions of {L}ie algebroids.
\newblock Manuscript, 1996.

\bibitem[K{\"a}l98]{kallstrom:2smooth}
Rolf K{\"a}llstr{\"o}m.
\newblock Smooth modules over lie algebroids, {II}.
\newblock {\em Manuscript}, 1998.

\bibitem[Kas]{kashiwara:master}
Masaki Kashiwara.
\newblock Algebraic study of systems of partial differential equations.
\newblock {\em M\'em. Soc. Math. France (N.S.)}, (63).

\bibitem[Kas80]{kashiwara:constr}
Masaki Kashiwara.
\newblock Faisceaux constructibles et syst\`emes holon\^omes d'\'equations aux
  d\'eriv\'ees partielles lin\'eaires \`a points singuliers r\'eguliers.
\newblock In {\em S\'eminaire Goulaouic-Schwartz, 1979--1980 (French)}, pages
  Exp. No. 19, 7. \'Ecole Polytech., Palaiseau, 1980.

\bibitem[Kas84]{kashiwara:riemann}
Masaki Kashiwara.
\newblock The {R}iemann- {H}ilbert problem for holonomic systems.
\newblock {\em Publ. Res. Inst. Math. Sci.}, 20(2):319--365, 1984.

\bibitem[Kas89]{kashiwara:rep}
Masaki Kashiwara.
\newblock Representation theory and {D}-modules on flag varieties.
\newblock {\em Soci\'et\' e Math. de France, Ast\' erisque}, pages 55--190,
  1989.

\bibitem[Kat87]{katz:galois}
Nicholas~M. Katz.
\newblock On the calculation of some differential {G}alois groups.
\newblock {\em Invent. Math.}, 87(1):13--61, 1987.

\bibitem[KK81]{kashiwara-ka}
Masaki Kashiwara and Takahiro Kawai.
\newblock On holonomic systems of microdifferential equations. {I}{I}{I}.
  {S}ystems with regular singularities.
\newblock {\em Publ. Res. Inst. Math. Sci.}, 17(3):813--979, 1981.

\bibitem[Lev75]{levelt}
A.~H.~M. Levelt.
\newblock Jordan decomposition for a class of singular differential operators.
\newblock {\em Ark. Mat.}, 13:1--27, 1975.

\bibitem[Lyu93]{lyubeznik}
Gennady Lyubeznik.
\newblock Finiteness properties of local cohomology modules (an application of
  ${D}$-modules to commutative algebra).
\newblock {\em Invent. Math.}, 113(1):41--55, 1993.

\bibitem[Mac95]{mackenzie-kirill}
Kirill C.~H. Mackenzie.
\newblock Lie algebroids and {L}ie pseudoalgebras.
\newblock {\em Bull. London Math. Soc.}, 27(2):97--147, 1995.

\bibitem[Mal96]{malgrange:reseaux}
Bernard Malgrange.
\newblock Connexions m\'eromorphes 2.
\newblock {\em Invent. Math}, 124:367--387, 1996.

\bibitem[Mat86]{matsumura}
Hideyuki Matsumura.
\newblock {\em Commutative Ring Theory}.
\newblock Cambridge University Press, 1986.

\bibitem[Meb80]{mebkhout:riemann}
Zoghman Mebkhout.
\newblock Sur le probl\`eme de {H}ilbert-{R}iemann.
\newblock In {\em Complex analysis, microlocal calculus and relativistic
  quantum theory (Proc. Internat. Colloq., Centre Phys., Les Houches, 1979)},
  volume 126 of {\em Lecture Notes in Phys.}, pages 90--110. Springer, Berlin,
  1980.

\bibitem[Meb82]{mebkhout:bidual}
Zoghman Mebkhout.
\newblock Th\'eor\`emes de bidualit\'e locale pour les ${\Dc}_{X}$-modules
  holonomes.
\newblock {\em Ark. Mat.}, 20(1):111--124, 1982.

\bibitem[Meb84]{mebkhout:equivalence}
Zoghman Mebkhout.
\newblock Une \'equivalence de cat\'egories, and une autre \'equivalence de
  cat\'egories.
\newblock {\em Compositio Math.}, 51(1):51--62, 63--88, 1984.

\bibitem[Meb89]{mebkhout:six}
Zoghman Mebkhout.
\newblock {\em Le formalisme des six op\'erations de {G}rothendieck pour les
  $\Dc_X$-modules coh\'erents}, volume~35 of {\em Travaux en Cours [Works in
  Progress]}.
\newblock Hermann, Paris, 1989.
\newblock With supplementary material by the author and L. Narv\'aez Macarro.

\bibitem[Oda83]{oda}
Tadao Oda.
\newblock Introduction to algebraic analysis on complex manifolds.
\newblock In {\em Algebraic varieties and analytic varieties (Tokyo, 1981)},
  volume~1 of {\em Adv. Stud. Pure Math.}, pages 29--48. North-Holland,
  Amsterdam, 1983.

\bibitem[Rei67]{reiffen}
Hans-J{\"o}rg Reiffen.
\newblock Das {L}emma von {P}oincar\'e f\"ur holomorphe {D}ifferential-formen
  auf komplexen {R}\"aumen.
\newblock {\em Math. Z.}, 101:269--284, 1967.

\bibitem[Rin63]{rinehart:difforms}
George~S. Rinehart.
\newblock Differential forms on general commutative algebras.
\newblock {\em Trans. Amer. Math. Soc.}, 108:195--222, 1963.

\bibitem[Sai80]{saito-kyoji:log}
Kyoji Saito.
\newblock Theory of logarithmic differential forms and logarithmic vector
  fields.
\newblock {\em J. Fac. Sci. Univ. Tokyo Sect. IA Math.}, 27(2):265--291, 1980.

\bibitem[Sch64]{scheja:fortsetzung}
G{\"u}nter Scheja.
\newblock Fortsetzungss\"atze der komplex-analytischen {C}ohomologie und ihre
  algebraische {C}harakterisierung.
\newblock {\em Math. Ann.}, 157:75--94, 1964.

\bibitem[Sei66]{seidenberg}
Abraham Seidenberg.
\newblock Derivations and integral closure.
\newblock {\em Pacific J. of Math.}, 16:167--73, 1966.

\bibitem[Ser56]{serre:gaga}
Jean-Pierre Serre.
\newblock G\'eom\'etrie alg\'ebrique et g\'eom\'etrie analytique.
\newblock {\em Ann. Inst. Fourier, Grenoble}, 6:1--42, 1955--1956.

\bibitem[Ser66]{serre:prolongement}
Jean-Pierre Serre.
\newblock Prolongement de faisceaux analytiques coherents.
\newblock {\em Ann. Inst. Fourier}, 16(1):363--374, 1966.

\bibitem[Siu70]{siu:local}
Yum-{T}ong Siu.
\newblock Analytic sheaves of local cohomology.
\newblock {\em Trans. Amer. Math. Soc.}, 148:347--366, 1970.

\bibitem[Tra69]{trautmann}
G{\"u}nther Trautmann.
\newblock Ein {E}ndlichkeitssatz in der analytischen {G}eometrie.
\newblock {\em Invent. Math.}, 8:143--174, 1969.

\bibitem[TUY89]{Tsuchiya-ueya}
Akihiro Tsuchiya, Kenji Ueno, and Yasuhiko Yamada.
\newblock Conformal field theory on universal family of stable curves with
  gauge symmetries.
\newblock In {\em Integrable systems in quantum field theory and statistical
  mechanics}, volume~19 of {\em Adv. Stud. Pure Math.}, pages 459--566.
  Academic Press, Boston, MA, 1989.

\bibitem[Uen97]{ueno:conformal}
Kenji Ueno.
\newblock Introduction to conformal field theory with gauge symmetries.
\newblock In {\em Geometry and physics (Aarhus, 1995)}, volume 184 of {\em
  Lecture Notes in Pure and Appl. Math.}, pages 603--745. Dekker, New York,
  1997.

\bibitem[Wei58]{weil:kahler}
Andr{\'e} Weil.
\newblock {\em Introduction \`a l'\'etude des vari\'et\'es k\"ahl\'eriennes}.
\newblock Hermann, Paris, 1958.
\newblock Publications de l'Institut de Math\'ematique de l'Universit\'e de
  Nancago, VI. Actualit\'es Sci. Ind. no. 1267.

\end{thebibliography}
\newcommand{\etalchar}[1]{$^{#1}$}

\end{document}